\documentclass[12pt]{svjour3}

\makeatletter
\def\ps@pprintTitle{%
 \let\@oddhead\@empty
 \let\@evenhead\@empty
 \def\@oddfoot{}%
 \let\@evenfoot\@oddfoot}
\makeatother

\usepackage{placeins,subcaption}
\usepackage{amsmath,amscd,wrapfig,amsfonts,amssymb,mathtools,mathrsfs}
\usepackage{booktabs,multirow,lscape,cool,bm}
\usepackage{url,parskip,caption}

\usepackage[top=1.0in, bottom=1.0in, left=1.0in, right=1.0in]{geometry}

\makeatletter
\g@addto@macro\normalsize{%
  \setlength\abovedisplayskip{.6em}
  \setlength\belowdisplayskip{.6em}
  \setlength\abovedisplayshortskip{.6em}
  \setlength\belowdisplayshortskip{.6em}
}

\begin{document}


\title{The Levin approach to the numerical calculation of phase functions}
\date{}
\author{Murdock Aubry         \and
        James Bremer
}

\institute{Murdock Aubry \at
              University of Toronto \\
              \email{murdock.aubry@mail.utoronto.ca}           
           \and
              James Bremer \at
              Department of Mathematics \\
              University of Toronto \\
              \email{bremer@math.toronto.edu}           
}

\maketitle
\vskip -7em

\begin{abstract}
When the eigenvalues of the coefficient matrix for a linear scalar ordinary differential equation
are of large magnitude, its solutions exhibit complicated behaviour, such as high-frequency oscillations,
rapid growth or rapid decay.  The cost of representing such solutions  using standard techniques typically 
grows with the magnitudes of the eigenvalues.
As a consequence, the  running times of standard solvers for ordinary differential equations also grow with these eigenvalues.
It is well known, however, that a large class of  scalar ordinary differential equations with slowly-varying coefficients
admit slowly-varying phase functions that can be  represented efficiently,
regardless of the magnitudes of the eigenvalues of the corresponding coefficient matrix.
  Here, we introduce two numerical algorithms 
for constructing slowly-varying phase functions for linear scalar ordinary differential equations
 inspired by the classical Levin method for evaluating oscillatory integrals.
The running times of our algorithms depend on the complexity of an equation's coefficients,
but are largely independent of the magnitudes of the  eigenvalues of its coefficient matrix.
Once these phase functions have been constructed, essentially
any reasonable initial or boundary value problem for the scalar equation can be easily solved,
essentially instantaneously.
The results of extensive numerical experiments  demonstrating the properties of our algorithms
are presented.

\end{abstract}










\begin{section}{Introduction}
The cost to represent the solutions of an $n^{th}$ order linear homogeneous ordinary differential equation
\begin{equation}
y^{(n)}(t) + q_{n-1}(t) y^{(n-1)}(t) + \cdots + q_1(t) y'(t) + q_0(t) y(t) = 0
\label{introduction:scalarode}
\end{equation}
using standard techniques, such as polynomial or trigonometric expansions,
increases with the magnitudes of the eigenvalues $\lambda_1(t), \ldots, \lambda_n(t)$
of the coefficient matrix
\begin{equation}
A(t) = \left(\begin{array}{cccccccc}
0 & 1 & 0 & \cdots & 0 & 0 \\
0 & 0 & 1 & \cdots & 0 & 0 \\
\vdots &  &  & \ddots &  & \vdots\\
0 & 0 & 0 & \cdots & 1 & 0 \\
0 & 0 & 0 & \cdots & 0 & 1 \\
-q_0(t) & -q_1(t) & -q_2(t) & \cdots&  -q_{n-2}(t) & -q_{n-1}(t)\\
\end{array}
\right)
\label{introduction:scalarcoef}
\end{equation}
obtained from  (\ref{introduction:scalarode}) in the usual way.  
As a consequence, when conventional solvers 
for ordinary differential equations  are applied to equations of this form, their  running times  
also increase with the magnitudes of the eigenvalues of (\ref{introduction:scalarcoef}).

It is well known, however, that many such equations admit phase functions whose
cost to represent via standard methods depends on the complexity of the coefficients $q_0,\ldots,q_{n-1}$,
but not the magnitudes of the eigenvalues of (\ref{introduction:scalarcoef}).
Indeed, if the $q_0,\ldots,q_{n-1}$ are slowly-varying on an interval $I$ and 
the coefficient matrix (\ref{introduction:scalarcoef})
corresponding to (\ref{introduction:scalarode}) 
has eigenvalues  $\lambda_1(t), \ldots, \lambda_n(t)$  which are
 distinct for all  $t$ in $I$, then it is possible to find slowly-varying functions 
$\psi_1,\ldots,\psi_n\colon I \to \mathbb{C}$
such that 
\begin{equation}
\left\{ \exp\left(\psi_j(t)\right)  :  j=1,\ldots,n \right\}
\label{introduction:phaserep}
\end{equation}
is a basis for the space of solutions of (\ref{introduction:scalarode}) given on the interval $I$.
This observation is the foundation of the WKB method, as well as many other techniques for the asymptotic analysis
 of ordinary differential equations (see, for instance, \cite{Miller}, \cite{Wasov}
and \cite{SpiglerPhase1,SpiglerPhase2,SpiglerZeros}).  

An equation of the form (\ref{introduction:scalarode}) is said to be nondegenerate on the interval $I$
provided the above condition holds; that is, if all of the eigenvalues $\lambda_1(t), \ldots, \lambda_n(t)$  of the corresponding
coefficient matrix (\ref{introduction:scalarcoef}) are distinct for each $t \in I$.
Moreover, $t_0$ is a turning point for (\ref{introduction:scalarode}) if
the eigenvalues 
of (\ref{introduction:scalarcoef})
are distinct in a deleted neighbourhood of $t_0$ but coalesce at $t_0$. 
The derivatives of the phase functions $\psi_1,\ldots,\psi_n$, which we will 
denote by $r_1,\ldots,r_n$, satisfy  an $(n-1)^{st}$ order  nonlinear inhomogeneous ordinary differential equation,
the general form of which is quite complicated.   When $n=2$, it is the Riccati equation
\begin{equation}
r'(t) + (r(t))^2 + q_1(t) r(t) + q_0(t) = 0;
\label{introduction:riccati1}
\end{equation}
when $n=3$, the nonlinear equation is
\begin{equation}
r''(t) + 3 r'(t) r(t) + (r(t))^3 + q_2(t) r'(t) + q_2(t) (r(t))^2 + q_1(t) r(t) + q_0(t) = 0;
\label{riccati2}
\end{equation}
and, for $n=4$, we have
\begin{equation}
\begin{aligned}
r'''(t) &+ 4 r''(t) r(t) + 3 (r'(t))^2 + 6  r'(t) (r(t))^2 + (r(t))^4 +q_3(t) (r(t))^3 + q_3(t) r''(t) + 3 q_3(t) r'(t) r(t) \\
&+q_2(t) (r(t))^2 +q_2(t) r'(t) + q_1(t) r(t) +q_0(t) = 0.
\end{aligned}
\label{riccati3}
\end{equation}
By a slight abuse of terminology, we will refer to the $(n-1)^{st}$ order nonlinear
equation obtained by inserting the representation 
\begin{equation}
y(t) = \exp\left(\int r(t)\, dt\right)
\end{equation}
into (\ref{introduction:scalarode}) as the $(n-1)^{st}$ order Riccati equation, or, alternatively, 
the Riccati equation for (\ref{introduction:scalarode}).


An obvious approach to initial and boundary value boundary problems for (\ref{introduction:scalarode}) calls for constructing
a suitable collection of slowly-varying phase functions by solving the corresponding Riccati equation numerically.
Doing so is not as straightforward as it sounds, however.  The principal difficulty is that 
most solutions of the Riccati equation for (\ref{introduction:scalarode}) are rapidly-varying when 
the eigenvalues $\lambda_1(t),\ldots,\lambda_n(t)$ are of large magnitude, and some mechanism is needed to 
select the slowly-varying solutions.

The article \cite{BremerPhase} introduces such a technique for the case of 
second order linear ordinary differential equations, all of which can be put in the form
\begin{equation}
y''(t) + q(t) y(t) = 0.
\label{introduction:two}
\end{equation}
 The eigenvalues of the corresponding coefficient matrix (\ref{introduction:scalarcoef}) are
given by
\begin{equation}
\lambda_1(t) =  \sqrt{-q(t)} \ \ \ \mbox{and} \ \  \ \lambda_2(t) =  -\sqrt{-q(t)}.
\label{introduction:eigentwo}
\end{equation}
In particular,  (\ref{introduction:two}) is nondegenerate on any interval on which $q$ does not vanish,
and $t_0$ is a turning point  if and only if  $q$ has an isolated zero there.
The algorithm of \cite{BremerPhase} can only be applied on intervals $I$ on which (\ref{introduction:two})
is nondegenerate.  It uses a continuation scheme of sorts to compute the values
of the derivatives $r_1$ and $r_2$ of the desired slowly-varying phase functions $\psi_1$ and $\psi_2$
at a particular point $c$ in $I$.  The Riccati equation is then solved using 
$r_1(c)$ and $r_2(c)$ as initial conditions in order to calculate $r_1$ and $r_2$
over the whole interval.  The slowly-varying phase functions are obtained by integration.

A simple modification of the algorithm of \cite{BremerPhase}, described in \cite{BremerPhase2},  extends it to the case in 
which $q$ is nonzero on an interval $[a,b]$, except at a finite number of isolated zeros;
that is, when (\ref{introduction:two}) is nondegenerate over $[a,b$] with a finite number of turning points.
The algorithm operates by introducing a  partition $a=\xi_1 < \xi_2 < \ldots < \xi_k=b$ of 
$[a,b]$ such that $\xi_2,\ldots,\xi_{k-1}$ are the roots of $q$ in the open
interval $(a,b)$.  Then, for each $j=1,\ldots,k-1$, the algorithm constructs
two slowly-varying phase functions $\psi_1^{j}$ and $\psi_2^{j}$
such that 
\begin{equation}
\exp\left(\psi_1^{j}(t)\right) \ \ \ \mbox{and}\ \ \ \exp\left(\psi_2^{j}(t)\right)
\end{equation}
comprise a basis in the space of solutions of (\ref{introduction:two}) given on the interval
$\left[\xi_j,\xi_{j+1}\right]$ using the approach of \cite{BremerPhase}.  It is necessary to partition $[a,b]$ because slowly-varying phase 
functions cannot always be extended across a turning point.
It is relatively straightforward to generalize the method of \cite{BremerPhase,BremerPhase2} to 
higher order scalar equations.  However, that is not the approach we take in this article. 
The authors have found that while the resulting algorithm  is highly-effective in the case of a large collection
of equations of the form (\ref{introduction:scalarode}), it has properties very similar 
to another approach discussed in this paper, but with slightly inferior performance.  

In this article, we describe two methods for 
constructing a collection $\psi_1,\ldots,\psi_n$ of slowly-varying phase functions  such that  (\ref{introduction:phaserep}) 
is a basis in the space of solutions of (\ref{introduction:scalarode}).
Both  are inspired by the classical Levin approach to the  numerical evaluation of oscillatory integrals.
Introduced in  \cite{Levin}, the Levin method is based on the observation that inhomogeneous equations of the form
\begin{equation}
y'(t) + p_0(t) y(t)  = f(t)
\end{equation}
admit solutions whose complexity depends on that of $p_0$ and $f$, but not on the magnitude
of $p_0$.  This principle extends to the case of equations of the form
\begin{equation}
y^{(n)}(t) + p_{n-1}(t) y^{(n-1)}(t) + \cdots + p_1(t) y'(y) + p_0(t) y(t) = f(t).
\label{introduction:inhom}
\end{equation}
That is, such equations admit solutions whose complexity depends on that of the right-hand side $f$ and of
the coefficients $p_0,\ldots,p_{n-1}$, but not on the magnitudes of the $p_0,\ldots,p_{n-1}$.
We exploit this principle by applying Newton's method to the Riccati equation for (\ref{introduction:scalarode}).
Starting with a slowly-varying initial guess ensures that each of 
linearized equations which arise are of the form (\ref{introduction:inhom}), and so admit a slowly-varying
solution.    Consequently, a slowly-varying solution of the Riccati equation can be constructed
via Newton's method as long as an appropriate slowly-varying initial guess is known.
Conveniently enough, there is an obvious mechanism for generating $n$ slowly-varying
initial guesses for the solution of the $(n-1)^{st}$ order Riccati equation.  In particular,
the eigenvalues $\lambda_1(t),\ldots,\lambda_n(t)$ of the matrix (\ref{introduction:scalarcoef}),
which are often used as low-accuracy approximations of solutions
of the Riccati equation in asymptotic methods, 
are suitable as  initial guesses for the Newton procedure.

Complicating matters is the fact that the differential operator 
\begin{equation}
D\left[y\right](t) = y^{(n)}(t) + p_{n-1}(t) y^{(n-1)}(t) + \cdots + p_1(t) y'(y) + p_0(t) y(t)
\label{introduction:dop}
\end{equation}
appearing on the left-hand side of   (\ref{introduction:inhom}) 
admits a nontrivial nullspace comprising all solutions of the homogeneous equation
\begin{equation}
y^{(n)}(t) + p_{n-1}(t) y^{(n-1)}(t) + \cdots + p_1(t) y'(y) + p_0(t) y(t) = 0.
\label{introduction:hom}
\end{equation}
This means, of course, that (\ref{introduction:inhom}) is not uniquely solvable.
But it also implies that most solutions of (\ref{introduction:inhom}) are rapidly-varying when the 
coefficients  $p_0,\ldots,p_{n-1}$  are of large magnitude
since the homogeneous equation (\ref{introduction:hom}) admits rapidly-varying solutions in such cases.
It is observed in \cite{Levin} that when the solutions of (\ref{introduction:hom})
are all  rapidly-varying but (\ref{introduction:inhom}) admits a slowly-varying solution $y_0$,
a simple spectral collocation method can be used to compute $y_0$ provided some care is taken in choosing  the discretization
grid.   In particular, if the collocation grid is sufficient to resolve the slowly-varying solution $y_0$, 
but not the solutions of (\ref{introduction:hom}), then the 
matrix discretizing (\ref{introduction:dop}) will be well-conditioned and 
inverting it will yield $y_0$.


We refer to the first of the algorithms introduced here as the global Levin method.    
It operates by subdividing $[a,b]$
until, on each resulting subinterval, every one of the slowly-varying phase functions
$\psi_1,\ldots,\psi_n$ it constructs is  accurately represented by a  Chebyshev expansion of a fixed order. The derivatives
$r_1,\ldots,r_n$  of the phase functions $\psi_1,\ldots,\psi_n$
are calculated on each subinterval by applying Newton's method to the Riccati equation and solving the resulting linearized equations via a 
spectral collocation method.  The method is highly accurate and effective, as long as the eigenvalues
 $\lambda_1(t), \ldots, \lambda_n(t)$  of the coefficient matrix (\ref{introduction:scalarcoef})
are of large magnitude on  the entire interval $[a,b]$.    In this case, each slowly-varying
solution of the Riccati equation is isolated from the others  so that 
the linear operators which arise during the Newton procedure have no slowly-varying
elements in their nullspaces.  Consequently, the spectral collocation method used to solve
the linearized equations always yields a unique solution and the Newton iterates
converge with no difficulty to a slowly-varying phase function.

When one or more of the eigenvalues $\lambda_j(t)$ is of small magnitude on some part of $[a,b]$, 
the global Levin method can fail.  In this event, there are multiple slowly-varying
solution of the Riccati equation near each initial guess.  Accordingly, the Newton iterates
can converge to different slowly-varying solutions of the Riccati equation, and because
the global Levin method is an adaptive approach which subdivides the interval $[a,b]$,
it can obtain  different solutions on  different subintervals.  In particular, the functions obtained from the global Levin method
in such cases can be discontinuous across the boundaries of the discretization subintervals, even
though they locally satisfy the Riccati equation on each particular subinterval (this effect can clearly be seen 
in Figure~\ref{experiments6:figure2} of Subsection~\ref{section:experiments:6}.)

The second algorithm of this paper, which we call the local Levin method,
is usually  slower and slightly less accurate than the global Levin method in cases in which both apply, but it 
overcomes the difficulties  which arise when the coefficient matrix (\ref{introduction:scalarcoef}) admits
eigenvalues  
of  small magnitude.  The local method operates in a manner  very similar to the algorithm of \cite{BremerPhase}
---  it first computes the values of  the derivatives $r_1, r_2,\ldots, r_n$ of the
desired slowly-varying phase functions $\psi_1,\ldots,\psi_n$  at a point $c$ and then 
solves the Riccati equation numerically
with those values used as initial conditions
in order to  construct $r_1,r_2, \ldots,r_n$ over the whole interval $[a,b]$. 
Rather than the continuation method of  \cite{BremerPhase}, however, the values
of the $r_1, r_2,\ldots, r_n$ at $c$ are computed by applying the   Levin approach to a single 
small subinterval of $[a,b]$ containing $c$.  Because the Levin approach is only used on a single subinterval of $[a,b]$, the nonuniqueness
of the slowly-varying phase functions is no longer a concern --- any set of slowly-varying phase
functions suffices for our purposes.
For the sake of simplicity, we describe the local Levin method
in the case in which  (\ref{introduction:scalarode})
is nondegenerate on the entire interval $[a,b]$ of interest.  However, it can easily be extended
to equations of the form (\ref{introduction:scalarode}) that are nondegenerate on $[a,b]$
except at a finite number of turning points by partitioning $[a,b]$ and applying it
on each of the resulting subintervals in the manner of  \cite{BremerPhase2}.

Our algorithms bear some superficial similarities to Magnus expansion methods, which
are another class of numerical solvers for systems of ordinary differential equations 
that rely on the  exponential representation of solutions.  
Introduced in \cite{Magnus},    Magnus expansions are certain series of the form
\begin{equation}
\sum_{k=1}^\infty \Omega_k(t)
\label{introduction:Magnus}
\end{equation}
such that  $\exp\left(\sum_{k=1}^\infty \Omega_k(t)\right)$
locally represents a fundamental matrix for a system of  differential equations
\begin{equation}
\mathbf{y}'(t) = A(t) \mathbf{y}(t).
\label{introduction:system}
\end{equation}
The first few terms for the series around $t=0$ are given by 
\begin{equation}
\begin{aligned}
\Omega_1(t) &=  \int_{0}^t A(s)\ ds,\\
\Omega_2(t) &= \frac{1}{2} \int_{0}^t \int_{0}^{t_1} \left[A(t_1),A(t_2)\right]\, dt_2dt_1 \ \ \ \mbox{and} \\
\Omega_3(t) &= \frac{1}{6} \int_{0}^t \int_{0}^{t_1} \int_{0}^{t_2} 
\left[A(t_1), \left[A(t_2),A(t_3)\right]\right] +  \left[A(t_3), \left[A(t_2),A(t_1)\right]\right]\,  dt_3dt_2dt_1.
\end{aligned}
\label{introduction:magnusterms}
\end{equation}
The straightforward evaluation of the $\Omega_j$ is nightmarishly expensive; however, a clever
technique which renders the calculations manageable is introduced in \cite{Iserles101}
and it paved the way for the development of a class of  numerical solvers 
which represent a fundamental matrix for 
(\ref{introduction:system}) over an interval $I$ via a collection of truncated
Magnus expansions. While the entries of the $\Omega_j$ are slowly-varying whenever the entries of $A(t)$ are slowly-varying, the 
radius of convergence of the series  in (\ref{introduction:Magnus})
depends on the magnitude of the coefficient matrix $A(t)$, which is, in turn,
related to the magnitudes of the eigenvalues of $A(t)$.
Of course, this means that the number of Magnus expansions
which are needed, and hence the cost of the method,
depends on the  magnitudes of the eigenvalues of $A(t)$.
See, for instance, \cite{Iserles102} which gives for estimates
of the growth in the running time of Magnus expansion methods in 
the case of an equation of the form (\ref{introduction:two}) as 
a function of the magnitude of the coefficient $q$.

The algorithms of this paper, by contrast,
represent a  fundamental matrix for (\ref{introduction:system}) with
$A(t)$ given by (\ref{introduction:scalarcoef})
in the form 
\begin{equation}
\exp\left(
\begin{array}{ccccc}
\psi_1(t)          &          &             &               & \\
                & \psi_2(t)   &             &               & \\
                &          & \ddots      &               & \\
                &          &             &   \psi_n(t)             & \\
\end{array}
\right),
\label{introduction:expdiag}
\end{equation}
where each $\psi_j$ is slowly-varying and can be represented at a cost independent of the 
magnitudes of the eigenvalues of $A(t)$.
Interestingly,  as a consequence of standard uniqueness results for ordinary
differential equations, in the case in which $A(t)$ is of the scalar form (\ref{introduction:scalarcoef})
and $\psi_1(0)=\psi_2(0)=\cdots=\psi_n(0)=1$,
the Magnus expansion
\begin{equation}
\exp\left(\sum_{k=1}^\infty \Omega_k(t) \right)
\end{equation}
for (\ref{introduction:system}) around $t=0$ must converge  to  (\ref{introduction:expdiag}).
It follows that 
\begin{equation}
\sum_{k=1}^\infty \Omega_k(t)
\end{equation}
is the sum of the  diagonal matrix whose nonzero entries are $\psi_1(t),\ldots,\psi_n(t)$
and a logarithm of the identity matrix.  In particular, Magnus expansions converge to  a matrix which can be represented
at a cost independent of the magnitudes even though the  individual terms in the expansion 
do not have this property.

Ostensibly, Magnus expansion methods apply to a much large class of systems of differential
equations than the algorithms we describe here.  However, since essentially any system
of ordinary differential equations can be converted into scalar form (see, for instance, \cite{Kolchin}),
it is possible to use our methods to solve a large class of systems of linear ordinary
differential equations in time independent of the magnitudes of their coefficient matrices.
We give one such example in the numerical experiments of this paper (see Section~\ref{section:experiments:7})
and leave a description of this approach to a future  work.  We note, however, that 
an approach of this type involves representing a fundamental matrix for a system
of differential equations (\ref{introduction:system}) in the form
\begin{equation}
\Phi(t) 
\exp\left(
\begin{array}{ccccc}
\psi_1(t)          &          &             &               & \\
                & \psi_2(t)   &             &               & \\
                &          & \ddots      &               & \\
                &          &             &   \psi_n(t)             & \\
\end{array}
\right),
\end{equation}
where $\Phi$ is an appropriately chosen slowly-varying transformation matrix
and the $\psi_1(t),\ldots,\psi_n(t)$ are slowly-varying phase functions
for a scalar equation.

The remainder of this article is organized as follows.  In Section~\ref{section:oneint}
we detail the procedure used by the local and global Levin method to construct the derivatives of the slowly-varying
phase functions over a single subinterval of the interval $[a,b]$ over which (\ref{introduction:scalarode})
is given.    A description of the global Levin method appears
in Section~\ref{section:global}, while Section~\ref{section:local} details the local Levin method.    
The results of numerical experiments demonstrating 
the properties of these algorithms are discussed in Section~\ref{section:experiments}.  We comment on the algorithms of  this article and give
a few brief suggestions for future work in Section~\ref{section:conclusion}.  Appendix~\ref{section:appendix}
details a standard adaptive spectral solver for ordinary differential equations which is used by the local Levin
method, and to construct reference solutions in our numerical experiments.

\end{section}

\begin{section}{The Levin procedure for a single subinterval}
\label{section:oneint}

In this section, we describe the procedure used by the local and  global Levin methods to
construct phase functions on a single subinterval $[c,d]$ of the domain $[a,b]$
over which the equation (\ref{introduction:scalarode}) is given.
The procedure takes as input the following:
\begin{enumerate}
\item
the subinterval $[c,d]$;

\item
an external  subroutine for evaluating the coefficients $q_0,\ldots,q_{n-1}$ in (\ref{introduction:scalarode}); and

\item
an integer $k$ which controls the order of the Chebyshev expansions used to represent phase functions over $[c,d]$.


\end{enumerate}
For each $j=1,\ldots,n$,  it outputs a $(k-1)^{st}$ order Chebyshev expansion which represents
the derivative $r_j$ of the phase function $\psi_j$ 
over the interval $[c,d]$.  The procedure operates  as follows:

\begin{enumerate}

\item
Construct the $k$-point extremal Chebyshev grid $t_1,\ldots,t_{k}$ on the
interval $[c,d]$
and the corresponding $k \times k$ Chebyshev spectral differentiation matrix $D$.
The nodes are given by the formula
\begin{equation}
t_j = \frac{d-c}{2} \cos\left(\pi \frac{n-j}{n-1}\right) + \frac{d+c}{2}.
\label{oneint:nodes}
\end{equation}
The matrix $D$ takes the vector of values 
\begin{equation}
\left(
\begin{array}{c}
f(t_1)\\
f(t_2)\\
\vdots\\
f(t_k)\\
\end{array}
\right)
\end{equation}
of a Chebyshev expansion of the form
\begin{equation}
f(t) = \sum_{j=0}^{k-1} a_j T_j\left(\frac{2}{d-c} t + \frac{d+c}{d-c} \right)
\end{equation}
to the vector 
\begin{equation}
\left(
\begin{array}{c}
f'(t_1)\\
f'(t_2)\\
\vdots\\
f'(t_k)\\
\end{array}
\right)
\end{equation}
of the values of its derivatives at the nodes $t_1,\ldots,t_j$.

\item
Evaluate the coefficients $q_0,\ldots,q_{n-1}$ at the points $t_1,\ldots,t_n$
by calling the external subroutine supplied by the user.

\item
Calculate the values of $n$ initial guesses $r_{1},\ldots,r_{n}$ for the Newton procedure 
at the nodes $t_1,\ldots,t_n$ by computing the eigenvalues
 of the coefficient matrices
\begin{equation}
A_j = \left(
\begin{array}{cccccccc}
0 & 1 & 0 & \cdots & 0 & 0 \\
0 & 0 & 1 & \cdots & 0 & 0 \\
\vdots &  &  & \ddots &  & \vdots\\
0 & 0 & 0 & \cdots & 1 & 0 \\
0 & 0 & 0 & \cdots & 0 & 1 \\
-q_0(t_j) & -q_1(t_j) & -q_2(t_j) & \cdots&  -q_{n-2}(t_j) & -q_{n-1}(t_j)\\
\end{array}
\right),\ \ \ j=1,\ldots,k.
\label{oneint:aj}
\end{equation}
More explicitly, the eigenvalues of $A_j$ give the values of $r_1(t_j),\ldots,r_n(t_j)$.

\item
Perform Newton iterations in order to refine each of the initial guesses $r_1,\ldots,r_n$.
Because the general form of the Riccati equation is quite complicated, we illustrate
the procedure when $n=2$, in which case
the Riccati equation is simply
\begin{equation}
r'(t) + (r(t))^2 + q_1(t) r(t) + q_0(t) = 0.
\end{equation}
In each iteration, we perform the following steps:

\begin{enumerate}

\item
Compute the residual 
\begin{equation}
\xi(t) = 
r'(t) + (r(t))^2 + q_1(t) r(t) + q_0(t) 
\end{equation}
of the current guess at the nodes $t_1,\ldots,t_k$.

\item
Form a spectral discretization of the linearized operator
\begin{equation}
L\left[\delta\right](t) = \delta'(t) + 2 r(t) \delta(t) + q_1(t) \delta(t).
\end{equation}
That is, form the $k\times k$ matrix
\begin{equation}
B = D + 
\left(\begin{array}{ccccc}
2r(t_1) + q_1(t_1) &                 &  &  \\
       & 2r(t_2) + q_1(t_2) &        &  & \\
       &        & \ddots &        &  \\  
       &         &         & 2r(t_k) + q_1(t_k)
\end{array}\right).
\end{equation}

\item
Solve the $k\times k$ linear system
\begin{equation}
B 
\left(
\begin{array}{c}
\delta(t_1)\\
\delta(t_2)\\
\vdots\\
\delta(t_k)
\end{array}
\right)
 = 
-\left(
\begin{array}{c}
\xi(t_1)\\
\xi(t_2)\\
\vdots\\
\xi(t_k)
\end{array}
\right)
\label{oneint:system}
\end{equation}
and update the current guess:
\begin{equation}
\left(
\begin{array}{c}
r(t_1)\\
r(t_2)\\
\vdots\\
r(t_k)
\end{array}
\right)
=
\left(
\begin{array}{c}
r(t_1)\\
r(t_2)\\
\vdots\\
r(t_k)
\end{array}
\right)
+
\left(
\begin{array}{c}
\delta(t_1)\\
\delta(t_2)\\
\vdots\\
\delta(t_k)
\end{array}
\right).
\end{equation}

\end{enumerate}

We perform a maximum of $8$ Newton iterations and the procedure is terminated if the value of
\begin{equation}
\sum_{j=1}^n \left|\delta(t_j)\right|^2
\end{equation}
is smaller than 
\begin{equation}
\left(100 \epsilon_0\right)^2\, 
\sum_{j=1}^n \left|r(t_j)\right|^2,
\end{equation}
where $\epsilon_0\approx 2.220446049250313\times 10^{-16}$ denotes machine zero for the IEEE double precision number
system.

\item
Form the Chebyshev expansions of the functions $r_1,\ldots,r_n$, which constitute
the outputs of this procedure.
\end{enumerate}

Standard eigensolvers often produce inaccurate results in the case of 
matrices of the form (\ref{oneint:aj}), particularly when the entries
are of large magnitude.  Fortunately, there are specialized techniques
available for companion matrices, and the matrices appearing in (\ref{oneint:aj})
are simply the transposes of such matrices.  Our implementation of the procedure
of this section uses the stable and highly-accurate technique of \cite{AURENTZ1,AURENTZ2}
to compute the eigenvalues of the matrices (\ref{oneint:aj}).  

Care must also be taken when solving the linear system (\ref{oneint:system})
since the associated operator has a nontrivial nullspace.  Most of the time,
the discretization being used is insufficient to resolve any part of that nullspace,
with the consequence that the matrix $B$ is well-conditioned.  However,
when elements of the nullspace are sufficiently slowly-varying, they can be captured
by the discretization, in which case the matrix $B$ will have small singular values.
Fortunately, it is known that this does not cause numerical difficulties
in the solution of (\ref{oneint:system}), provided a truncated singular value
decomposition is used to invert the system.  Experimental evidence to this effect was presented in \cite{LevinLi,LiImproved}
and a careful analysis of the phenomenon appears in  \cite{SerkhBremerLevin}.
Because the truncated singular value decomposition is quite expensive, we actually use
a rank-revealing QR decomposition to solve the linear system (\ref{oneint:system})
in our implementation of the procedure of this section.  This was found to be about five times faster,
and it lead to no apparent loss in accuracy.

Rather than computing the eigenvalues of each of the matrices (\ref{oneint:aj}) in order
to construct initial guesses for the Newton procedure, one could accelerate
the algorithm slightly by computing the eigenvalues of only one  $A_j$ and use
the constant functions $r_1(t) = \lambda_1(t_j), \ldots, r(t) = \lambda_n(t_j)$ 
as initial guesses instead.  We did not make use of this optimization in our
implementation of the algorithm of this paper.  Our aim was to produce
a reference code which is as robust as possible, not to produce the fastest
code possible.

\end{section}

\begin{section}{The global Levin method}
\label{section:global}

We now describe the global Levin method for the construction of 
a collection of slowly-varying phase functions $\psi_1,\ldots,\psi_n$
such that (\ref{introduction:phaserep}) is a basis in the space
of solutions of an equation of the form (\ref{introduction:scalarode}).
It applies in  the case in which the eigenvalues
$\lambda_1(t), \ldots, \lambda_n(t)$  of the coefficient
matrix (\ref{introduction:scalarcoef}) are of large magnitude
on the entire interval $[a,b]$ of interest.    Throughout this section, we denote
the derivatives of the phase functions $\psi_1,\ldots,\psi_n$
by $r_1,\ldots,r_n$.

The algorithm takes as input the following:
\begin{enumerate}
\item
the interval $[a,b]$ over which the equation is given;

\item
an external subroutine for evaluating the coefficients $q_0,\ldots,q_{n-1}$ in (\ref{introduction:scalarode});

\item
a point $\eta$ on the interval $[a,b]$ and the desired values $\psi_1(\eta),\ldots,\psi_n(\eta)$ for the phase
functions at that point;

\item
an integer $k$ which controls the order of the piecewise Chebyshev expansions used to represent phase functions;  and

\item

a parameter $\epsilon$ which specifies the desired accuracy for the phase functions.
\end{enumerate}

It outputs $(k-1)^{st}$ order  piecewise Chebyshev expansions representing the phase functions $\psi_1,\ldots,\psi_n$ 
on the interval $[a,b]$. By a  $(k-1)^{st}$ order piecewise Chebyshev 
expansions  on the interval $[a,b]$, we mean  a sum of the form
\begin{equation}
\begin{aligned}
&\sum_{i=1}^{m-1} \chi_{\left[x_{i-1},x_{i}\right)} (t) 
\sum_{j=0}^{k-1} \lambda_{ij}\ T_j\left(\frac{2}{x_{i}-x_{i-1}} t + \frac{x_{i}+x_{i-1}}{x_{i}-x_{i-1}}\right)\\
+
&\chi_{\left[x_{m-1},x_{m}\right]} (t) 
\sum_{j=0}^{k-1} \lambda_{mj}\ T_j\left(\frac{2}{x_{m}-x_{m-1}} t + \frac{x_{m}+x_{m-1}}{x_{m}-x_{m-1}}\right),
\end{aligned}
\label{global:chebpw}
\end{equation}
where $a = x_0 < x_1 < \cdots < x_m = b$ is a partition of $[a,b]$,
$\chi_I$ is the characteristic function on the interval $I$ and 
$T_j$ is the Chebyshev polynomial of degree $j$.
We note that terms appearing in the first line of (\ref{global:chebpw}) involve
the characteristic function of a half-open interval, while that appearing
in the second involves  the characteristic function of a closed interval.
This ensures that exactly one term in  (\ref{global:chebpw}) is nonzero for each point $t$ in $[a,b]$.

The algorithm maintains two lists of subintervals of $[a,b]$: one of 
``accepted subintervals'' and one of subintervals which need to be 
processed.  Initially, the list of accepted subintervals is empty 
and the list  of subintervals to be processed contains $[a,b]$.
The following steps are repeated as long as the list of subintervals
to process is nonempty:
\begin{enumerate}

\item
Remove a subinterval $[c,d]$ from the list of subintervals
to process.

\item
Apply the procedure of Section~\ref{section:oneint} in order to construct
 $(k-1)^{st}$ order Chebyshev expansions 
\begin{equation}
r_j(t) \approx  \sum_{i=0}^{k-1} a_{i}^j T_i\left(\frac{2}{d-c} t + \frac{d+c}{d-c} \right),\ \ \ \ j=1,\ldots,n,
\label{global:rj}
\end{equation}
which purportedly represent $r_1,\ldots,r_n$ over $[c,d]$.  The user-supplied external subroutine and integer $k$
are passed as input parameters to the procedure of Section~\ref{section:oneint}.

\item
For each $j=1,\ldots,n$, calculate the quantity
\begin{equation}
\xi_j = 
\frac{ 
\sum_{i=\left\lceil \frac{k+1}{2}\right\rceil}^{k-1} \left|a_i^j \right|^2
}
{
\sum_{i=0}^{k-1} \left|a_i^j \right|^2.
}
\end{equation}

\item
If
\begin{equation}
\xi=\max\left\{\xi_1,\ldots,\xi_n\right\}
\end{equation}
is less than the user-specified parameter $\epsilon$, put the interval $[c,d]$ into the list of accepted intervals
and use (\ref{global:rj}) as the local Chebyshev expansions of the derivatives $r_1,\ldots,r_n$
over the subinterval $[c,d]$.

\item
Otherwise,  if $\xi \geq \epsilon$, put the intervals
\begin{equation}
\left[c,\frac{c+d}{2}\right]
\ \ \ \mbox{and}\ \ \ 
\left[\frac{c+d}{2},d\right]
\end{equation}
into the list of intervals to process.

\end{enumerate}

Upon termination of this procedure, we have $(k-1)^{st}$ order piecewise Chebyshev expansions
representing the derivatives $r_1,\ldots,r_n$ of the phase functions.  The  list of accepted
subintervals determines the partition of $[a,b]$ used by the piecewise expansions.  The phase functions $\psi_1,\ldots,\psi_n$ 
themselves are constructed via spectral integration with the particular choice of 
antiderivatives determined through the parameter $\eta$ and the  values of 
$\psi_1(\eta),\ldots,\psi_n(\eta)$ that are specified as inputs to the algorithm.

\end{section}

\begin{section}{The local Levin method}
\label{section:local}

In this section, we describe the local Levin method for the construction of 
a collection of slowly-varying phase functions $\psi_1,\ldots,\psi_n$
such that (\ref{introduction:phaserep}) is a basis in the space
of solutions of an equation of the form (\ref{introduction:scalarode}).
We describe it in the case in which the equation is nondegenerate on an interval
$(a,b)$ of interest.  In contrast to the global Levin method,
which requires that the eigenvalues  $\lambda_1(t), \ldots, \lambda_n(t)$  of the coefficient
matrix (\ref{introduction:scalarcoef}) are all of large magnitudes across the whole  interval $[a,b]$,
it functions perfectly when one or more of the eigenvalues
is of small magnitude. 

The algorithm takes as input the following:
\begin{enumerate}
\item
the interval $[a,b]$ over which the equation is given;

\item
an external subroutine for evaluating the coefficients $q_0,\ldots,q_{n-1}$ in (\ref{introduction:scalarode});

\item
a point $\eta$ on the interval $[a,b]$ and the desired values $\psi_1(\eta),\ldots,\psi_n(\eta)$ for the phase
functions at that point;

\item
an integer $k$ which controls the order of the piecewise Chebyshev expansions used to represent phase functions;

\item

a parameter $\epsilon$ which specifies the desired accuracy for the phase functions; and

\item
a subinterval $[a_0,b_0]$ of $[a,b]$ over which the Levin procedure is to be applied and a point $\sigma$
in that interval.

\end{enumerate}

The algorithm proceeds by first applying the procedure of Section~\ref{section:oneint} on the subinterval 
$[a_0,b_0]$ of $[a,b]$ in order to approximate the values of the derivatives $r_1,\ldots,r_n$ of the slowly-varying
phase functions at the point $\sigma \in [a_0,b_0]$.  The parameter $k$ and the user-supplied external subroutine
are passed to that procedure.    

Next, for each $j=1,\ldots,n$, the Riccati equation is solved using the value of $r_j(\sigma)$ to specify the desired solution.
These calculations are performed via the adaptive  spectral method described in Appendix~\ref{section:appendix}.
The parameters $k$ and $\epsilon$ are passed to that procedure.  Since most solutions of the Riccati equation are rapidly-varying
and we are seeking a slowly-varying solution, these problems are extremely stiff.  The solver of Appendix~\ref{section:appendix}
is well-adapted to such problems; however, essentially any solver for stiff ordinary differential equations would serve in its place.
The result is a collection of $(k-1)^{st}$ order piecewise Chebyshev expansions representing the derivatives $r_1,\ldots,r_n$ 
of the phase functions $\psi_1,\ldots,\psi_n$. 

Finally, spectral integration is used to construct the phase functions $\psi_1,\ldots,\psi_n$ from their  derivatives
$r_1,\ldots,r_n$.  The particular  antiderivatives are determined by the values 
$\psi_1(\eta), \ldots, \psi_n(\eta)$ specified as inputs to the algorithm.

\end{section}

\begin{section}{Numerical experiments}
\label{section:experiments}

In this section, we present the results of numerical experiments which were conducted to illustrate the properties
of the algorithms of this paper.  The code for these experiments was written in Fortran and compiled with version 13.1.1
of the GNU Fortran compiler.     
They were performed on a desktop computer equipped  with an AMD 3900X processor and 32GB  of memory.  
No attempt  was made to parallelize our code.

Our algorithms call for computing the eigenvalues of matrices of the
form (\ref{introduction:scalarcoef}).  Unfortunately, standard eigensolvers
lose significant accuracy when applied to many matrices of this type. However, 
because the  transpose of (\ref{introduction:scalarcoef}) is a companion matrix,
we were able to use the backward stable algorithm of  \cite{AURENTZ1,AURENTZ2} 
for computing the eigenvalues of companion matrices to perform these calculations.

We employed two different methods to assess the accuracy of the slowly-varying phase functions obtained by our algorithms.  
When possible, we used them to solve an initial or boundary value  problem for an equation of the form (\ref{introduction:scalarode})
and measured the absolute accuracy of the result by comparison with the output of the standard solver described in 
Appendix~\ref{section:appendix}.    Absolute accuracy was measured rather
than relative accuracy because all of the solutions we calculated were oscillatory on at least some part of 
their domain of definition.
This first approach is not viable when the real parts of the eigenvalues 
are too large since initial and boundary value problems for 
(\ref{introduction:scalarode}) are highly ill-conditioned in this event.
Moreover, the solver of Appendix~\ref{section:appendix} becomes prohibitively expensive
when the sizes of the eigenvalues are excessive, regardless of whether
it is their real parts, their imaginary parts or both which are of large magnitude.

Because of these limitations, with few exceptions, we restricted our use of this approach
to cases in which the parameter $\omega$ used  to control the frequency of oscillation of solutions
was  no larger than $2^{16}$, and we never applied it to a  problem in which the eigenvalues of the coefficient
matrix have real parts of large magnitude.
Moreover, because the condition numbers of initial and boundary value problems for 
(\ref{introduction:scalarode}) grow  with the magnitude of the eigenvalues of (\ref{introduction:scalarcoef}),
it is expected that the accuracy of the solutions obtained via any numerical
approach will deteriorate as the eigenvalues of the coefficient matrix increase.  In particular,
the absolute error in the solutions of (\ref{introduction:scalarode}) obtained by our method
are limited principally by the conditioning of the problem 
and not by the accuracy with which  phase functions are constructed.

Our second method consisted of constructing slowly-varying phase functions
by running one of our algorithms using IEEE quadruple  precision arithmetic and
comparing the results to those obtained by running our algorithms using standard
double precision arithmetic.  This is a highly unsatisfying approach,
but the authors are not aware of another algorithm for the high-accuracy
approximation of slowly-varying phase functions and standard solvers
for ordinary differential equations perform quite poorly when applied to 
most of the problems we consider here.


In all of our experiments,  the input parameters for the algorithms of Sections~\ref{section:global}
and \ref{section:local} were set as follows. The value of $k$, which determines the order of the 
Chebyshev expansions used to represent phase functions, was taken to be $16$ and the parameter
$\epsilon$, which specifies the desired precision for the phase functions, was taken to be $10^{-12}$.
The domain for each equation we considered was $[-1,1]$, and the 
particular antiderivatives $\psi_1,\ldots,\psi_n$ of  the functions $r_1,\ldots,r_n$ 
were chosen through the requirement that  $\psi_1(0) = \psi_2(0) = \cdots = \psi_n(0) = 0$.
In the case of the local Levin method, the 
procedure of Section~\ref{section:oneint} used to determine the initial conditions
for the functions $r_1,\ldots,r_n$  was performed on the  subinterval $[-0.1,-0.0]$.


%
%
%
%
\begin{subsection}{An initial value problem for a second order equation }
\label{section:experiments:1}

The experiments described in this section concern the equation
\begin{equation}
y''(t) -  \frac{i\omega}{1+t^4} y'(t)  + \omega^3 \left( \frac{1+\cos^2(t)}{2+\omega \exp(t)}\right) y(t) = 0.
\label{experiments1:1}
\end{equation}
To give a sense of the dependence of the eigenvalues of the coefficient matrix
for (\ref{experiments1:1})  on the parameter $\omega$, we note that
\begin{equation}
\begin{aligned}
\lambda_1(0) &=
\frac{i \omega }{2}-\frac{i \sqrt{\omega ^2 (\omega +2) (9 \omega
   +2)}}{2 (\omega +2)} \sim -i\omega +\frac{4i}{3} + \mathcal{O}\left(\frac{1}{\omega}\right) \ \ \mbox{as} \ \ \omega\to\infty\ \ \ \mbox{and}\\
\lambda_2(0) &= 
\frac{i \omega }{2}+\frac{i \sqrt{\omega ^2 (\omega +2) (9 \omega
   +2)}}{2 (\omega +2)} \sim 2i\omega -\frac{4i}{3} + \mathcal{O}\left(\frac{1}{\omega}\right) \ \ \mbox{as} \ \ \omega\to\infty.\\\\
\end{aligned}
\end{equation}
Moreover, Figure~\ref{experiments1:figure2}(b) contains 
plots of the eigenvalues of the coefficient matrix for (\ref{experiments1:1}) 
when $\omega=2^{16}$.

Our first experiment proceeded as follows.
For each $\omega=2^8, 2^9, \ldots, 2^{16}$, we used both the global and local Levin methods
to solve (\ref{experiments1:1}) over the interval $[-1,1]$  subject to the condition
\begin{equation}
y(0) = 1 \ \ \ \mbox{and} \ \ \ y'(0) = i \omega.
\label{experiments1:2}
\end{equation}
We then measured the absolute errors in each obtained solution at $10,000$ equispaced points in
the interval $[-1,1]$  via comparison with a reference solution constructed 
using the standard solver described in Appendix~\ref{section:appendix}.  
%
The second experiment was conducted as follows.  For each $\omega=2^8,2^9,\ldots,2^{20}$, we constructed
slowly-varying phase functions for (\ref{experiments1:1}) over the interval $[-1,1]$ by running both the global
and local Levin methods using double precision arithmetic (as usual). 
We then measured the relative errors in each obtained phase function at $10,000$ equispaced points in
the interval $[-1,1]$  by comparison with phase functions constructed by applying the 
global Levin method to (\ref{experiments1:1}), but this time using extended precision arithmetic to perform the calculations.
Figure~\ref{experiments1:figure1} gives the results of these experiments.
Figure~\ref{experiments1:figure2}(a)
contains  plots of the derivatives of the slowly-varying phase functions produced by global Levin method
when $\omega=2^{16}$.

\begin{figure}[h]
\hfil
\includegraphics[width=.40\textwidth]{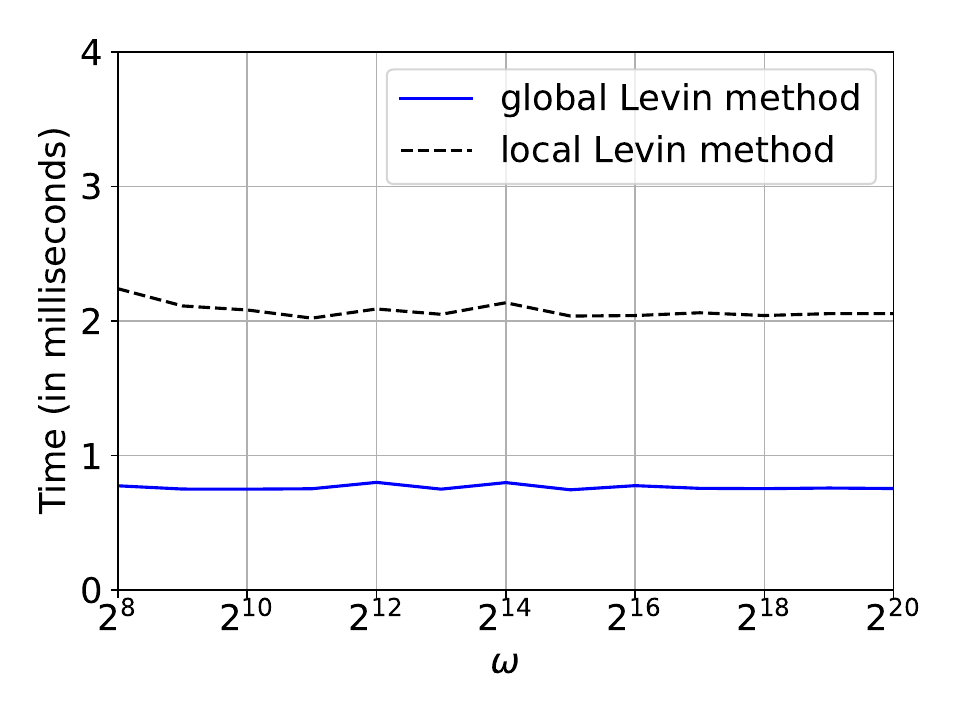}
\hfil
\includegraphics[width=.40\textwidth]{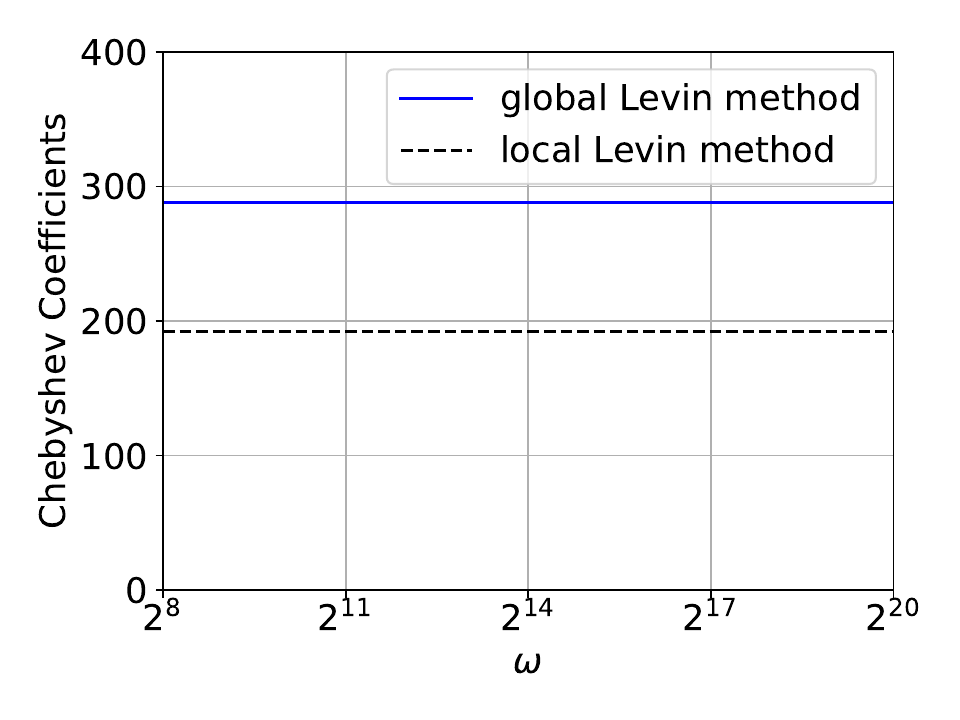}
\hfil

\hfil
\includegraphics[width=.40\textwidth]{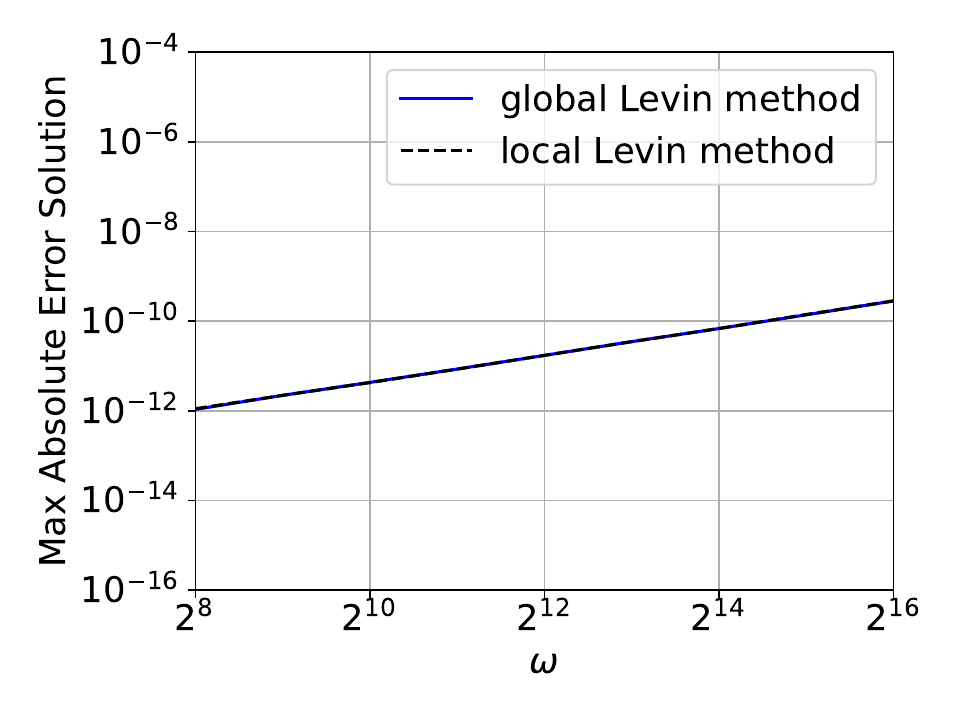}
\hfil
\includegraphics[width=.40\textwidth]{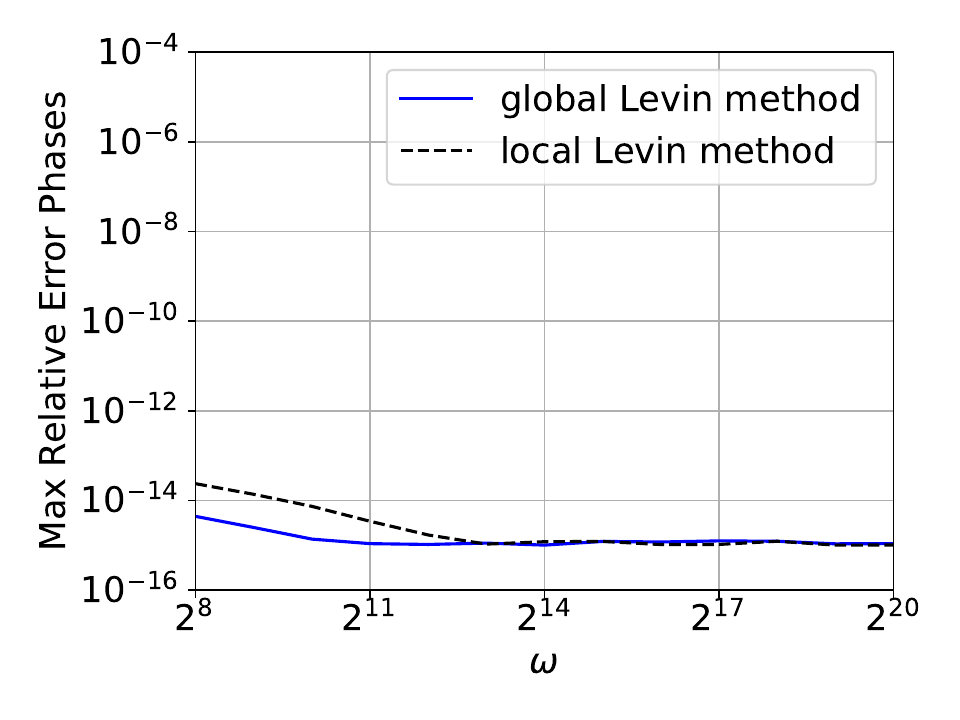}
\hfil

\caption{The results of the experiments of Subsection~\ref{section:experiments:1}.
The  upper-left plot gives the time required by each of our methods
as a function of the parameter $\omega$.  The upper-right plot 
reports the number  of Chebyshev coefficients required to represent the slowly-varying phase functions
as a function of $\omega$.  
The lower-left plot  reports the largest observed
absolute error in the solutions of the initial value problem for (\ref{experiments1:1}) obtained using
each of our methods.    The plot on the lower right gives
the largest observed relative error in the slowly-varying phase functions constructed by our algorithms.
The plots in the upper-right and lower-left appear to have only one line
because the solutions obtained by the local and global Levin methods closely coincide.
}
\label{experiments1:figure1}
\end{figure}

\begin{figure}[h]
\centering
\hfil
\begin{subfigure}{.45\textwidth}
\hfil
\includegraphics[width=.45\textwidth]{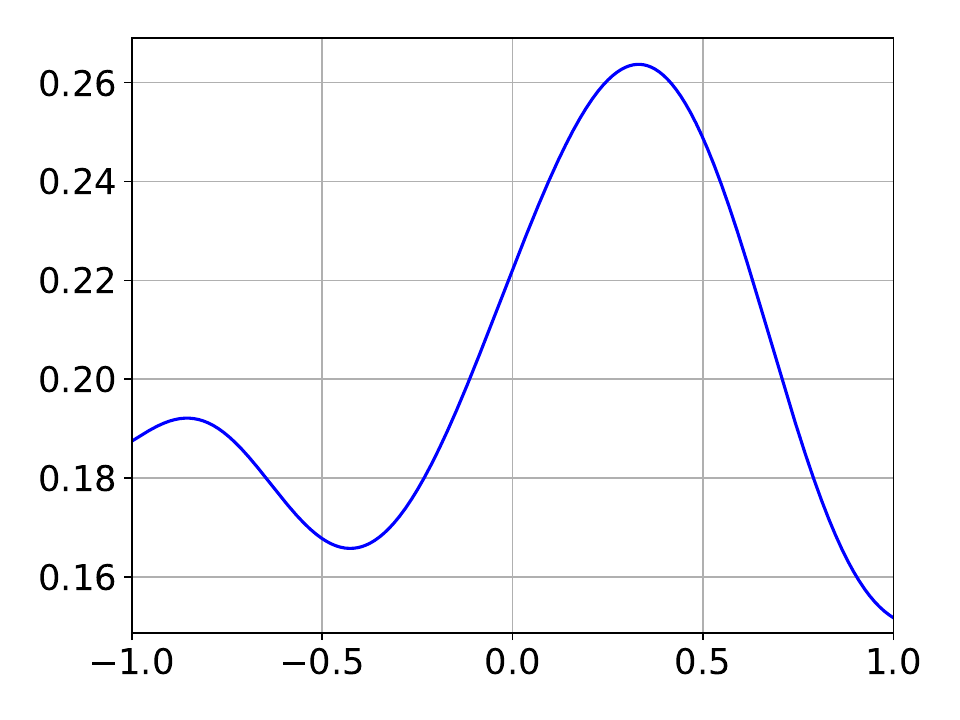}
\hfil
\includegraphics[width=.45\textwidth]{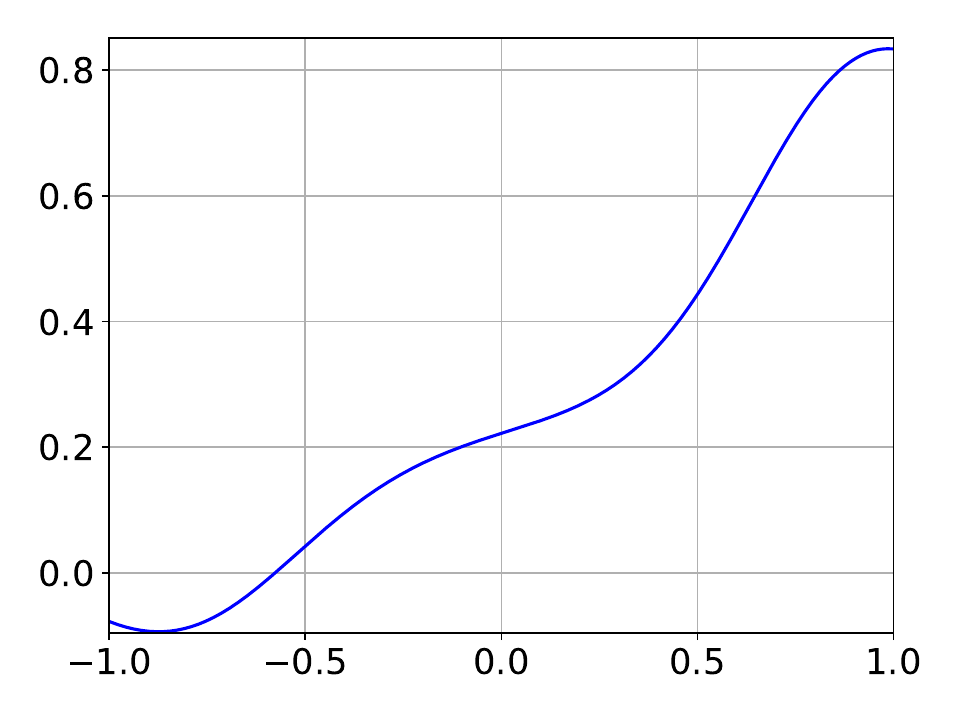}
\hfil

\hfil
\includegraphics[width=.45\textwidth]{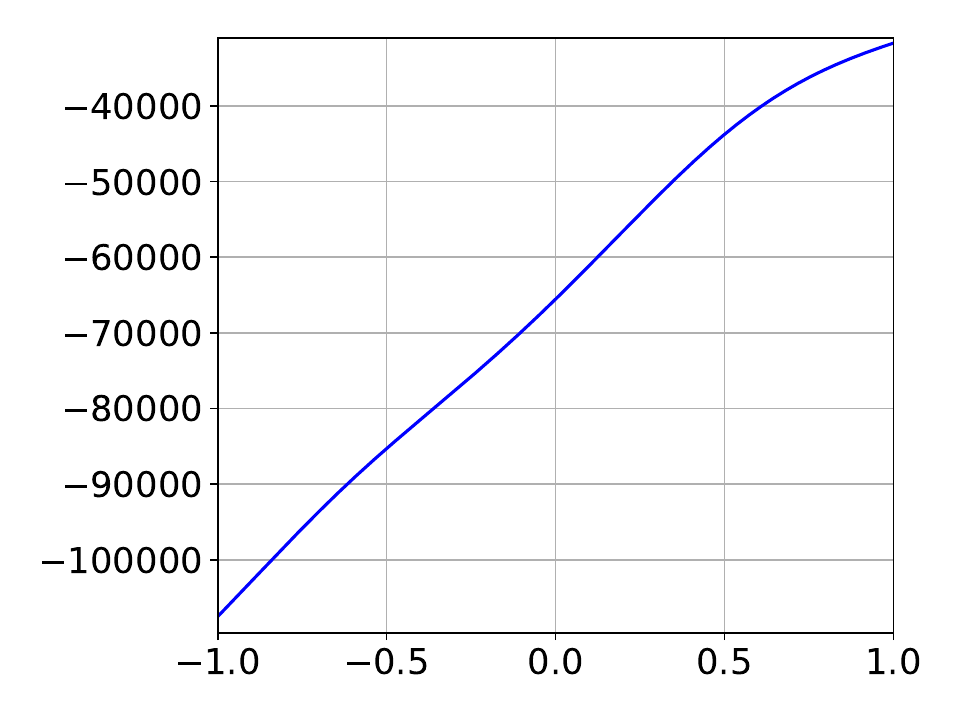}
\hfil
\includegraphics[width=.45\textwidth]{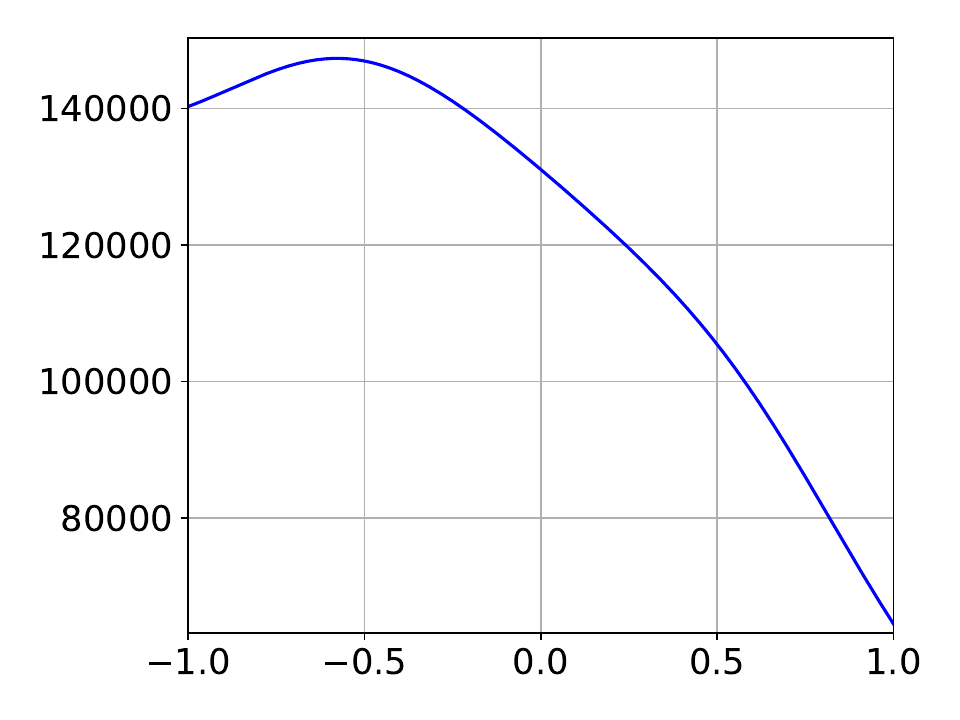}
\hfil
\caption{The derivatives of the two slowly-varying phase functions produced
by applying the global Levin method to Problem~(\ref{experiments1:1}) of Subsection~\ref{section:experiments:1}
when $\omega=2^{16}$.  
Each column corresponds to one of the phase functions,
with the real part appearing in the first row and the imaginary part in the second.}
\end{subfigure}
\hfill
\begin{subfigure}{.45\textwidth}
\hfil
\includegraphics[width=.45\textwidth]{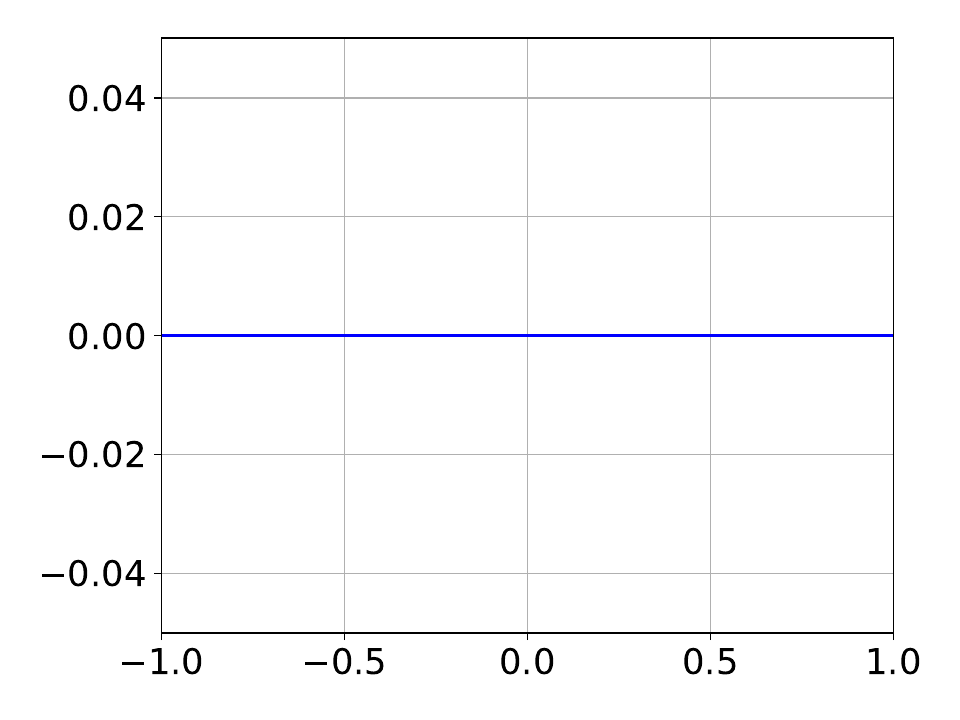}
\hfil
\includegraphics[width=.45\textwidth]{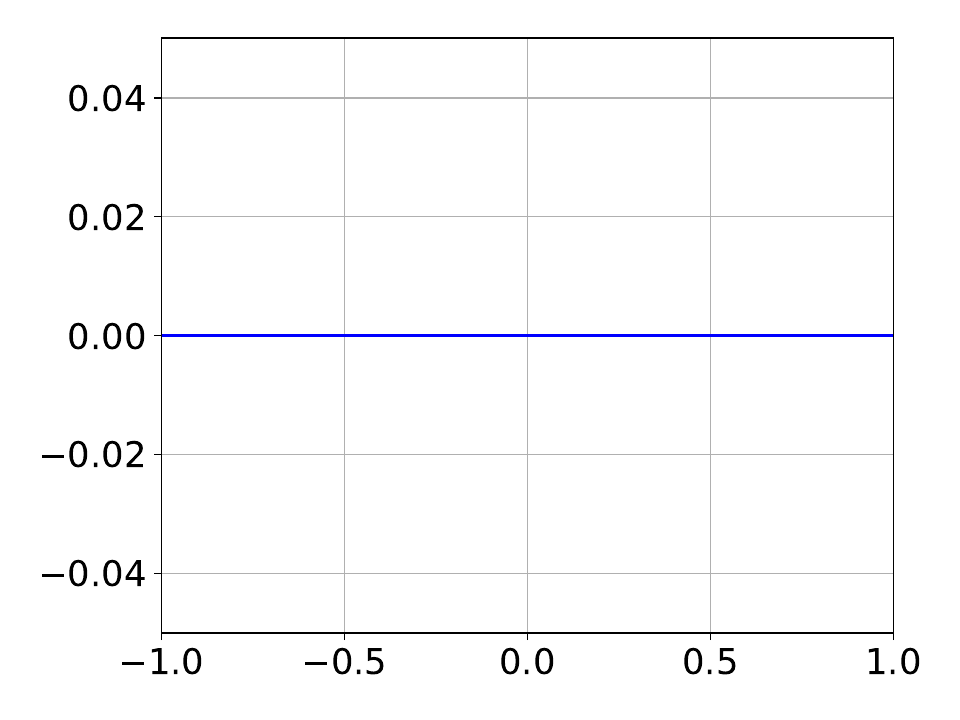}
\hfil

\hfil
\includegraphics[width=.45\textwidth]{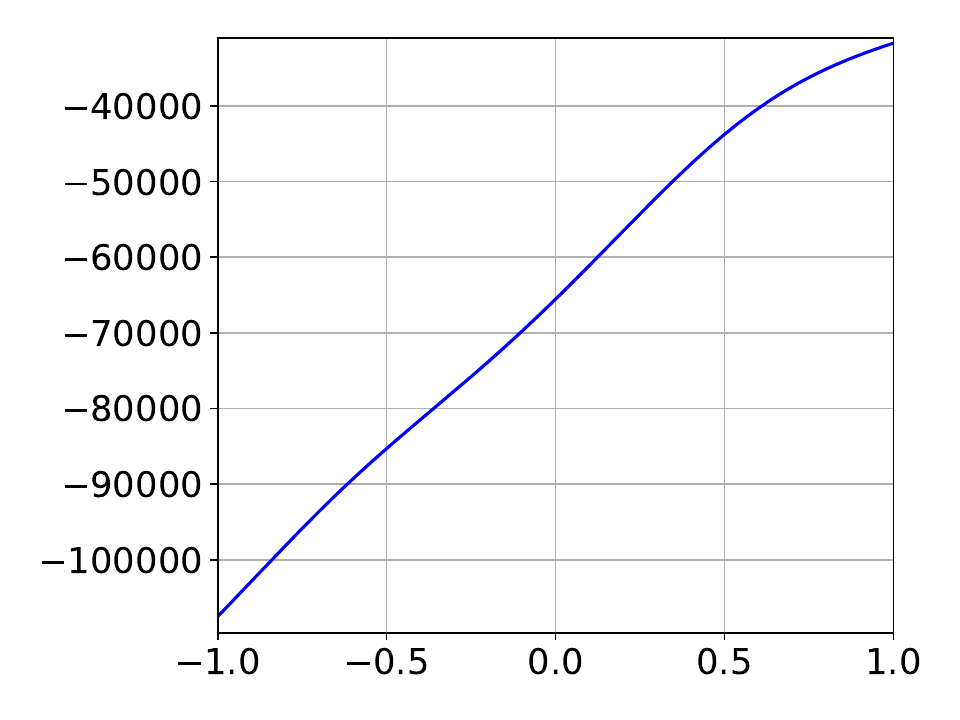}
\hfil
\includegraphics[width=.45\textwidth]{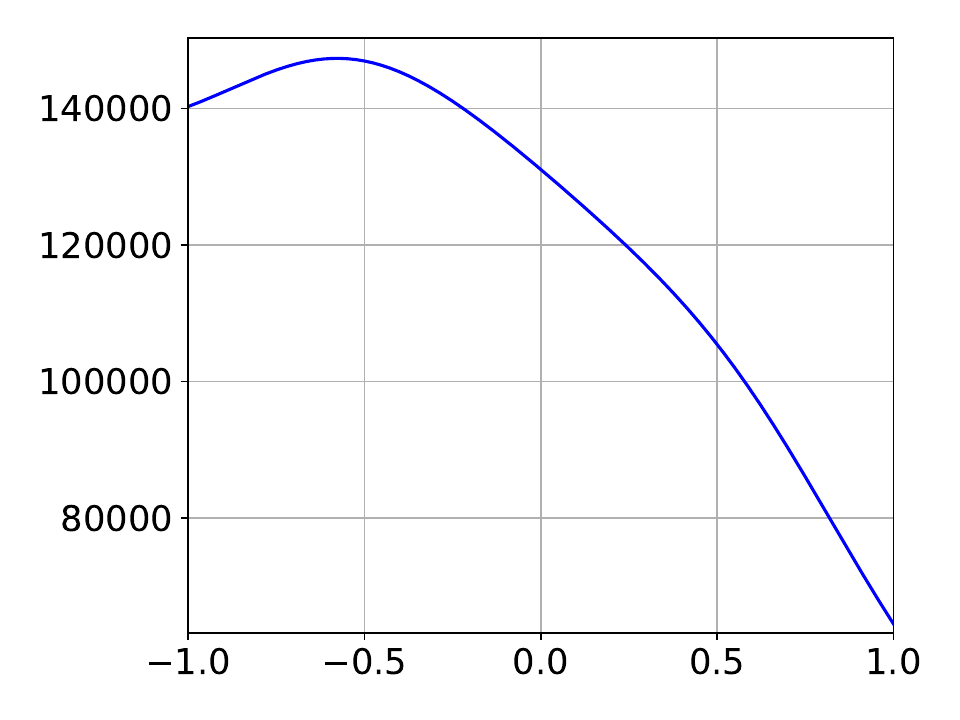}
\hfil
\caption{
The eigenvalues $\lambda_1(t), \lambda_2(t)$ of the 
coefficient  matrix corresponding to Equation~(\ref{experiments1:1})
of Subsection~\ref{section:experiments:1}
when $\omega=2^{16}$.  
Each column corresponds to one of the eigenvalues,
with the real part appearing in the first row and the imaginary part in the second.\\
}
\end{subfigure}
\hfil
\caption{Plots of some of the phase functions and eigenvalues
which arose in the course of the  experiments of Subsection~\ref{section:experiments:1}.}
\label{experiments1:figure2}
\end{figure}

\end{subsection}

%
%
%
\begin{subsection}{An initial value problem for a third order equation }
\label{section:experiments:2}

In this experiments, we considered the equation
\begin{equation}
y'''(t) + q_2(\omega,t) y''(t) + q_1(\omega,t) y'(t) + q_0(\omega,t) y(t) = 0,
\label{experiments2:1}
\end{equation}
where
\begin{equation}
\begin{aligned}
q_0(\omega,t) &=
\frac{4 \omega^3 \left(i \omega \sin ^2(t)+i \omega+\sin (t)\right)}{\left(t^2+1\right)
   \left(\omega e^t+1\right)},
\\
q_1(\omega,t) &= 
\frac{\omega \left(\omega \left(4 \omega t^2+\omega e^t+1\right)+\left(\omega \left(4
   t^2+e^t+4\right)+1\right) \sin (t) (\omega \sin
   (t)-i)\right)}{\left(t^2+1\right) \left(\omega e^t+1\right)}\ \ \ \mbox{and}\ \ \
\\
q_2(\omega,t) &= 
i \omega \left(\frac{4 \omega}{\omega e^t+1}+\frac{1}{t^2+1}-1\right)-i \omega \sin
   ^2(t)-\sin (t).
\end{aligned}
\end{equation}
At the point $0$, the eigenvalues of the coefficient matrix for (\ref{experiments2:1})
are given by the formulas
\begin{equation}
\begin{aligned}
\lambda_1(0) = -\frac{4i\omega^2}{1+\omega},\ \ \
\lambda_2(0) = -i\omega \ \ \ \mbox{and} \ \ \
\lambda_1(0) = i\omega.
\end{aligned}
\end{equation}
Plots of the eigenvalues $\lambda_1(t),\lambda_2(t),\lambda_3(t)$ of the coefficient matrix
for (\ref{experiments2:1}) when $\omega=2^{16}$ can be found in 
Figure~\ref{experiments2:figure3} .

In the first experiment, for each $\omega=2^8,2^9,\ldots,2^{16}$, we 
we used the local Levin method and the global
Levin method to compute solutions of (\ref{experiments2:1}).  We then
 measured the errors in each obtained solution at $10,000$ equispaced points in
the interval $[-1,1]$  by comparison with a reference solution constructed via the standard solver
described in Appendix~\ref{section:appendix}.  
%
In the second experiment, for each $\omega=2^8,2^9,\ldots,2^{20}$, we constructed
slowly-varying phase functions for (\ref{experiments2:1}) over the interval $[-1,1]$ by running both the global
and local Levin methods using double precision arithmetic.
We then measured the relative errors in each obtained phase function at $10,000$ equispaced points in
the interval $[-1,1]$  by comparison with phase functions constructed by applying the 
global Levin method to (\ref{experiments2:1}), but this time using extended precision arithmetic to perform the calculations.
Figure~\ref{experiments2:figure1} gives the results of these two experiments.
Figure~\ref{experiments2:figure2}
contains  plots of the derivatives of the slowly-varying phase functions produced by global Levin method
when $\omega=2^{16}$.

\begin{figure}[h]
\hfil
\includegraphics[width=.40\textwidth]{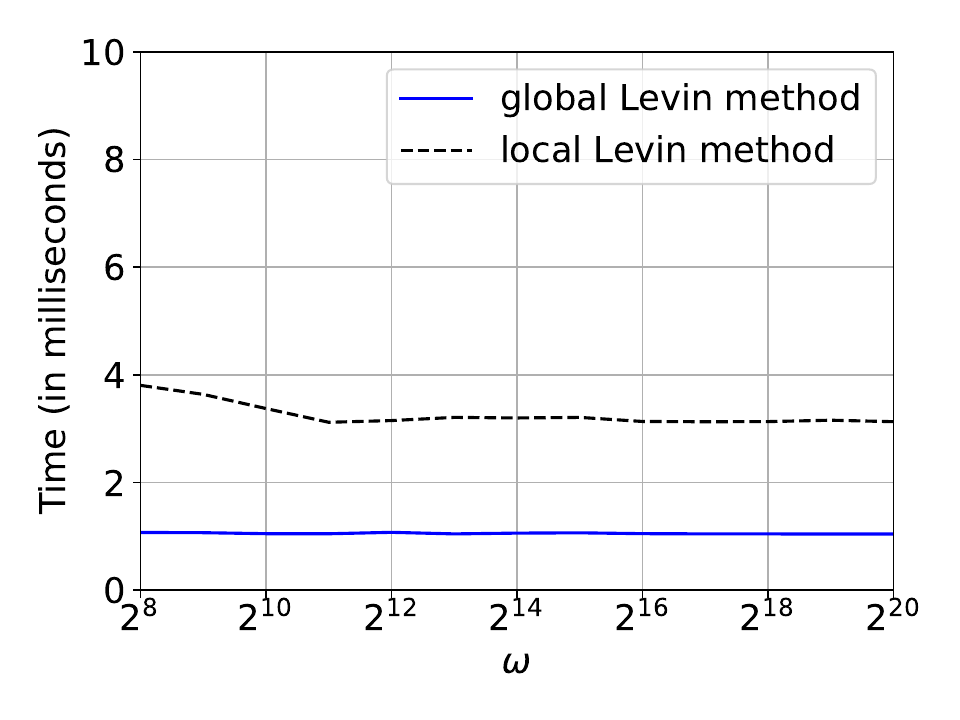}
\hfil
\includegraphics[width=.40\textwidth]{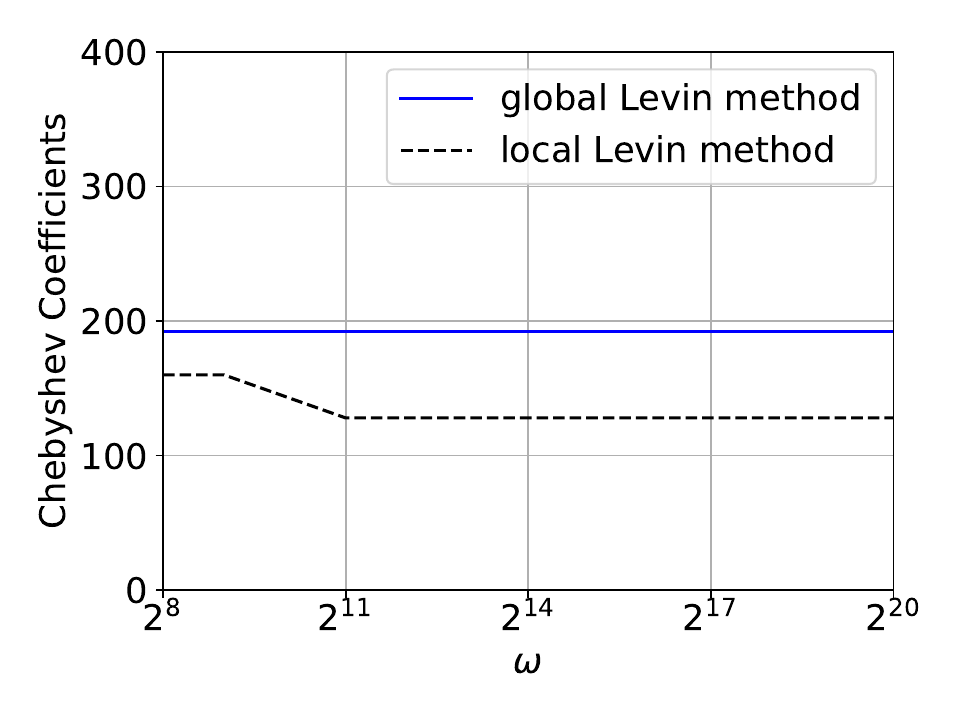}
\hfil

\hfil
\includegraphics[width=.40\textwidth]{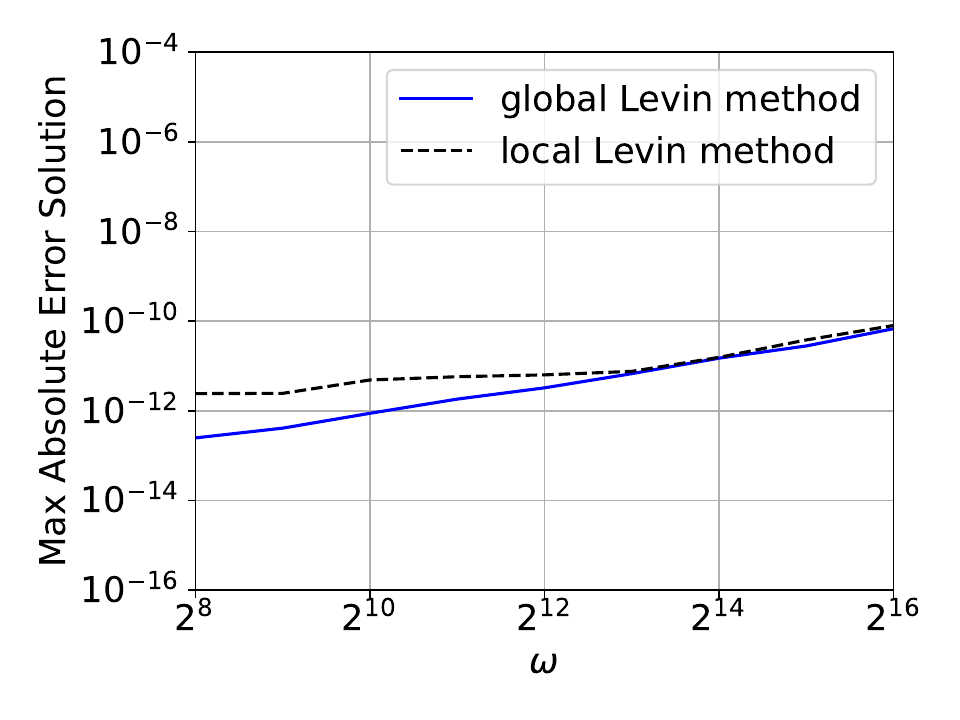}
\hfil
\includegraphics[width=.40\textwidth]{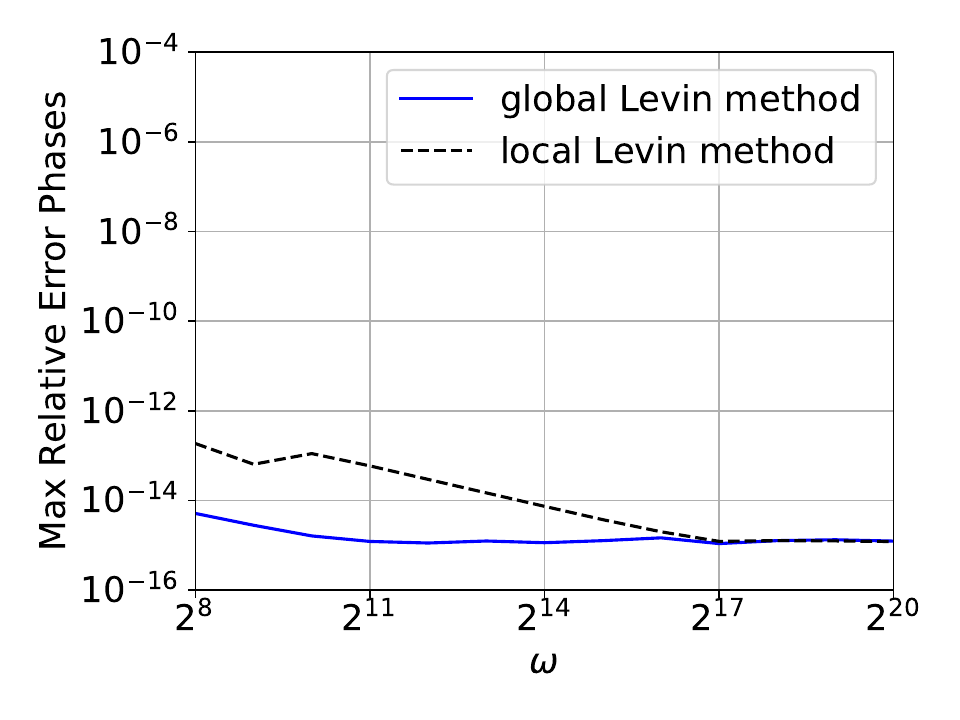}
\hfil

\caption{
The results of the experiments of Subsection~\ref{section:experiments:2}.
The  upper-left plot gives the time required by each of our methods
as a function of the parameter $\omega$.  The upper-right plot 
reports the number  of Chebyshev coefficients required to represent the slowly-varying phase functions
as a function of $\omega$.  The lower-left plot  reports the largest observed
absolute error in the solution of the initial value problem for (\ref{experiments2:1}), again as a function of $\omega$.
Finally, the plot on the lower right gives
the largest observed relative error in the slowly-varying phase functions constructed by our algorithms.
}
\label{experiments2:figure1}
\end{figure}

\begin{figure}[h]
\hfil
\includegraphics[width=.28\textwidth]{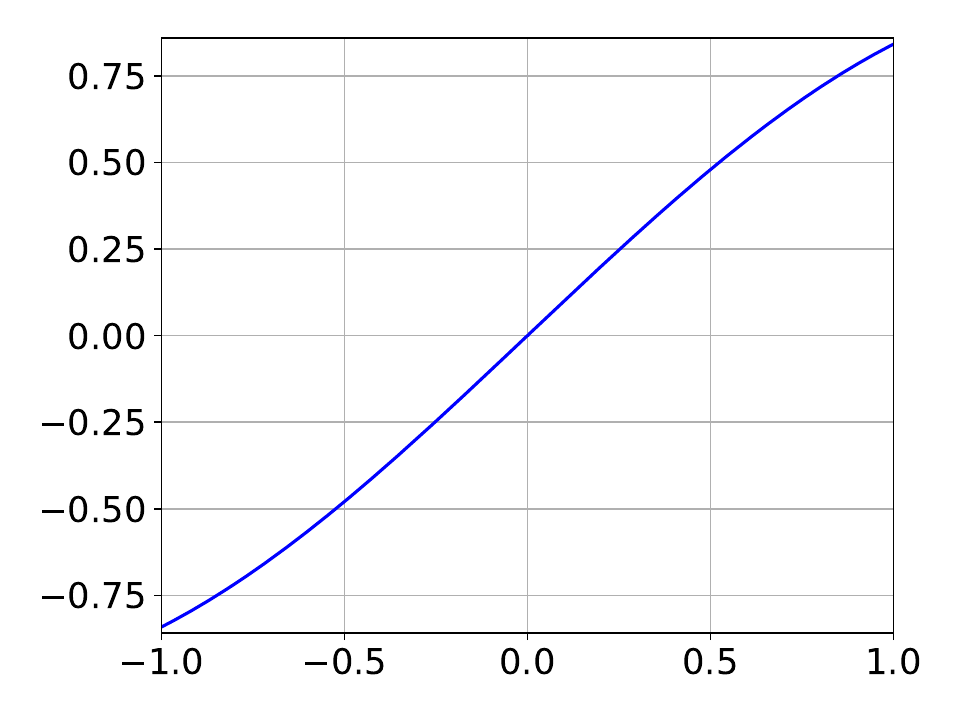}
\hfil
\includegraphics[width=.28\textwidth]{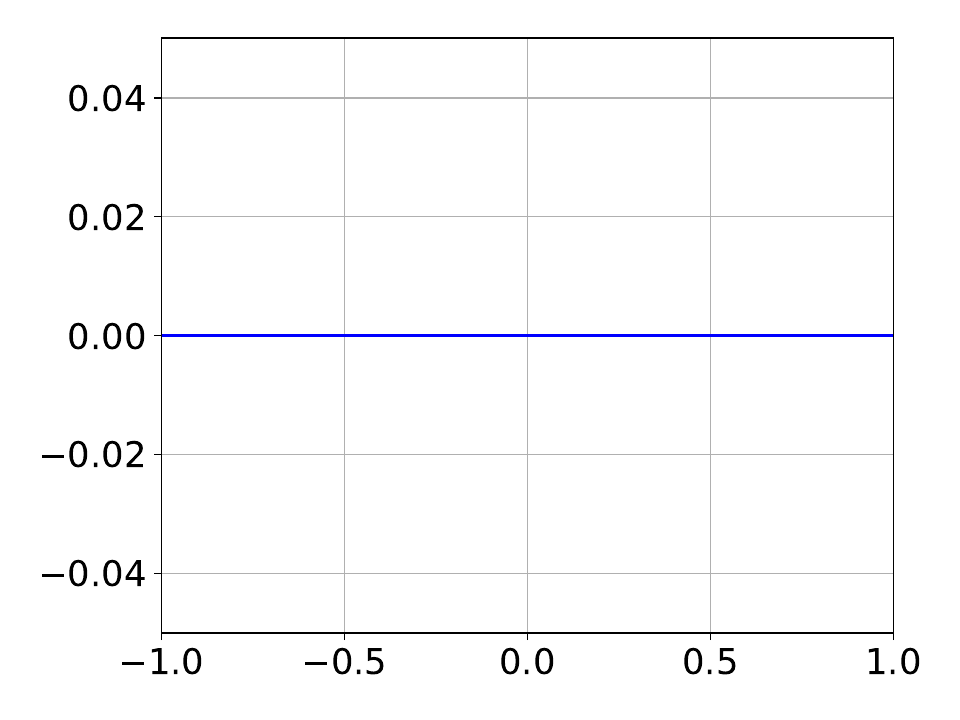}
\hfil
\includegraphics[width=.28\textwidth]{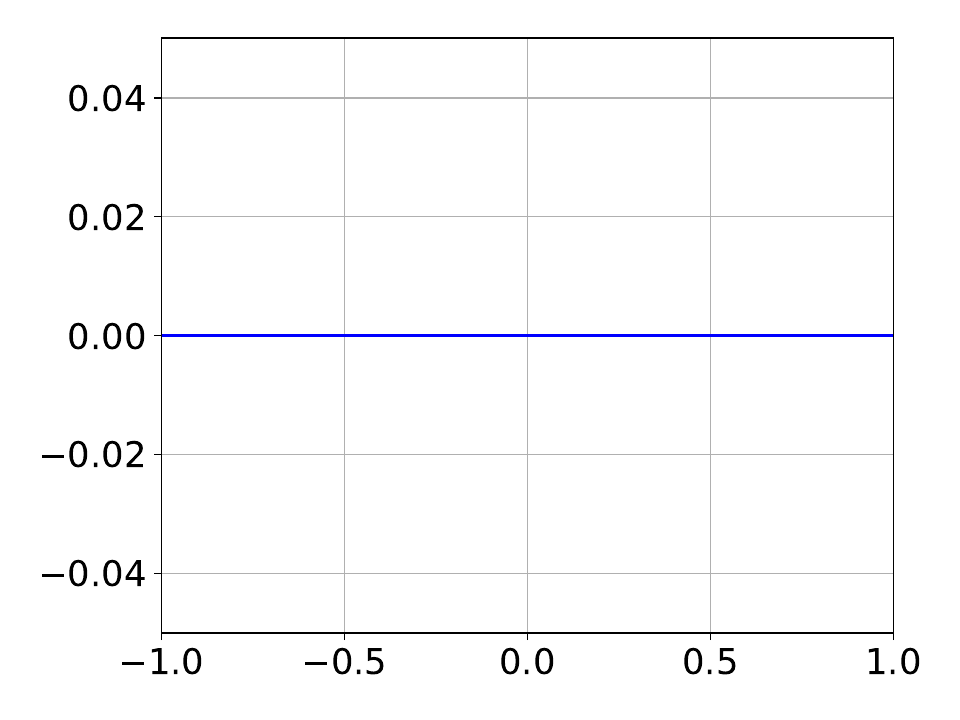}
\hfil

\hfil
\includegraphics[width=.28\textwidth]{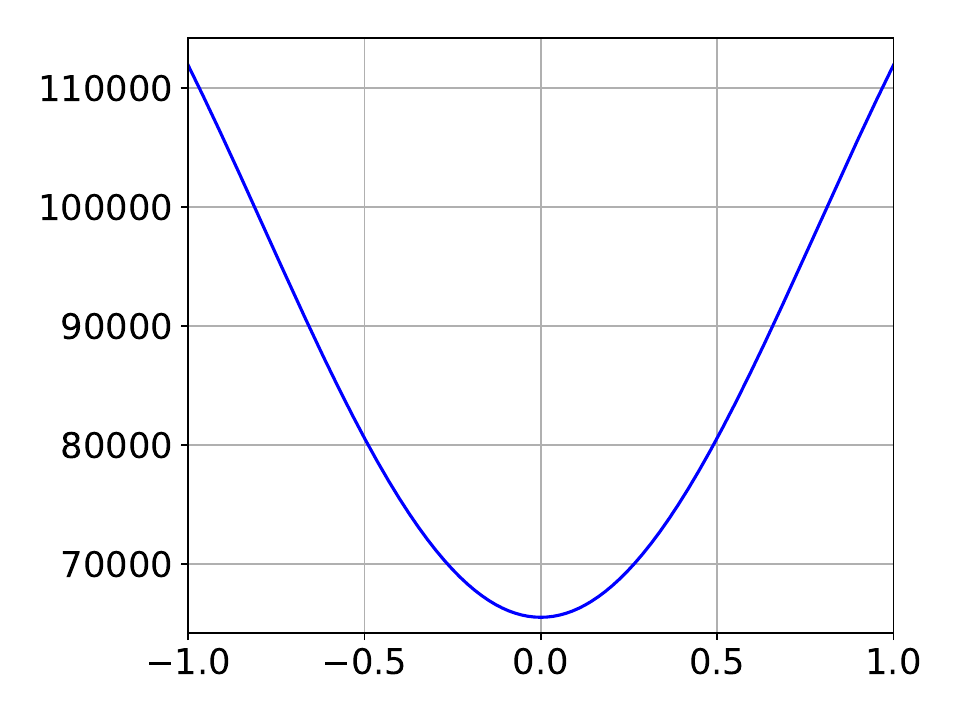}
\hfil
\includegraphics[width=.28\textwidth]{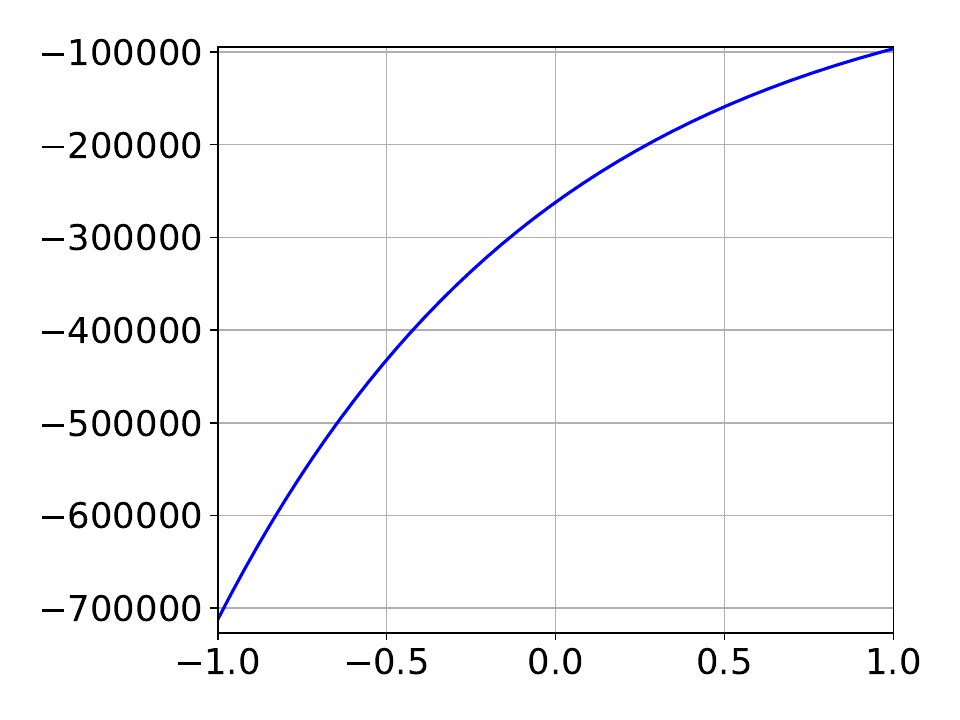}
\hfil
\includegraphics[width=.28\textwidth]{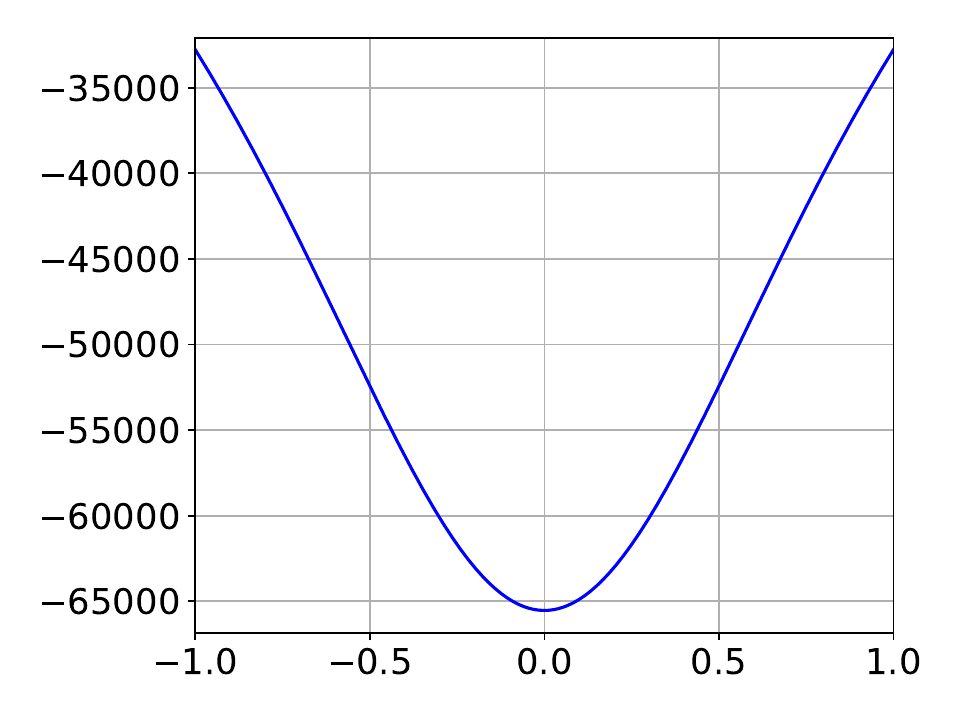}
\hfil

\caption{The eigenvalues $\lambda_1(t), \lambda_2(t), \lambda_3(t)$ of the coefficient 
matrix corresponding to Equation~(\ref{experiments2:1}) of Subsection~\ref{section:experiments:2}
when the parameter $\omega$ is equal to $2^{16}$.  Each column corresponds to one of the eigenvalues,
with the real part appearing in the first row and the imaginary part in the second.}
\label{experiments2:figure3}
\end{figure}

\begin{figure}[h]
\hfil
\includegraphics[width=.28\textwidth]{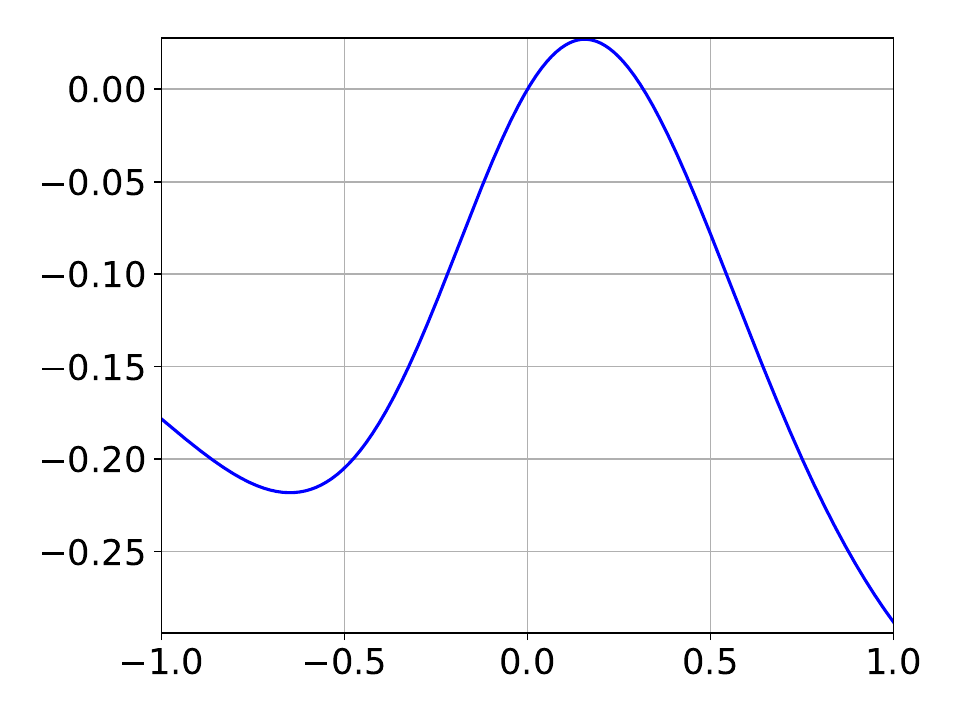}
\hfil
\includegraphics[width=.28\textwidth]{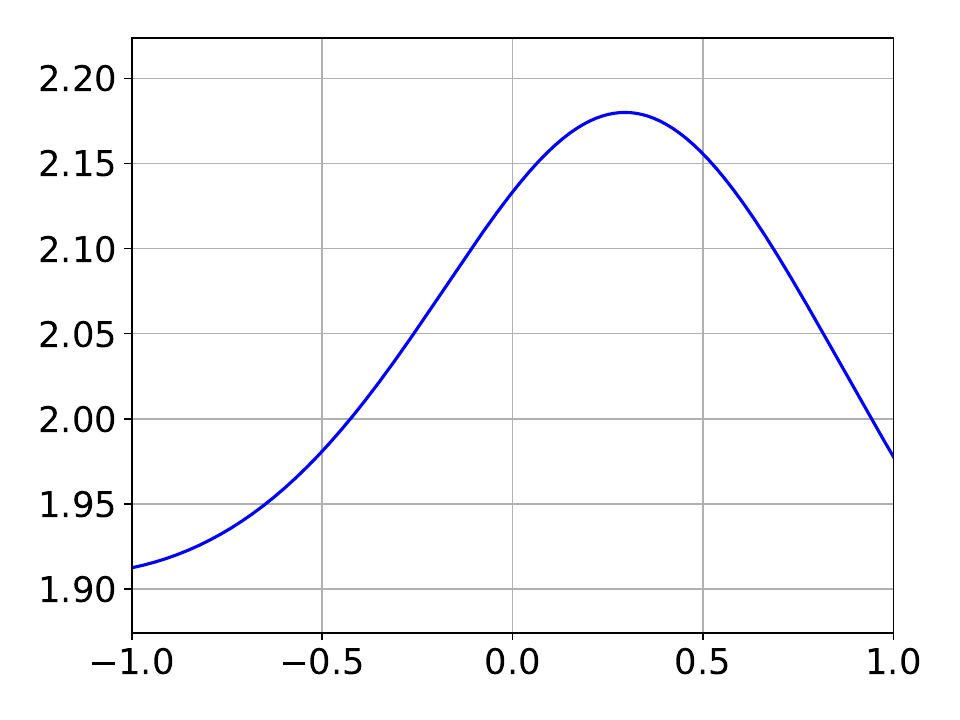}
\hfil
\includegraphics[width=.28\textwidth]{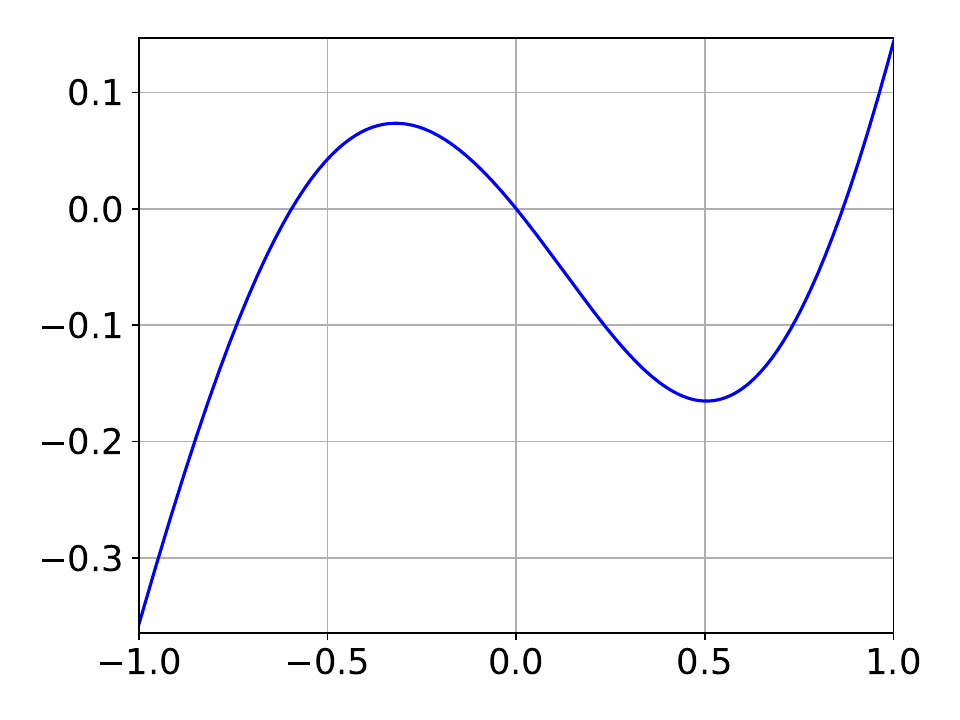}
\hfil

\hfil
\includegraphics[width=.28\textwidth]{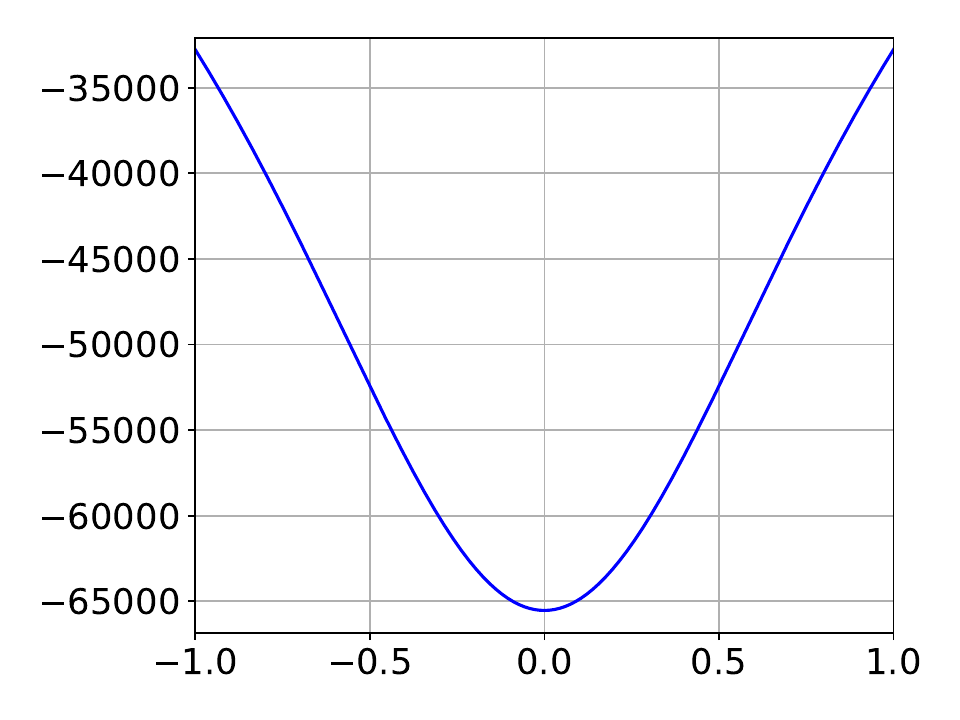}
\hfil
\includegraphics[width=.28\textwidth]{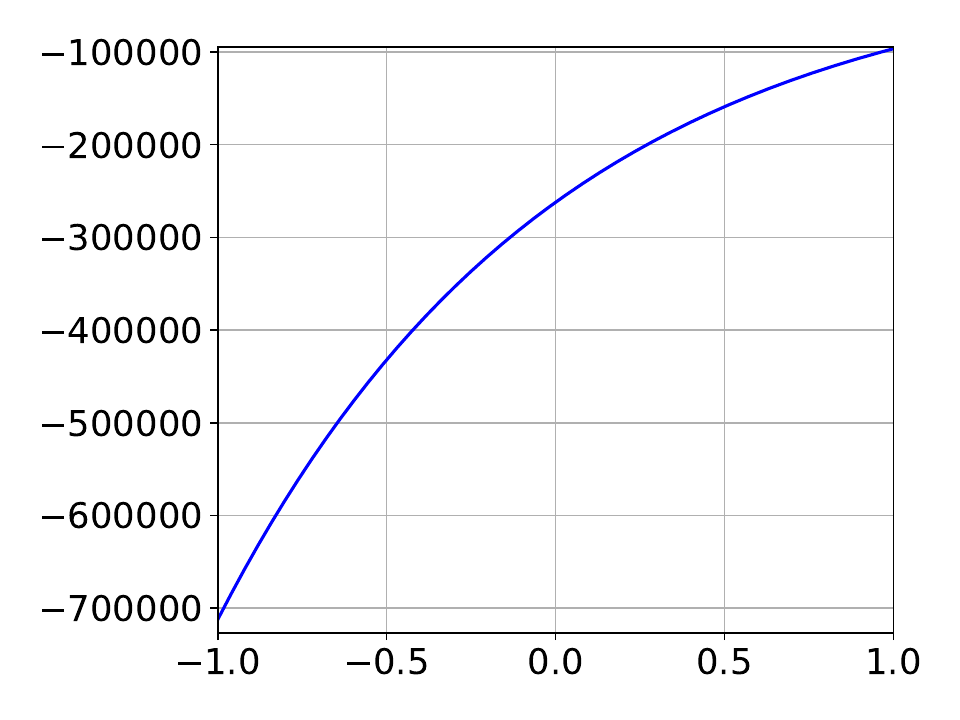}
\hfil
\includegraphics[width=.28\textwidth]{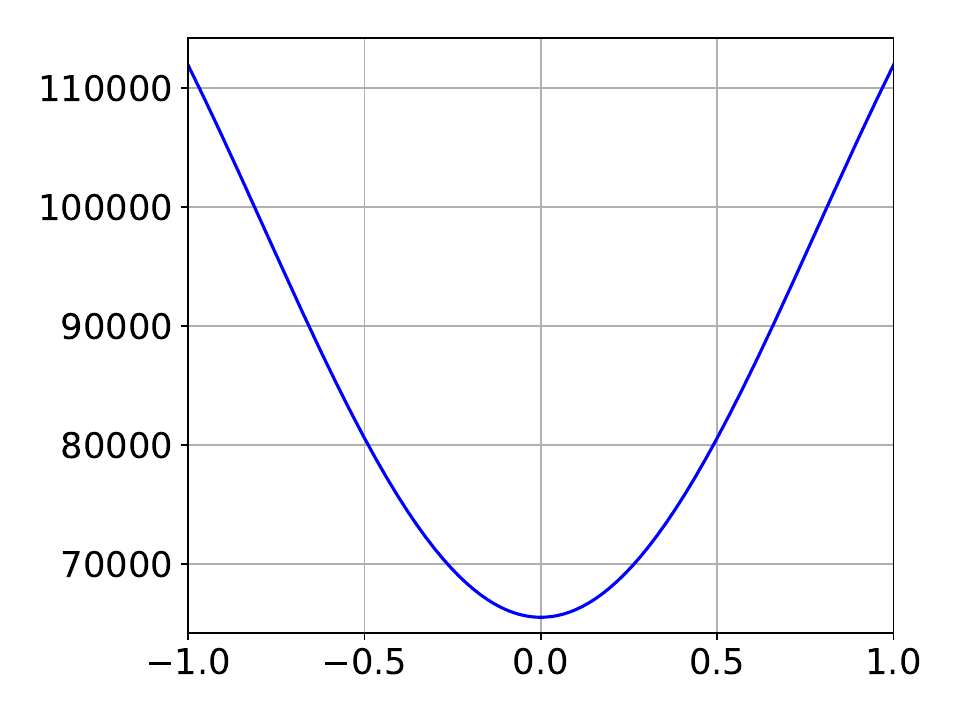}
\hfil

\caption{The derivatives of the three slowly-varying phase functions produced
by applying the global Levin method to Equation~(\ref{experiments2:1}) of Subsection~\ref{section:experiments:2}
when the parameter $\omega$ is equal to $2^{16}$.  Each column corresponds to one of the phase functions,
with the real part appearing in the first row and the imaginary part in the second.}
\label{experiments2:figure2}
\end{figure}

\end{subsection}

%
%
\begin{subsection}{A boundary value problem for a third order equation}
\label{section:experiments:3}

In the experiments described in this section, we considered the equation
\begin{equation}
y'''(t) - (1+2i\omega) (1+\sin^2(2t))\, y''(t) + (3i\omega+\omega^2)\frac{1}{1-\frac{t}{2}}\, y'(t) + 
(2\omega^2 -2i\omega^3) \frac{\exp(t)}{1+t^4}\, y(t) = 0.
\label{experiments3:1}
\end{equation}
The eigenvalues of the coefficient matrix for (\ref{experiments3:1}) at the point $0$ are
\begin{equation}
\begin{aligned}
\lambda_1(0) = i\omega,\ \ \ \lambda_2(0)=  2i\omega \ \ \ \mbox{and}\ \ \ \lambda_3(0)= 1-i\omega.
\end{aligned}
\end{equation}
Plots of the eigenvalues $\lambda_1(t),\lambda_2(t),\lambda_3(t)$ of the coefficient matrix
for (\ref{experiments3:1}) when $\omega=2^{16}$ can be found in 
Figure~\ref{experiments3:figure3} .

\begin{figure}[h]
\hfil
\includegraphics[width=.40\textwidth]{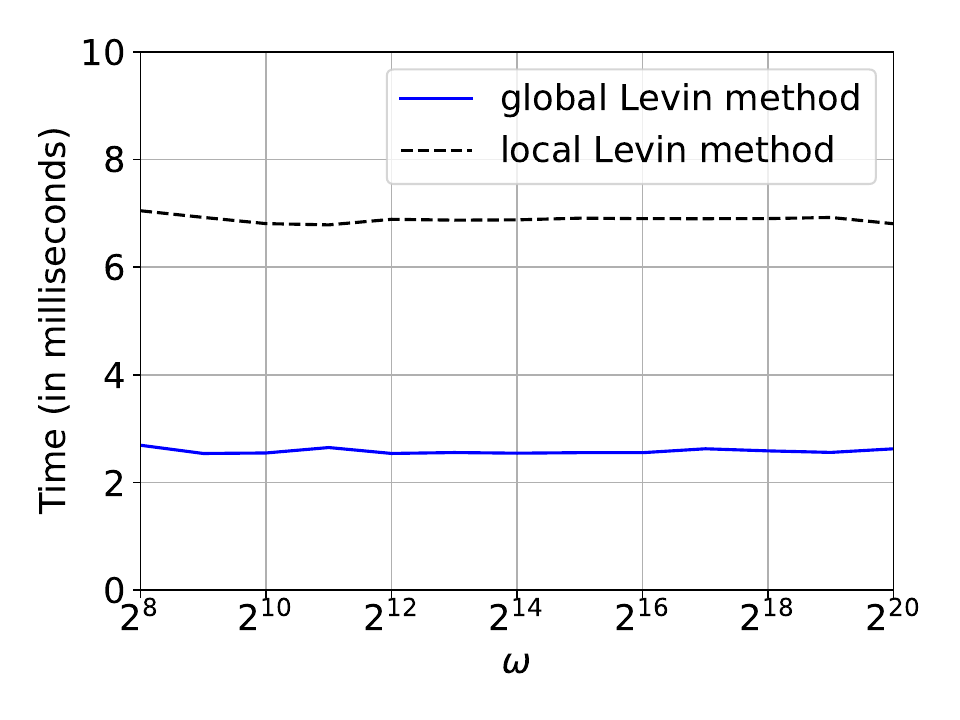}
\hfil
\includegraphics[width=.40\textwidth]{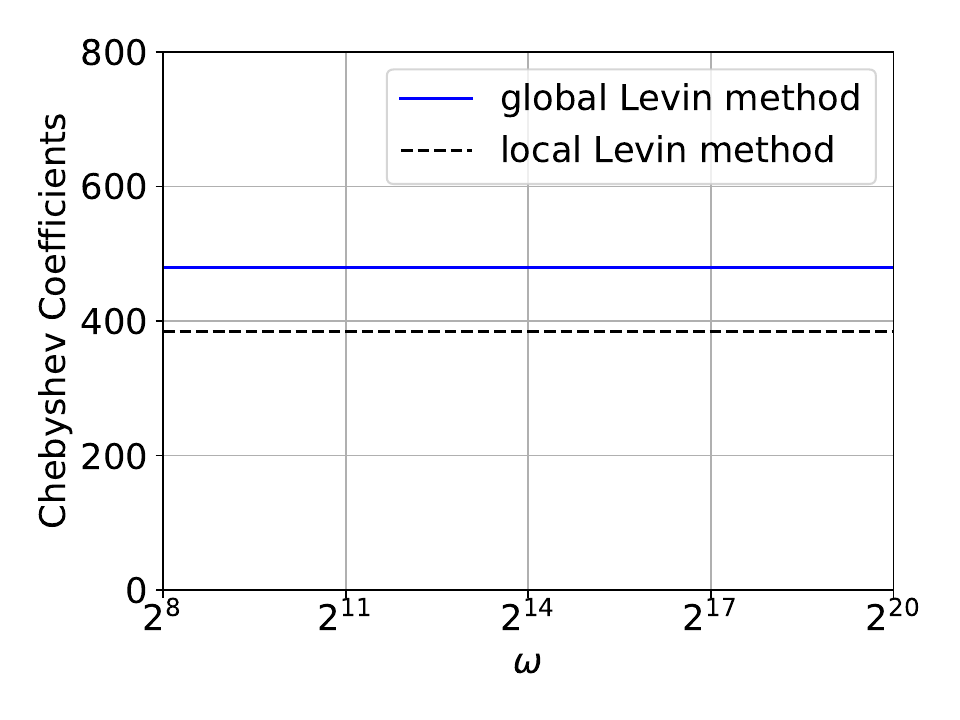}
\hfil

\hfil
\includegraphics[width=.40\textwidth]{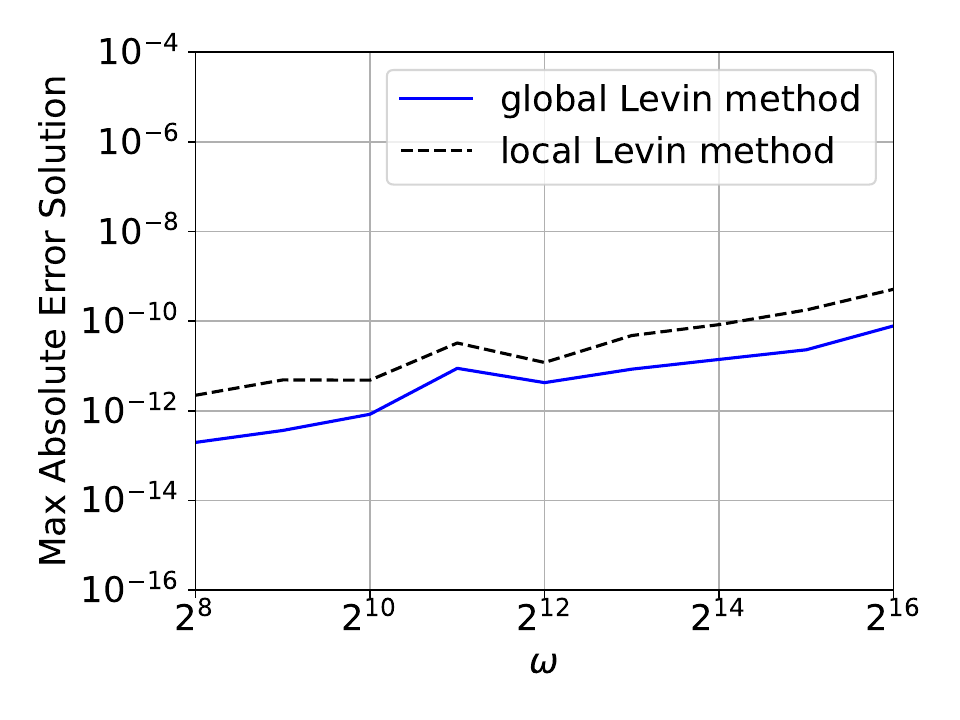}
\hfil
\includegraphics[width=.40\textwidth]{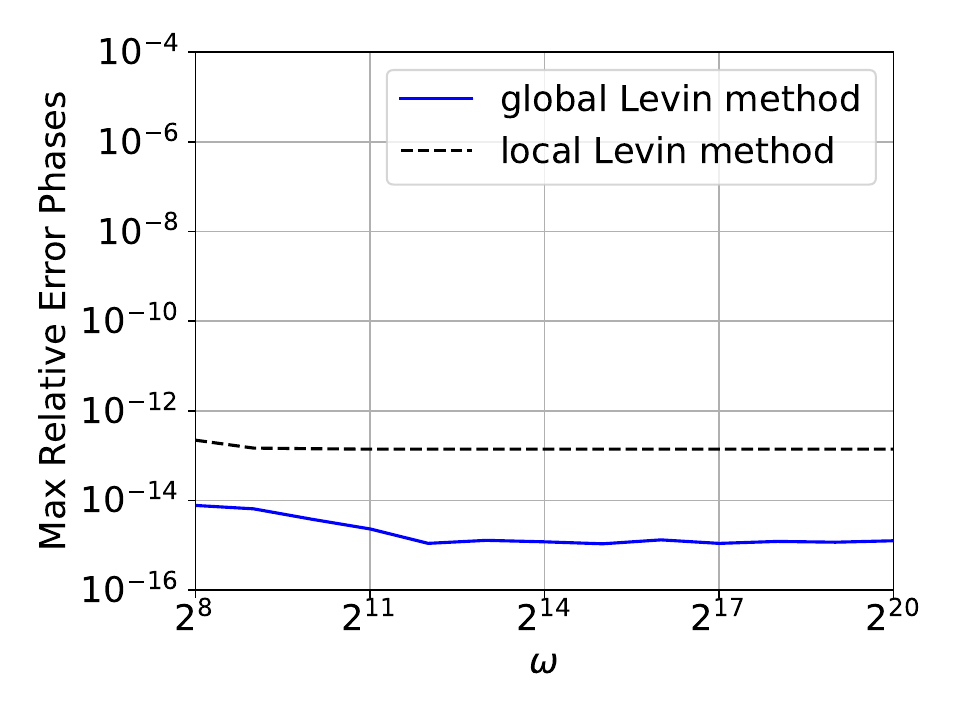}
\hfil

\caption{The results of the experiments of Subsection~\ref{section:experiments:3}.
The  upper-left plot gives the time required by each of our methods
as a function of the parameter $\omega$.  The upper-right plot 
reports the number  of Chebyshev coefficients required to represent the slowly-varying phase functions
as a function of $\omega$.  The lower-left plot  reports the largest observed
absolute error in the solution of a boundary value problem for (\ref{experiments3:1}), again as a function of $\omega$.
Finally, the plot on the lower right gives
the largest observed relative error in the slowly-varying phase functions constructed by our algorithms.
}
\label{experiments3:figure1}
\end{figure}

\begin{figure}[h]
\hfil
\includegraphics[width=.28\textwidth]{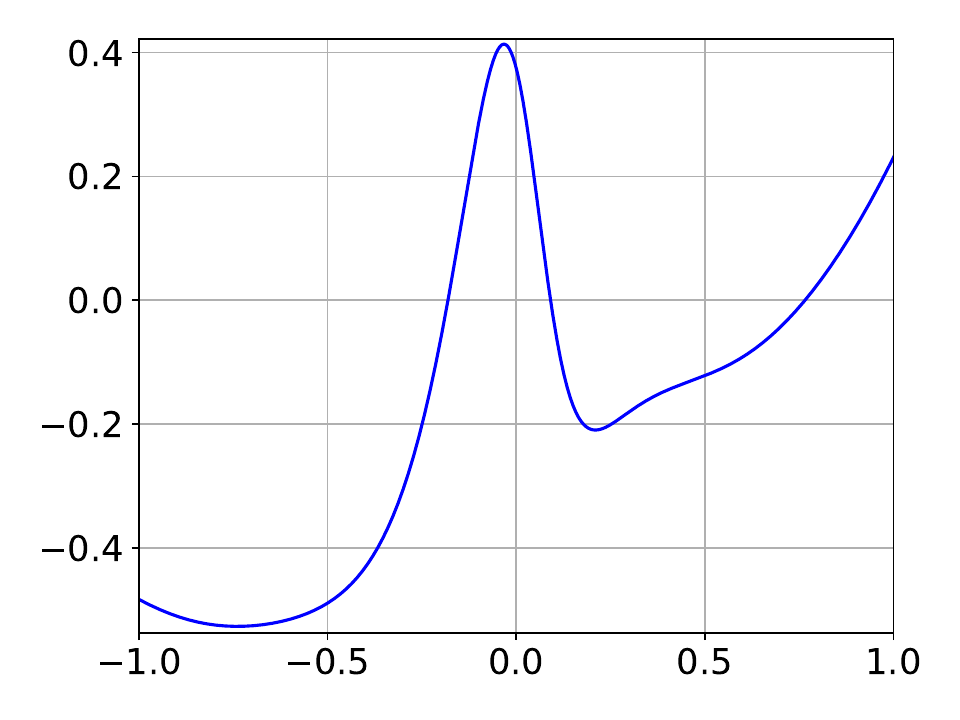}
\hfil
\includegraphics[width=.28\textwidth]{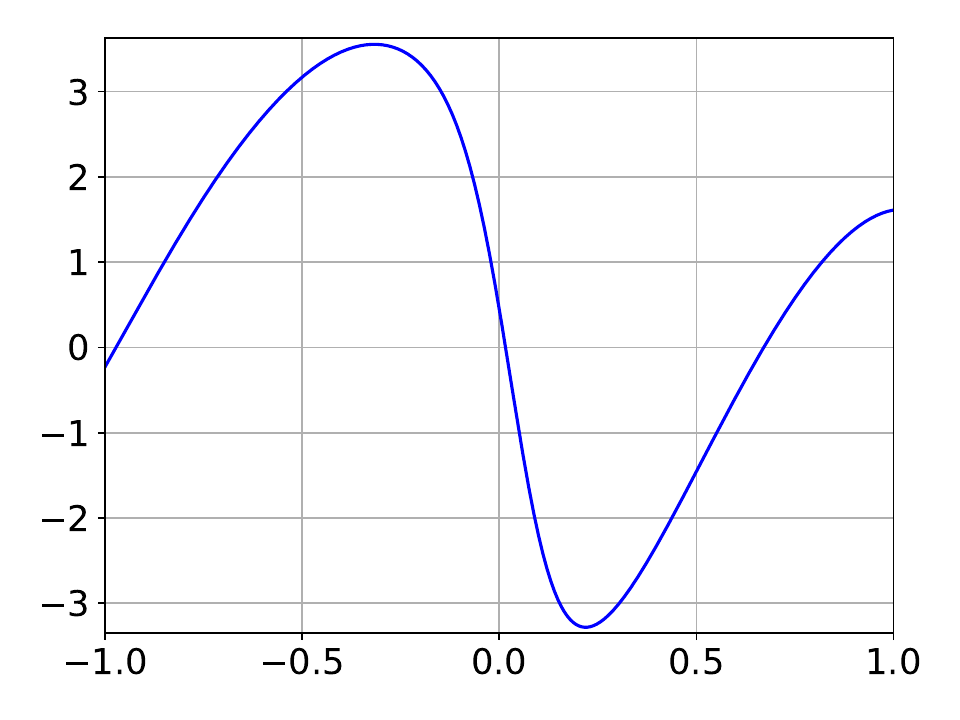}
\hfil
\includegraphics[width=.28\textwidth]{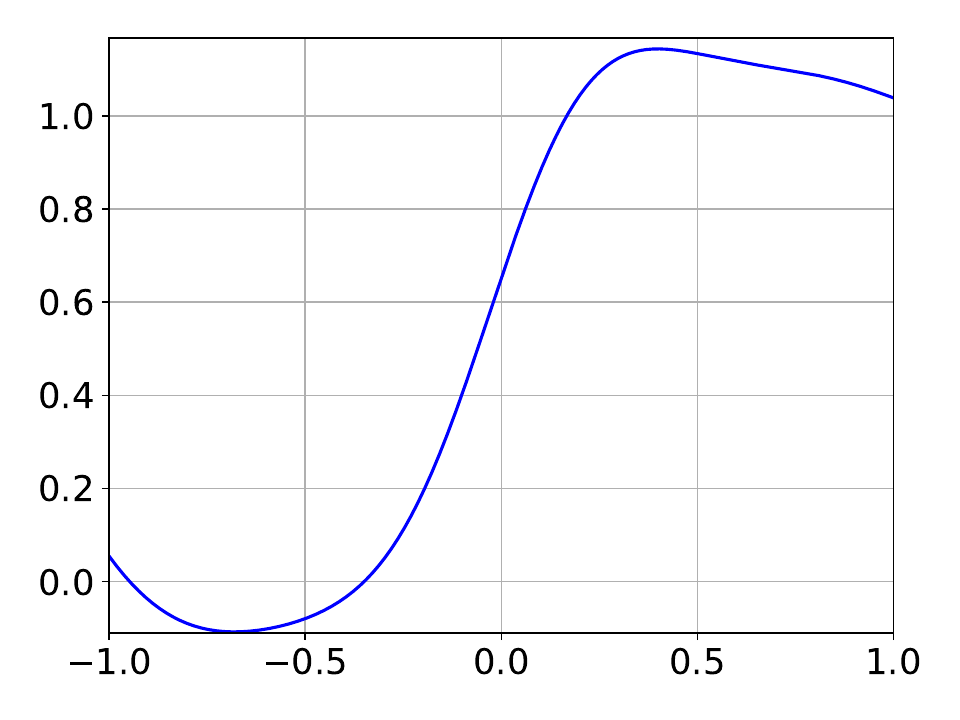}
\hfil

\hfil
\includegraphics[width=.28\textwidth]{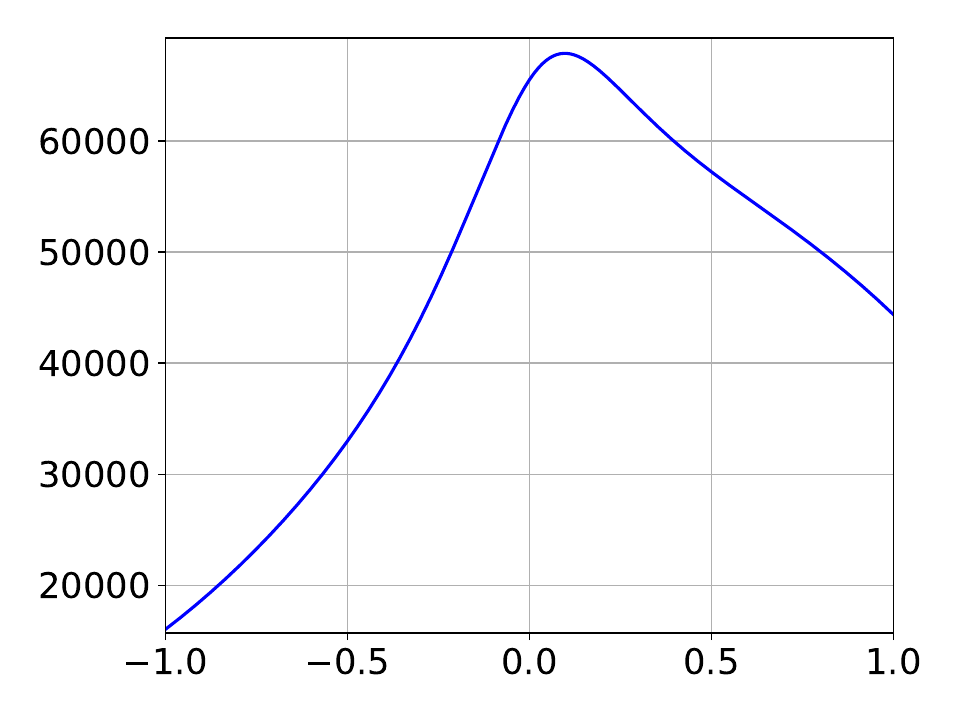}
\hfil
\includegraphics[width=.28\textwidth]{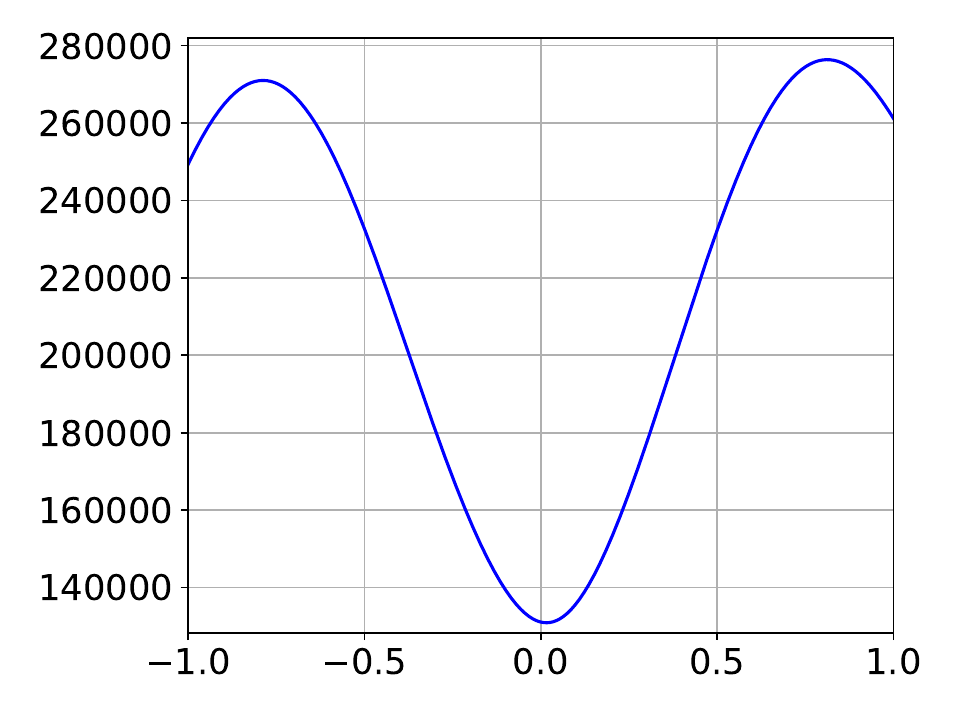}
\hfil
\includegraphics[width=.28\textwidth]{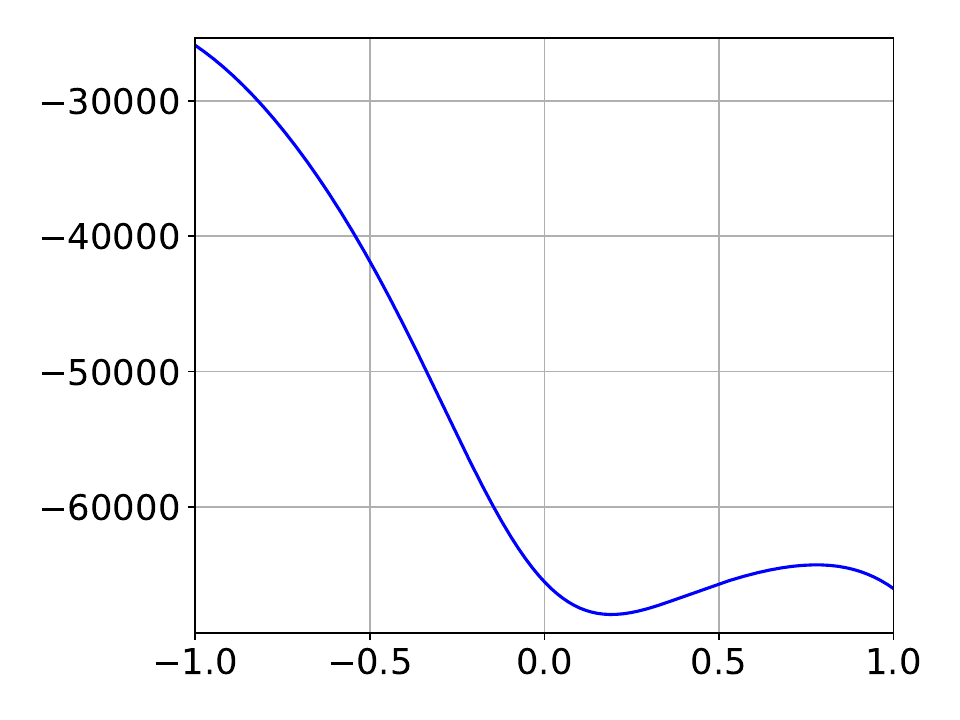}
\hfil

\caption{The derivatives of the three slowly-varying phase functions produced
by applying the global Levin method to Equation~(\ref{experiments3:1}) of Subsection~\ref{section:experiments:3}
when the parameter $\omega$ is equal to $2^{16}$.  Each column corresponds to one of the phase functions,
with the real part appearing in the first row and the imaginary part in the second.}
\label{experiments3:figure2}
\end{figure}

In the first experiment, for each $\omega=2^8,2^9,\ldots,2^{16}$, we 
we used the local Levin method and the global
Levin method to compute solutions of (\ref{experiments3:1})
over the interval $[-1,1]$ subject to the conditions
\begin{equation}
y(-1) =  y(1) = 1 \ \ \ \mbox{and} \ \ \ y'(-1) = 0.
\label{experiments3:2}
\end{equation}
We then
 measured the errors in each obtained solution at $10,000$ equispaced points in
the interval $[-1,1]$  by comparison with a reference solution constructed via the standard solver
described in Appendix~\ref{section:appendix}.
In the second experiment, for each $\omega=2^8,2^9,\ldots,2^{20}$, we constructed
slowly-varying phase functions for (\ref{experiments3:1}) over the interval $[-1,1]$ by running both the global
and local Levin methods using double precision arithmetic (as usual). 
We then measured the relative errors in each obtained phase function at $10,000$ equispaced points in
the interval $[-1,1]$  by comparison with phase functions constructed by applying the 
global Levin method to (\ref{experiments3:1}), but this time using extended precision arithmetic to perform the calculations.
Figure~\ref{experiments3:figure1} gives the results of these two experiments.
Figure~\ref{experiments3:figure2}
contains  plots of the derivatives of the slowly-varying phase functions produced by global Levin method
when $\omega=2^{16}$.

\end{subsection}

\begin{figure}[h]
\hfil
\includegraphics[width=.28\textwidth]{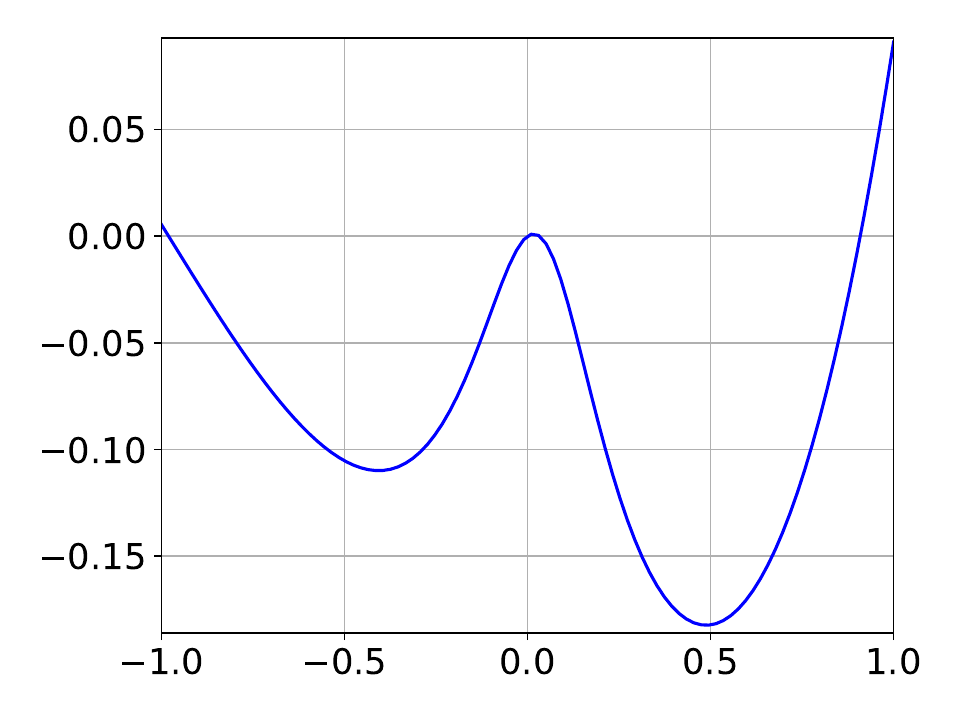}
\hfil
\includegraphics[width=.28\textwidth]{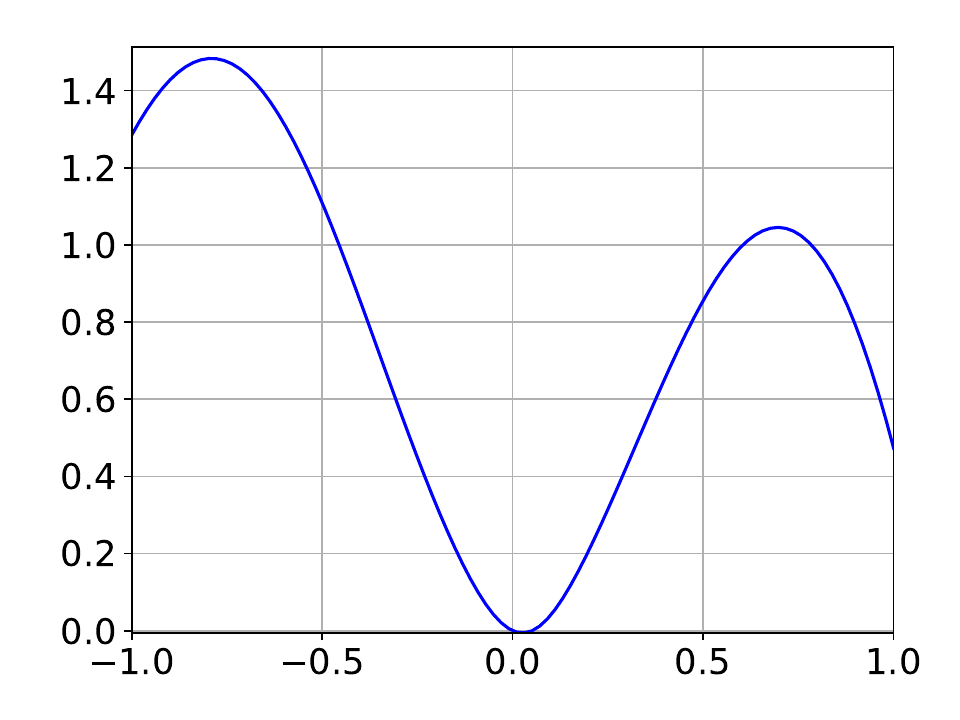}
\hfil
\includegraphics[width=.28\textwidth]{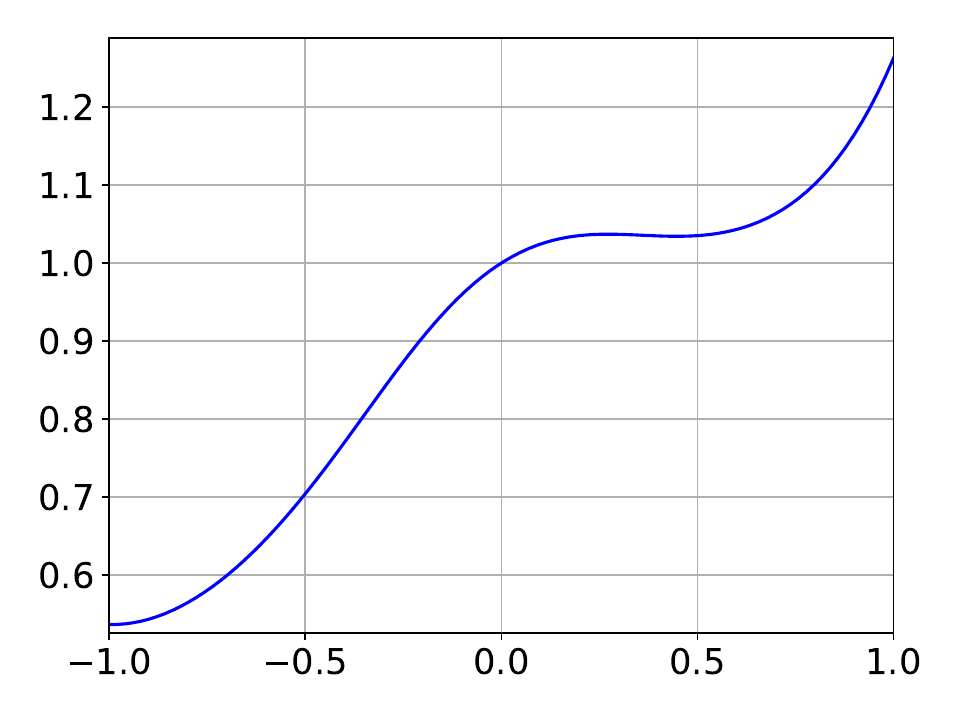}
\hfil

\hfil
\includegraphics[width=.28\textwidth]{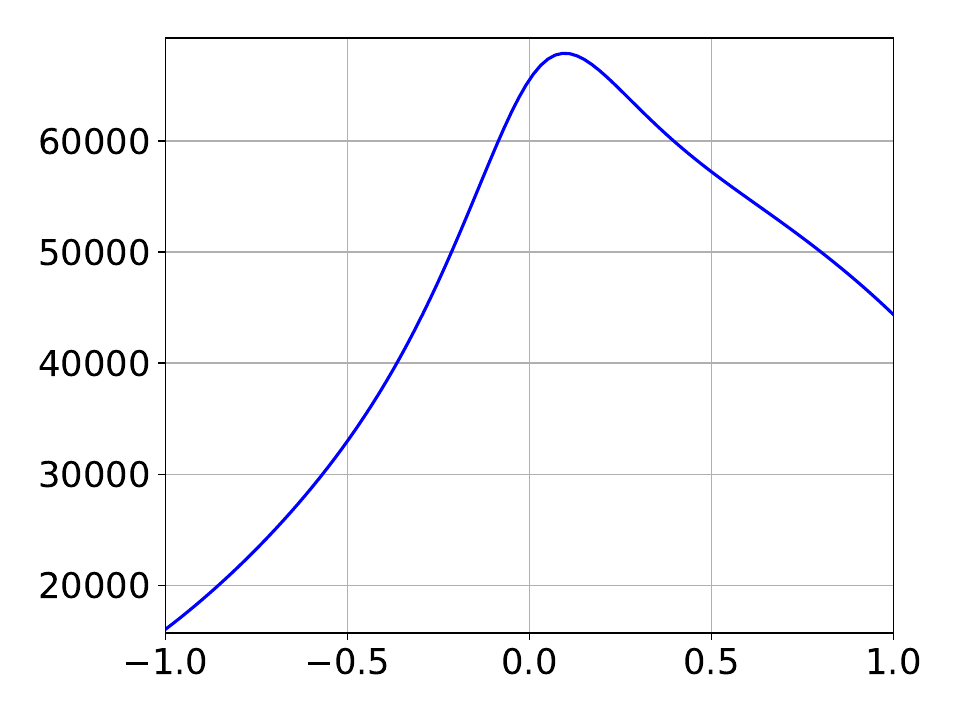}
\hfil
\includegraphics[width=.28\textwidth]{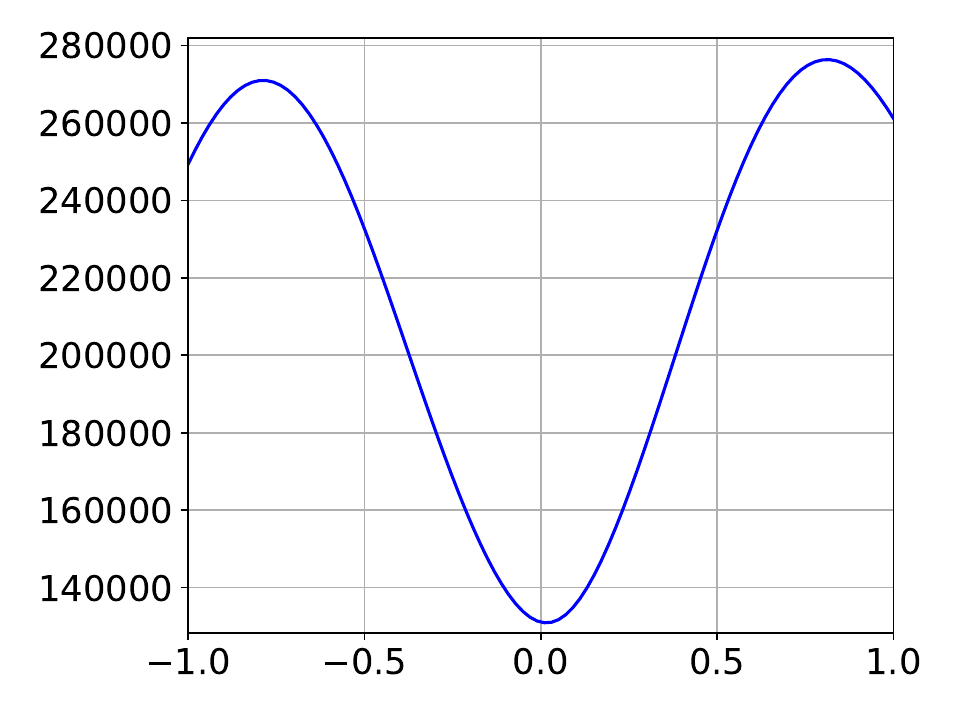}
\hfil
\includegraphics[width=.28\textwidth]{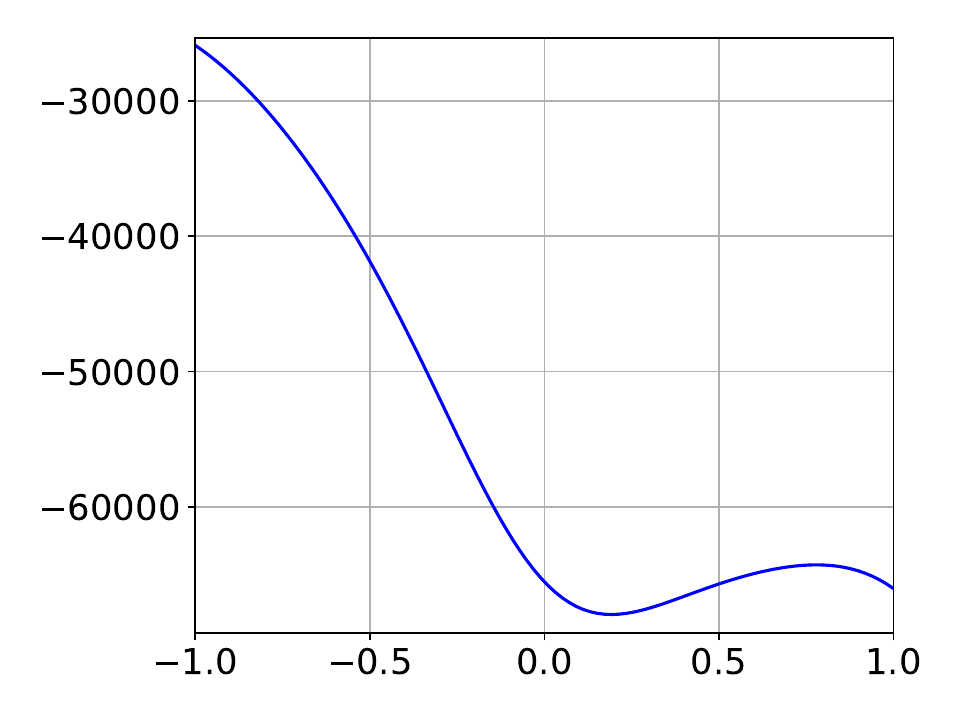}
\hfil

\caption{The eigenvalues $\lambda_1(t), \lambda_2(t), \lambda_3(t)$ of the coefficient 
matrix corresponding to Equation~(\ref{experiments3:1}) of Subsection~\ref{section:experiments:3}
when the parameter $\omega$ is equal to $2^{16}$.  Each column corresponds to one of the eigenvalues,
with the real part appearing in the first row and the imaginary part in the second.}
\label{experiments3:figure3}
\end{figure}

%
%
\begin{subsection}{An initial value problem for a fourth order equation}
\label{section:experiments:4}

In the experiments described in this section, we considered the equation
\begin{equation}
y''''(t) + \left(-5i \omega (1+t^2) + 5\omega^2\, \frac{8+\cos^4(3t)}{2+t^4} \right) y''(t) + 4 \omega^4\, \frac{2+\sin(3t)}{2+t} y(t) = 0.
\label{experiments4:1}
\end{equation}
We have the following formulas for the eigenvalues of the coefficient matrix at the point $0$:
\begin{equation*}
\begin{aligned}
\lambda_1(0) &= \sqrt{\frac{5i\omega-5\omega^2+\omega\sqrt{-25-50i\omega+9\omega^2}}{2}} \sim -i\omega + \frac{5}{6} + \mathcal{O}\left(\frac{1}{\omega}\right)\\
\lambda_2(0) &= -\sqrt{\frac{5i\omega-5\omega^2+\omega\sqrt{-25-50i\omega+9\omega^2}}{2}} \sim i\omega - \frac{5}{6} + \mathcal{O}\left(\frac{1}{\omega}\right)\\
\lambda_3(0) &= \sqrt{\frac{5i\omega-5\omega^2-\omega\sqrt{-25-50i\omega+9\omega^2}}{2}} \sim 2i\omega + \frac{5}{3} + \mathcal{O}\left(\frac{1}{\omega}\right)\\
\lambda_4(0) &= -\sqrt{\frac{5i\omega-5\omega^2-\omega\sqrt{-25-50i\omega+9\omega^2}}{2}} \sim -2i\omega - \frac{5}{3} + \mathcal{O}\left(\frac{1}{\omega}\right).
\end{aligned}
\end{equation*}
%

In the first experiment, for $\omega=2^8,2^9,\ldots,2^{16}$, we used both the global Levin method
and the local Levin method to solve (\ref{experiments4:1}) over the interval
$[-1,1]$ subject to the conditions
\begin{equation}
y(0) = 1, \ \ \ y'(0) = i\omega, \ \ \ y''(0) = -\omega^2 \ \ \ \mbox{and}\ \ \ y'''(0) = -i \omega^3.
\label{experiments4:2}
\end{equation}
Then we measured the errors in each obtained solution at $10,000$ equispaced points in
the interval $[-1,1]$  by comparison with reference solutions constructed via the standard solver
described in Appendix~\ref{section:appendix}.
In the second experiment, for each $\omega=2^8,2^9,\ldots,2^{20}$, we constructed
slowly-varying phase functions for (\ref{experiments4:1}) over the interval $[-1,1]$ by running both the global
and local Levin methods using double precision arithmetic (as usual). 
We then measured the relative errors in each obtained phase function at $10,000$ equispaced points in
the interval $[-1,1]$  by comparison with phase functions constructed by applying the 
global Levin method to (\ref{experiments4:1}), but this time using extended precision arithmetic to perform the calculations.
Figure~\ref{experiments4:figure1} gives the results of these two experiments.
Figure~\ref{experiments4:figure2}
contains  plots of the derivatives of the slowly-varying phase functions produced by global Levin method
when $\omega=2^{16}$.

\begin{figure}[h]
\hfil
\includegraphics[width=.40\textwidth]{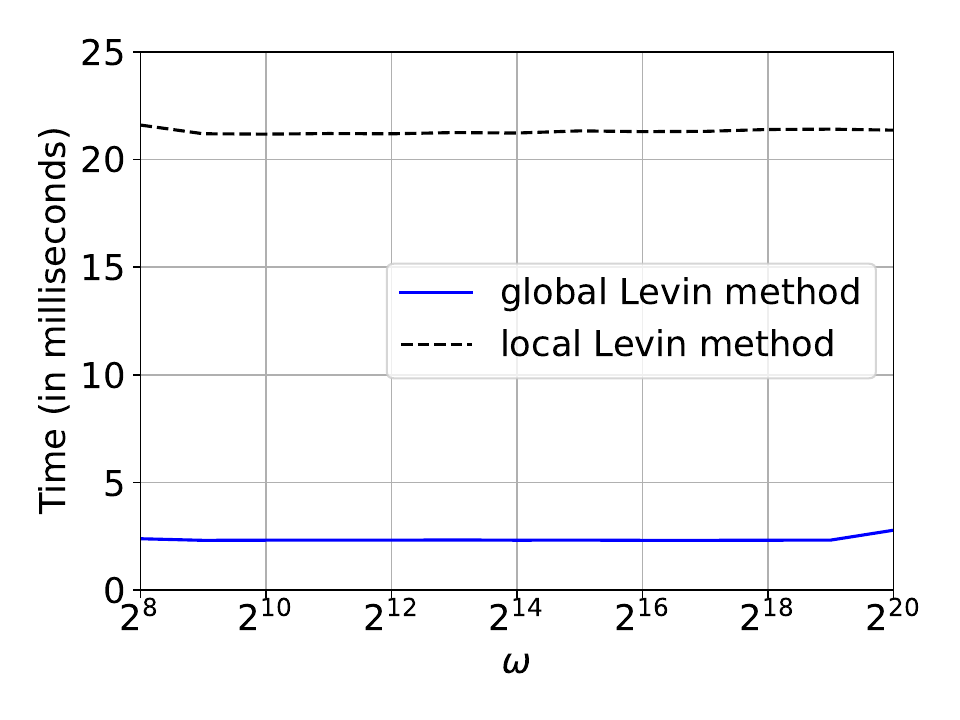}
\hfil
\includegraphics[width=.40\textwidth]{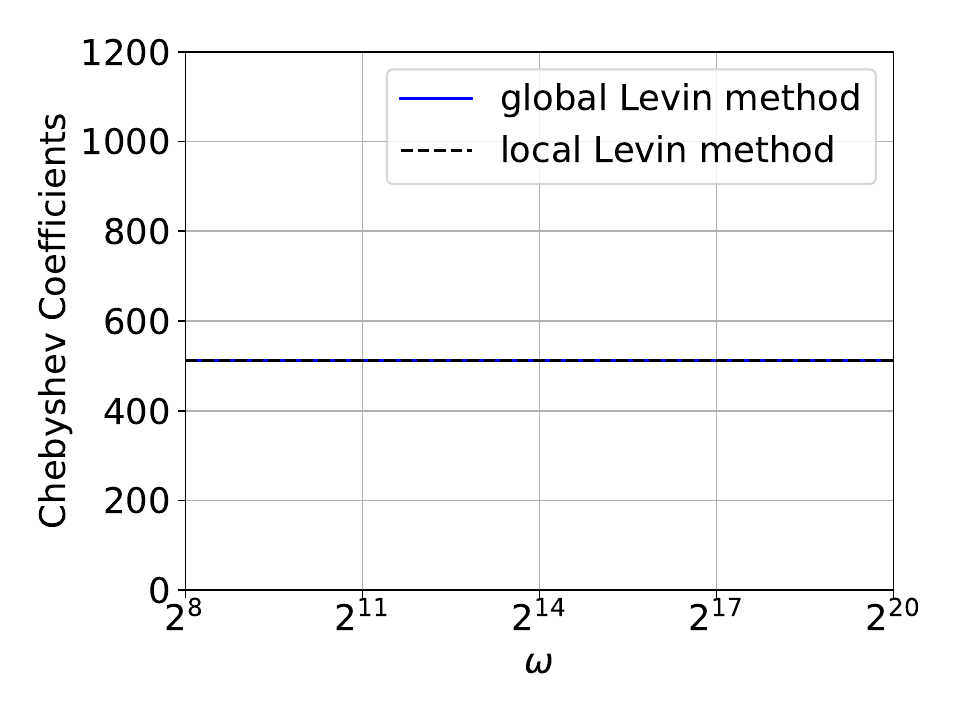}
\hfil

\hfil
\includegraphics[width=.40\textwidth]{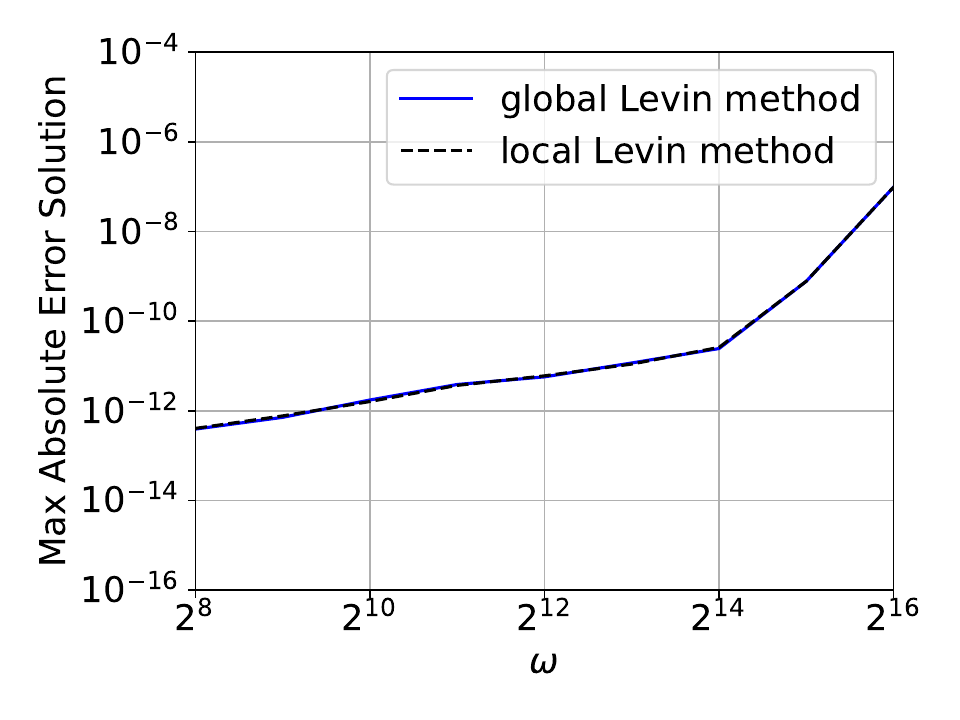}
\hfil
\includegraphics[width=.40\textwidth]{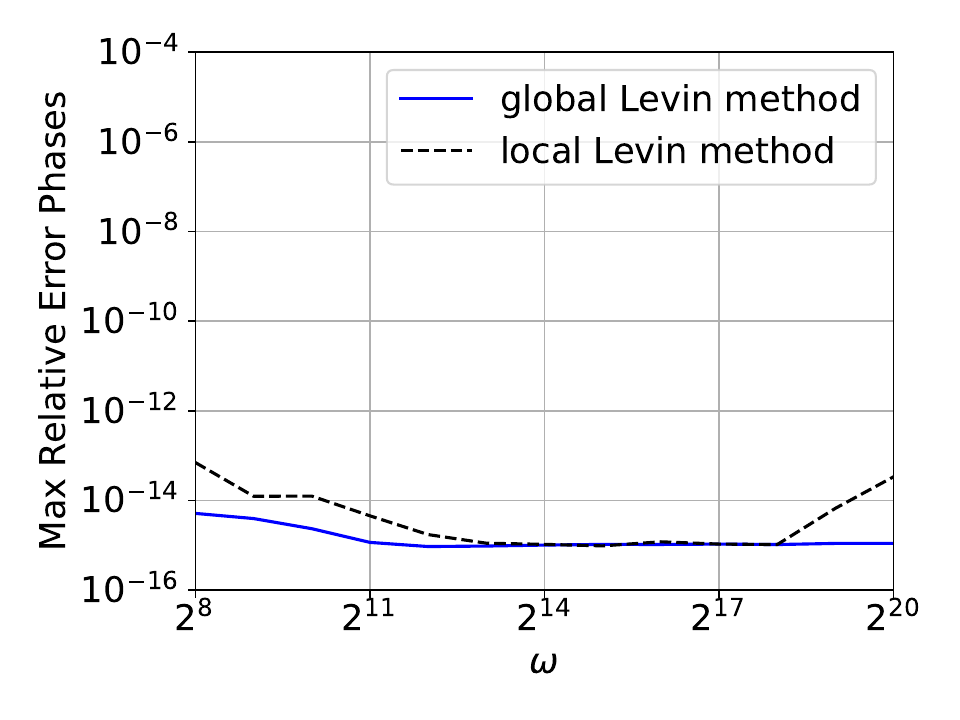}
\hfil

\caption{The results of the experiment of Subsection~\ref{section:experiments:4}.
The  upper-left plot gives the time required by each of our methods
as a function of the parameter $\omega$.  The upper-right plot 
reports the number  of Chebyshev coefficients required to represent the slowly-varying phase functions
as a function of $\omega$.  The lower-left plot  reports the largest observed
absolute error in the solution of the initial value problem for (\ref{experiments4:1}), again as a function of $\omega$.
The plot on the lower right gives
the largest observed relative error in the slowly-varying phase functions constructed by our algorithms.
The plot on the lower left appears to have only one line because the solutions
generated by the local and global Levin methods closely coincide.
}
\label{experiments4:figure1}
\end{figure}

\begin{figure}[h]
\hfil
\includegraphics[width=.22\textwidth]{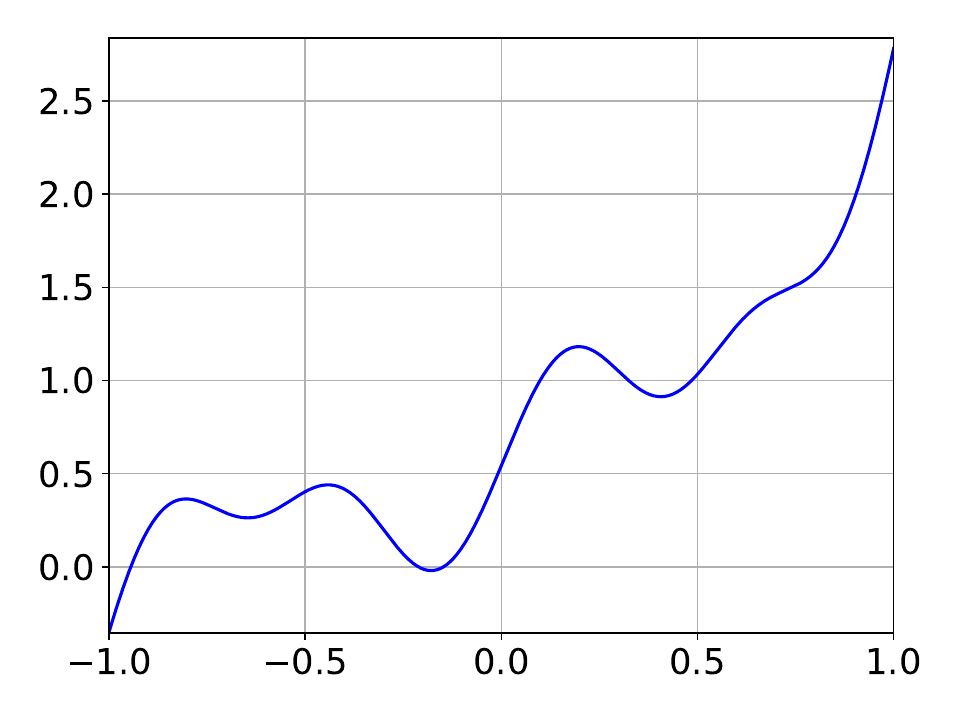}
\hfil
\includegraphics[width=.22\textwidth]{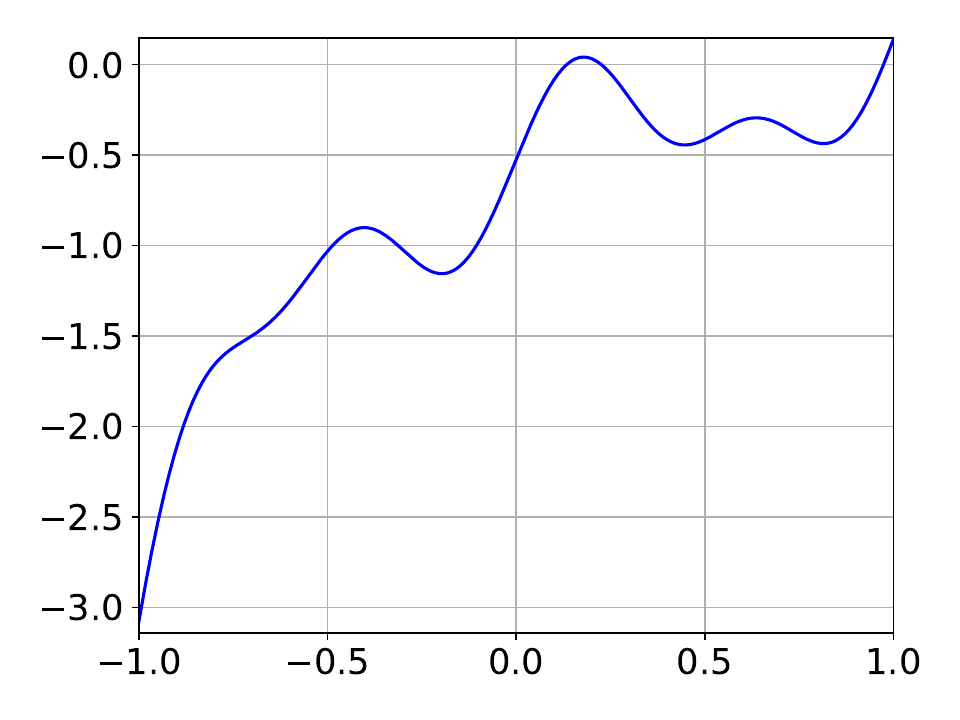}
\hfil
\includegraphics[width=.22\textwidth]{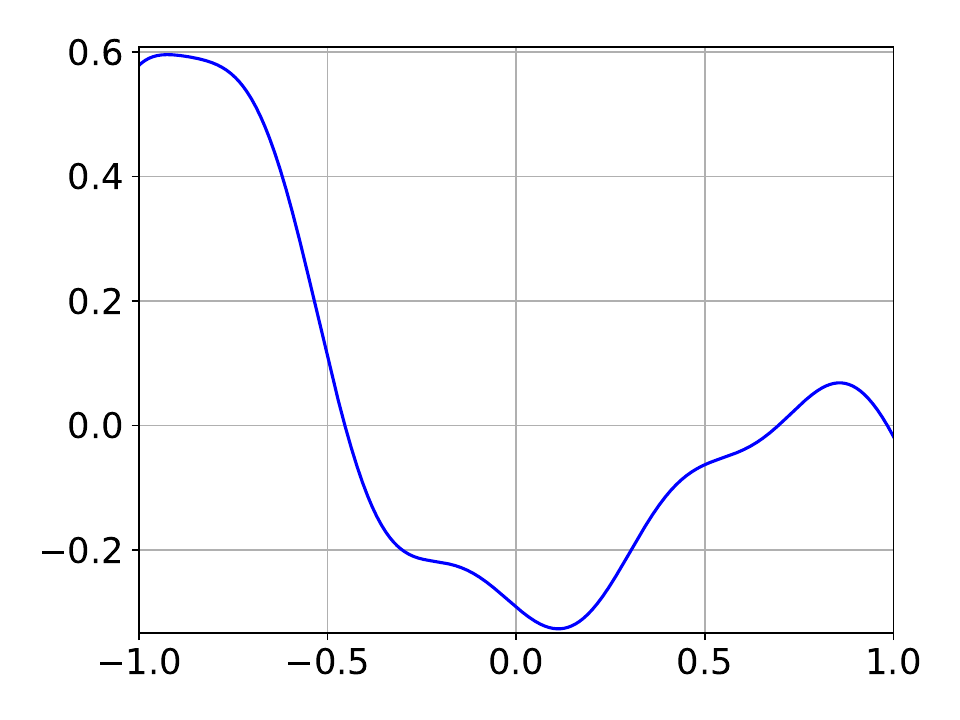}
\hfil
\includegraphics[width=.22\textwidth]{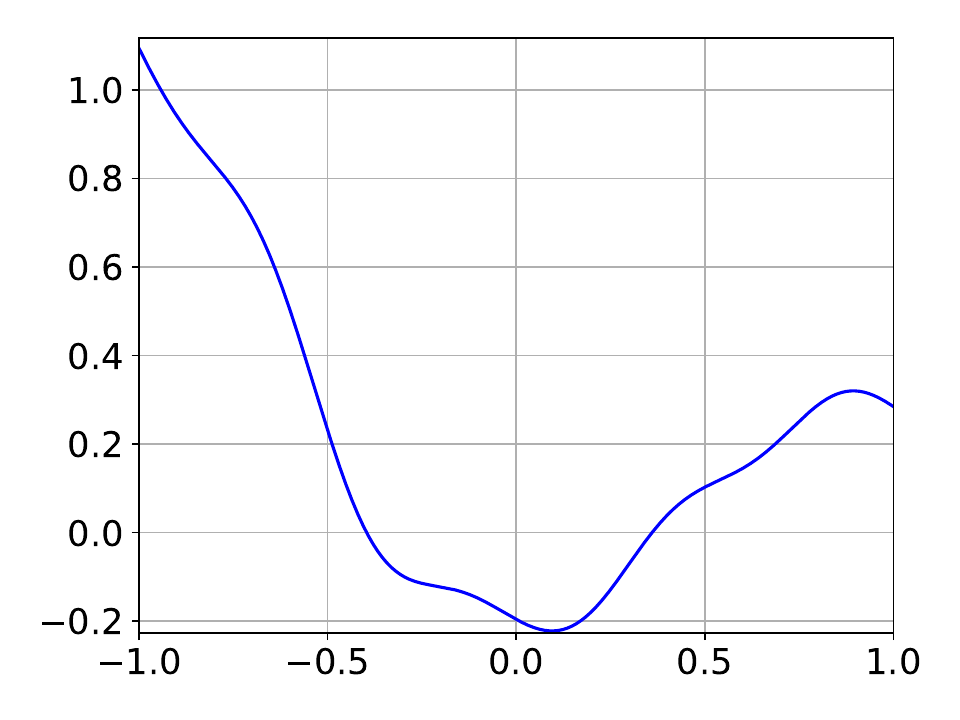}
\hfil

\hfil
\includegraphics[width=.22\textwidth]{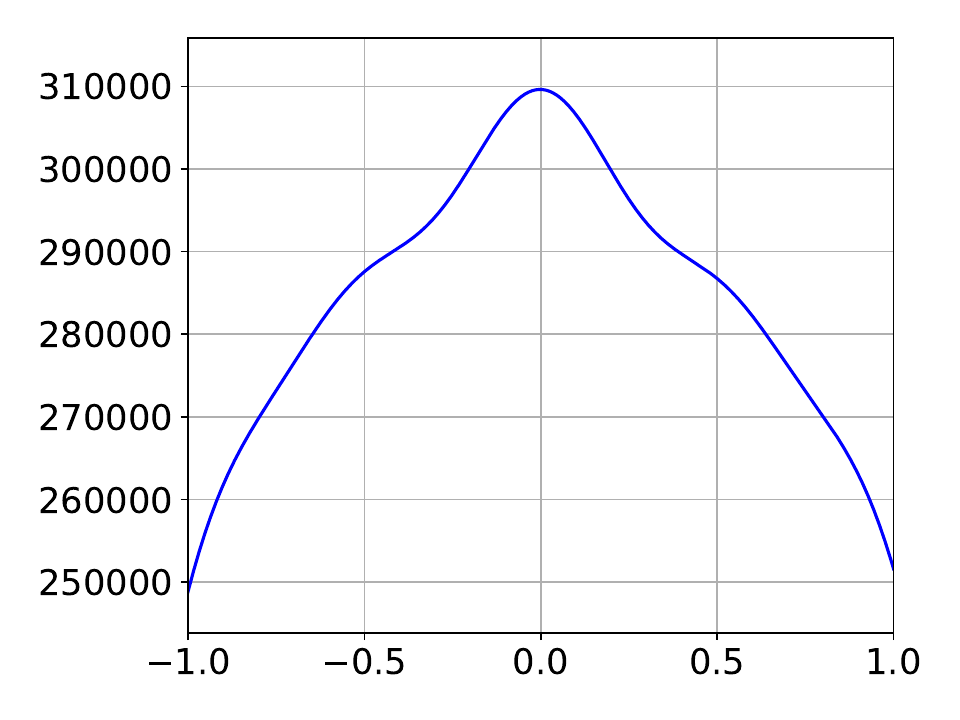}
\hfil
\includegraphics[width=.22\textwidth]{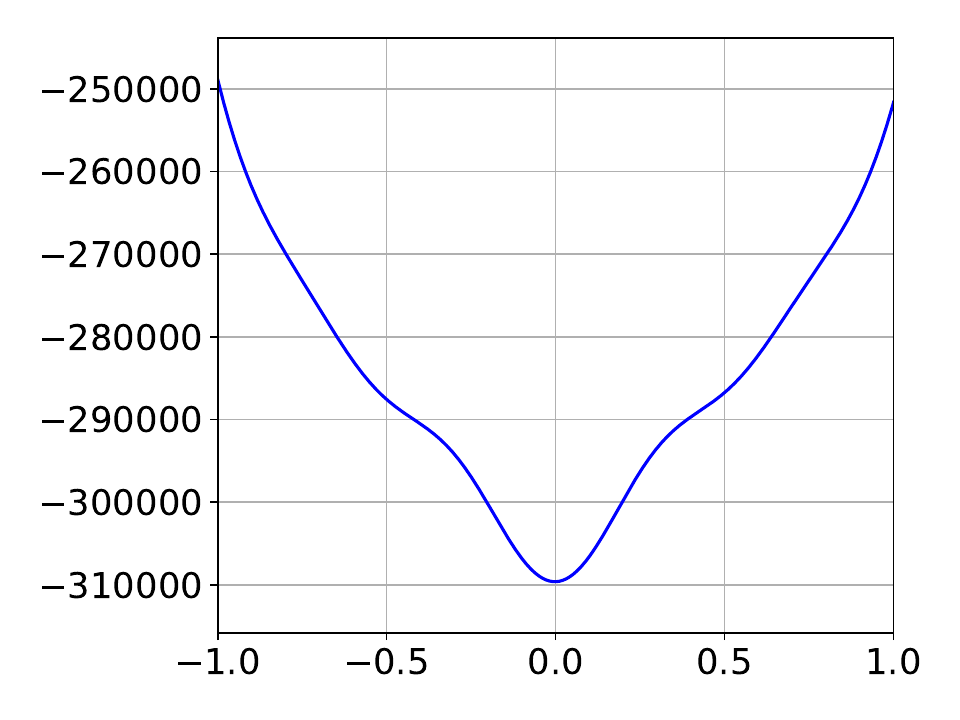}
\hfil
\includegraphics[width=.22\textwidth]{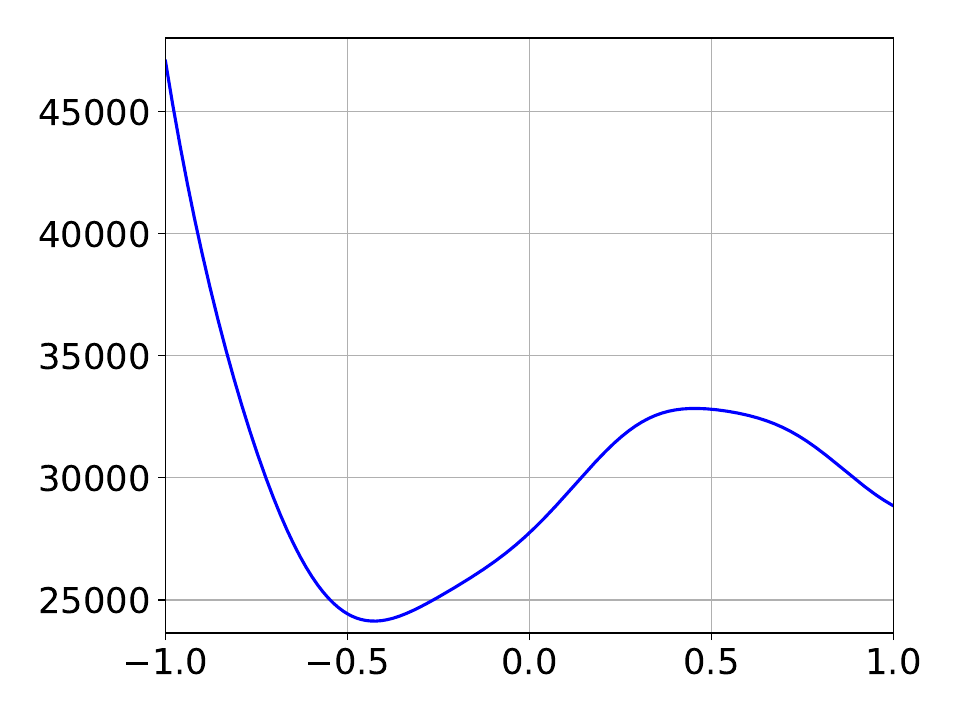}
\hfil
\includegraphics[width=.22\textwidth]{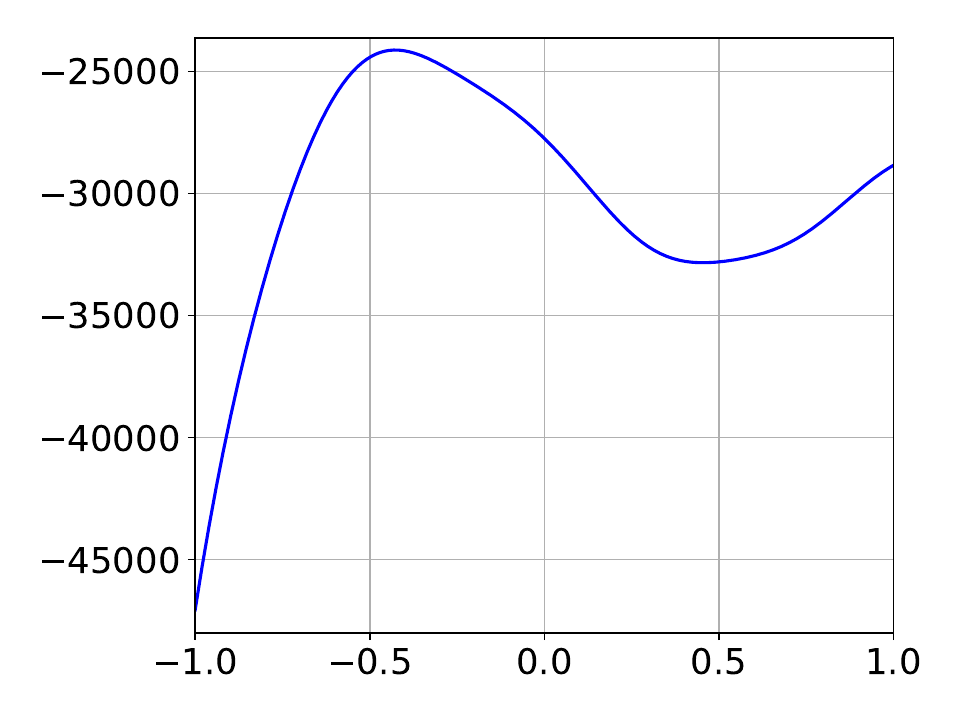}
\hfil

\caption{The derivatives of the four slowly-varying phase functions produced
by applying the global Levin method to Equation~(\ref{experiments4:1}) of Subsection~\ref{section:experiments:4}
when the parameter $\omega$ is equal to $2^{16}$.  Each column corresponds to one of the phase functions,
with the real part appearing in the first row and the imaginary part in the second.}
\label{experiments4:figure2}
\end{figure}




\end{subsection}

%
%
\begin{subsection}{An equation whose coefficient matrix has eigenvalues with large real parts}
\label{section:experiments:5}

The experiment of this section concerns the fourth order differential equation
\begin{equation}
y''''(t) + \omega^4\left(\frac{2+\cos(7t)^2}{1+t^4}\right) y(t) = 0.
\label{experiments5:1}
\end{equation}
When $\omega$ is large, some of the eigenvalues of the corresponding coefficient matrix have real parts of large magnitude,
with the consequence that almost all initial and boundary value problems for (\ref{experiments5:1}) are highly ill-conditioned.
This makes testing the accuracy of the obtained phase functions by using them to construct solutions
of (\ref{experiments5:1}) and comparing the results to reference solutions virtually impossible.
Instead, we assess the accuracy of the phase functions using our second approach, via comparison
with phase functions constructed by our algorithms using extended precision arithmetic.

More explicitly, for each $\omega=2^8,2^9,\ldots,2^{20}$, we constructed
slowly-varying phase functions for (\ref{experiments5:1}) over the interval $[-1,1]$ by running the global
method using double precision arithmetic.
We then measured the relative errors in each obtained phase function at $10,000$ equispaced points in
the interval $[-1,1]$  by comparison with phase functions constructed by applying the 
global Levin method to (\ref{experiments5:1}), but this time using extended precision arithmetic to perform the calculations.
Figure~\ref{experiments5:figure1} gives the results of these experiments.
Plots of the eigenvalues $\lambda_1(t),\lambda_2(t),\lambda_3(t),\lambda_4(t)$ of the coefficient matrix
for (\ref{experiments5:1}) when $\omega=2^{16}$ can be found in  Figure~\ref{experiments5:figure3} .

\begin{figure}[h]
\hfil
\includegraphics[width=.325\textwidth]{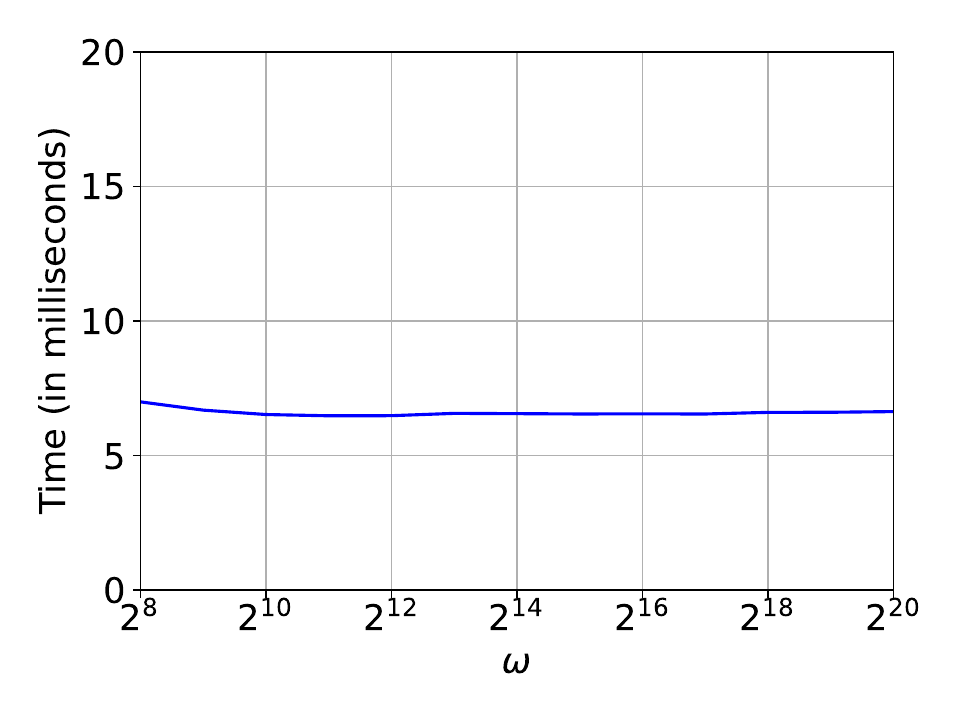}
\hfil
\includegraphics[width=.325\textwidth]{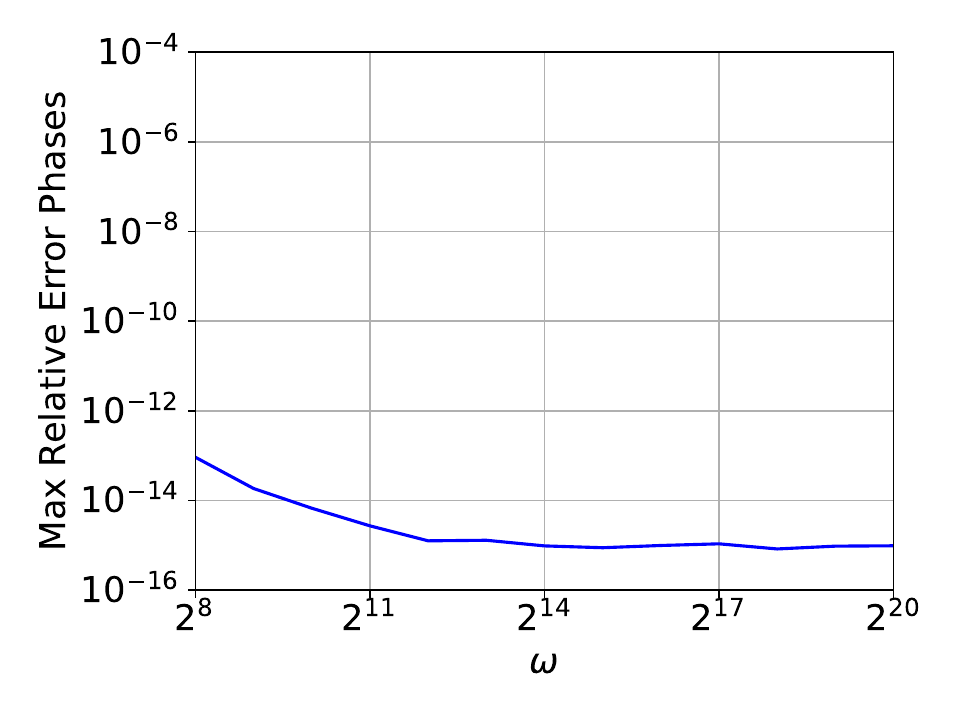}
\hfil
\includegraphics[width=.325\textwidth]{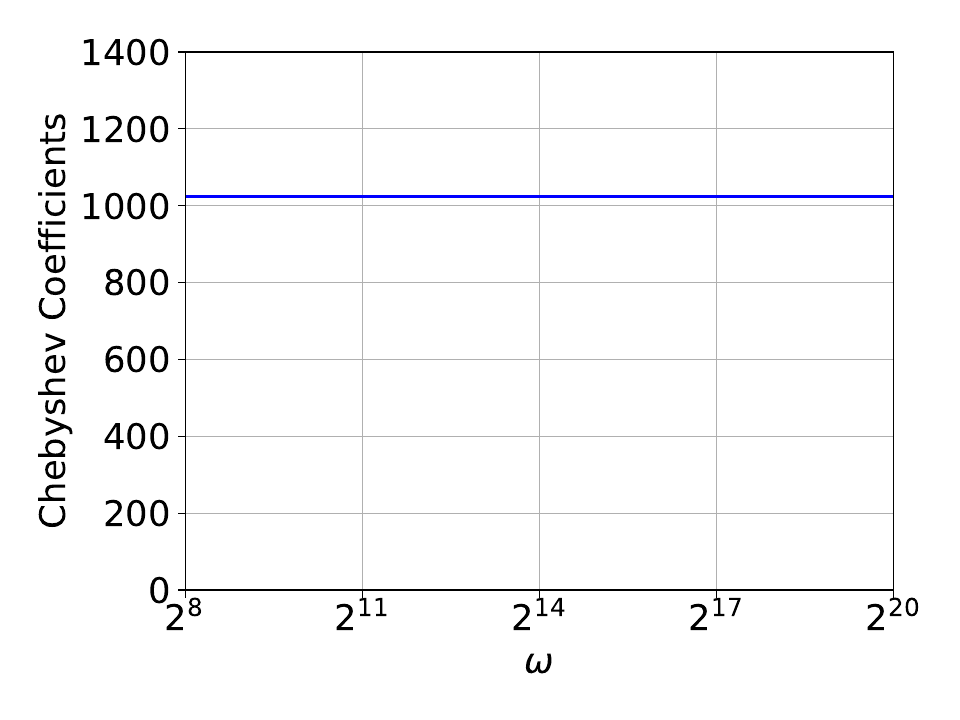}
\hfil
\caption{The results of the experiment of Subsection~\ref{section:experiments:5}.
The leftmost plot gives the time required by the global Levin method
as a function of the parameter $\omega$.   The middle plot reports the maximum relative error in the slowly-varying phase functions.
The plot on  the right shows the number  of Chebyshev coefficients required to represent the slowly-varying phase functions, again
as a function of $\omega$.}
\label{experiments5:figure1}
\end{figure}




\begin{figure}[h]

\hfil
\includegraphics[width=.24\textwidth]{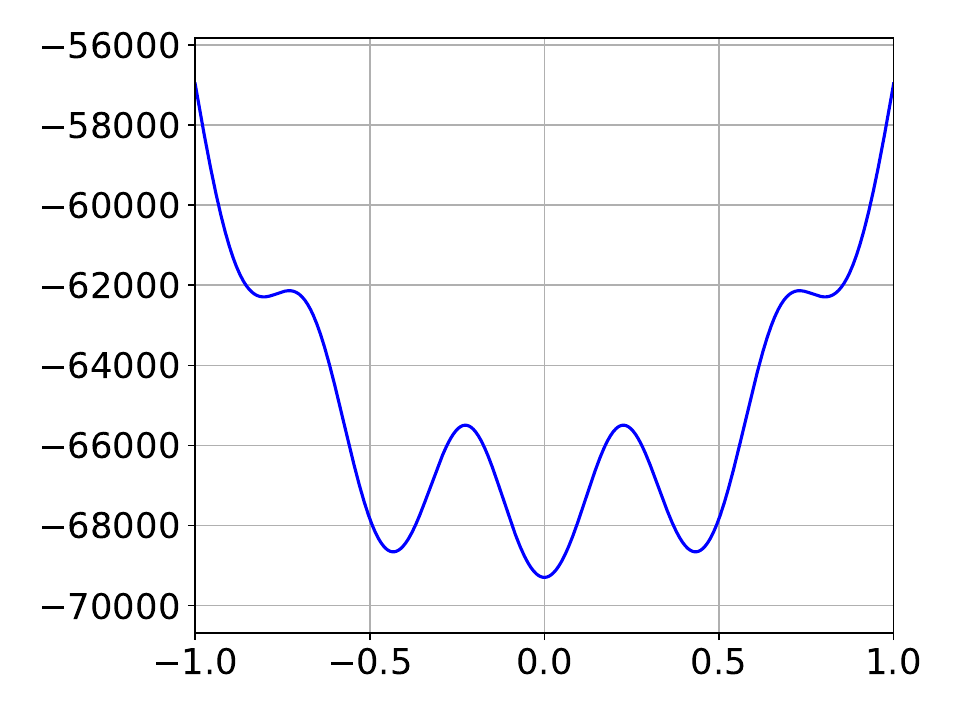}
\hfil
\includegraphics[width=.24\textwidth]{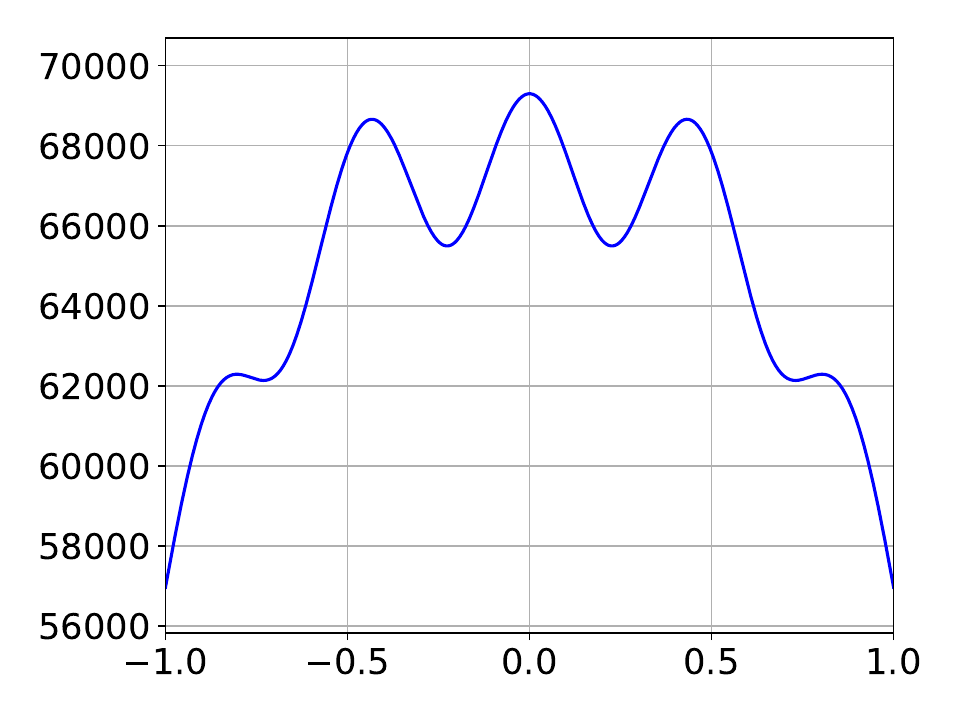}
\hfil
\includegraphics[width=.24\textwidth]{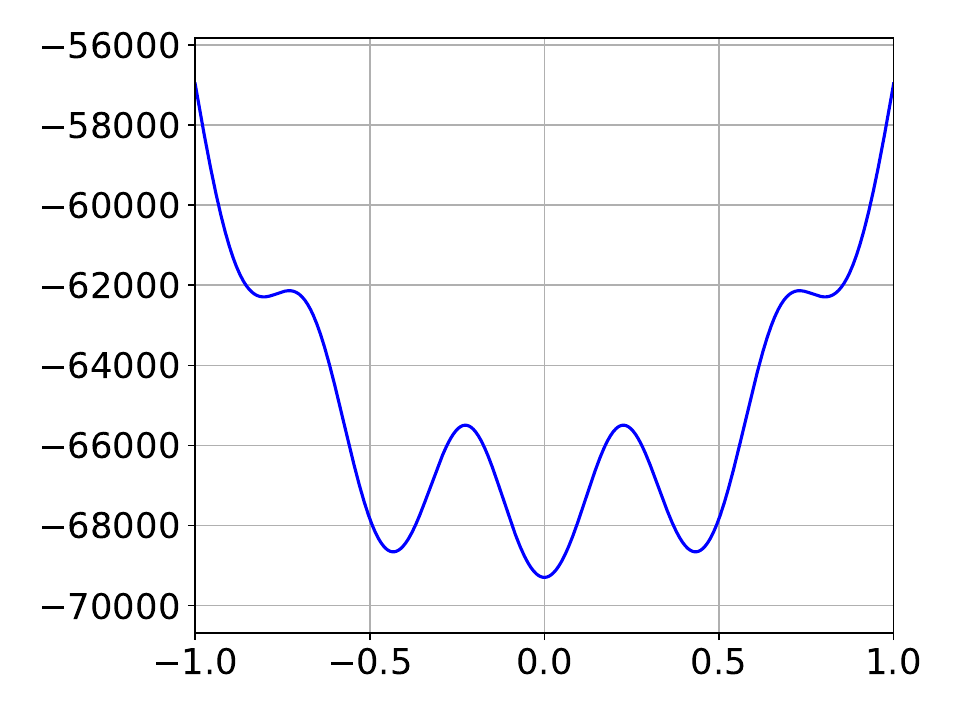}
\hfil
\includegraphics[width=.24\textwidth]{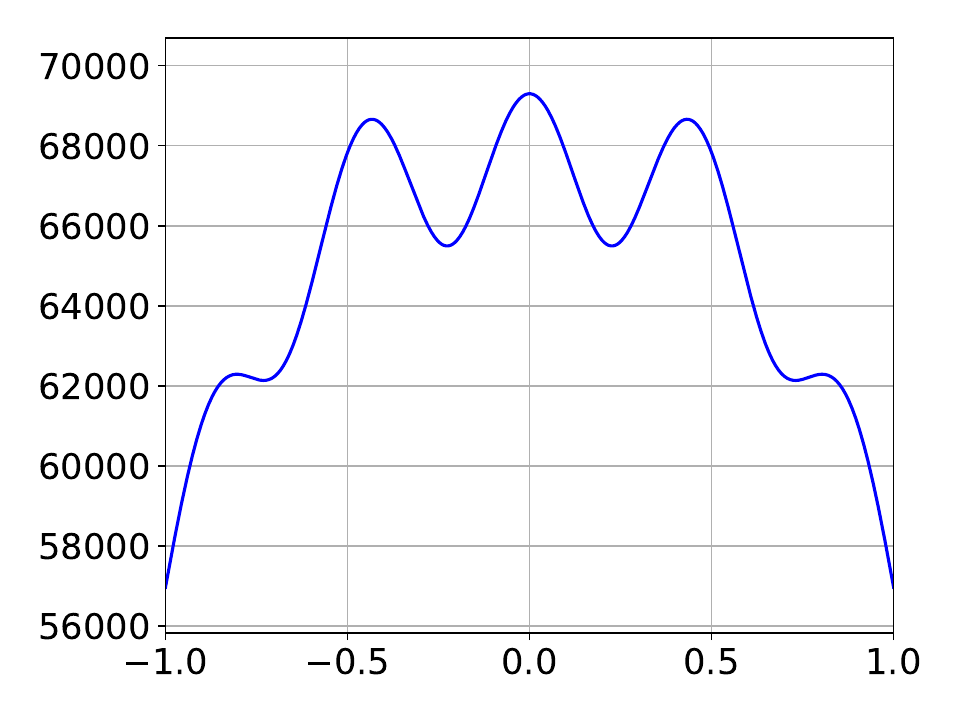}
\hfil

\hfil
\includegraphics[width=.24\textwidth]{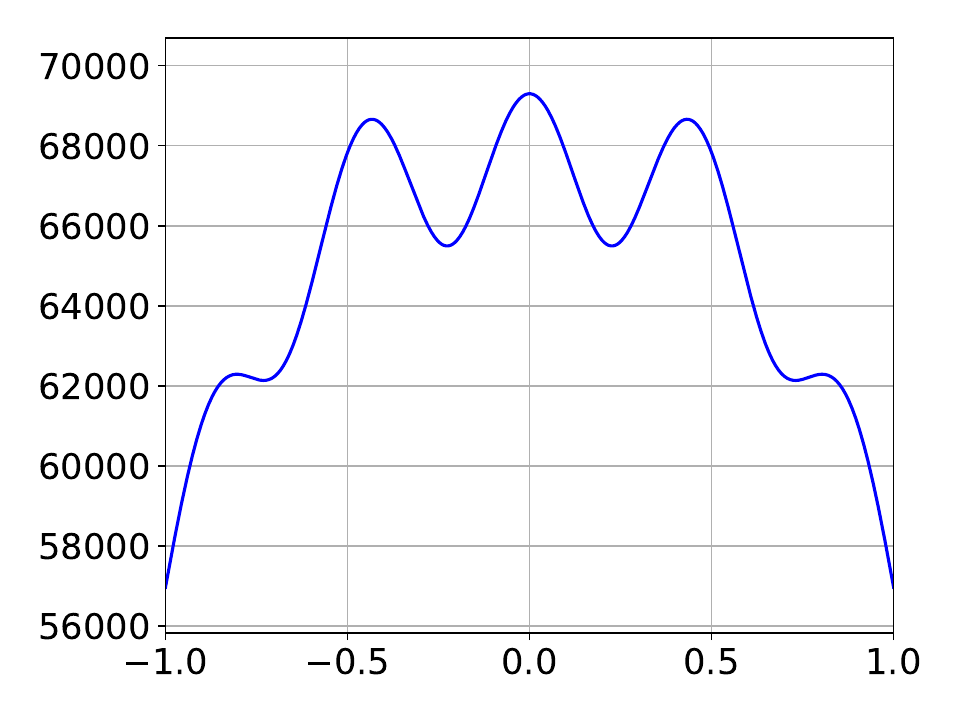}
\hfil
\includegraphics[width=.24\textwidth]{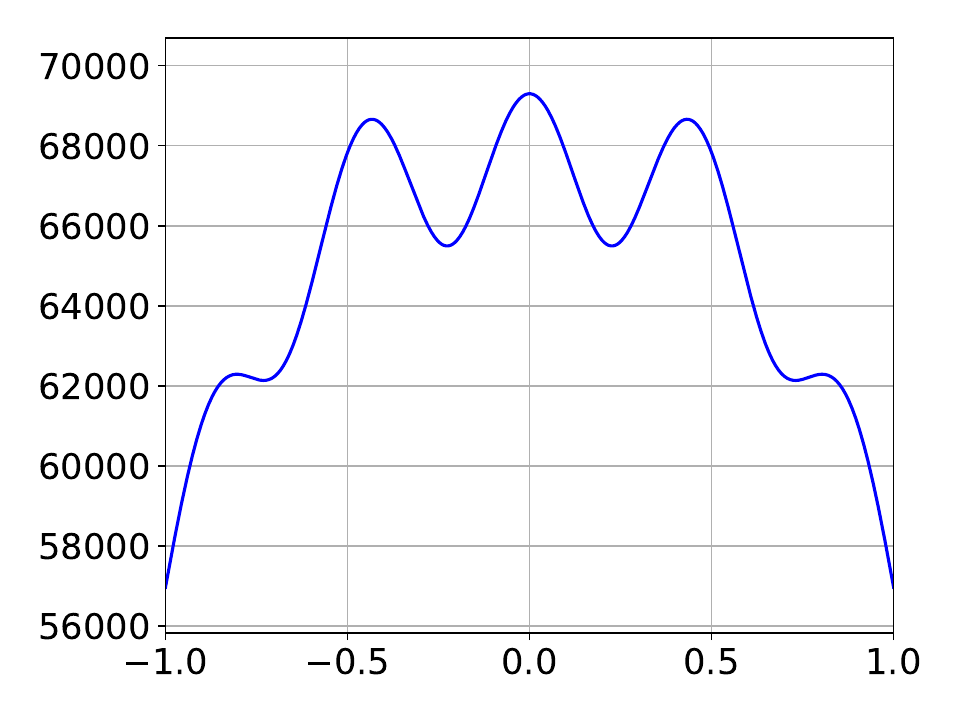}
\hfil
\includegraphics[width=.24\textwidth]{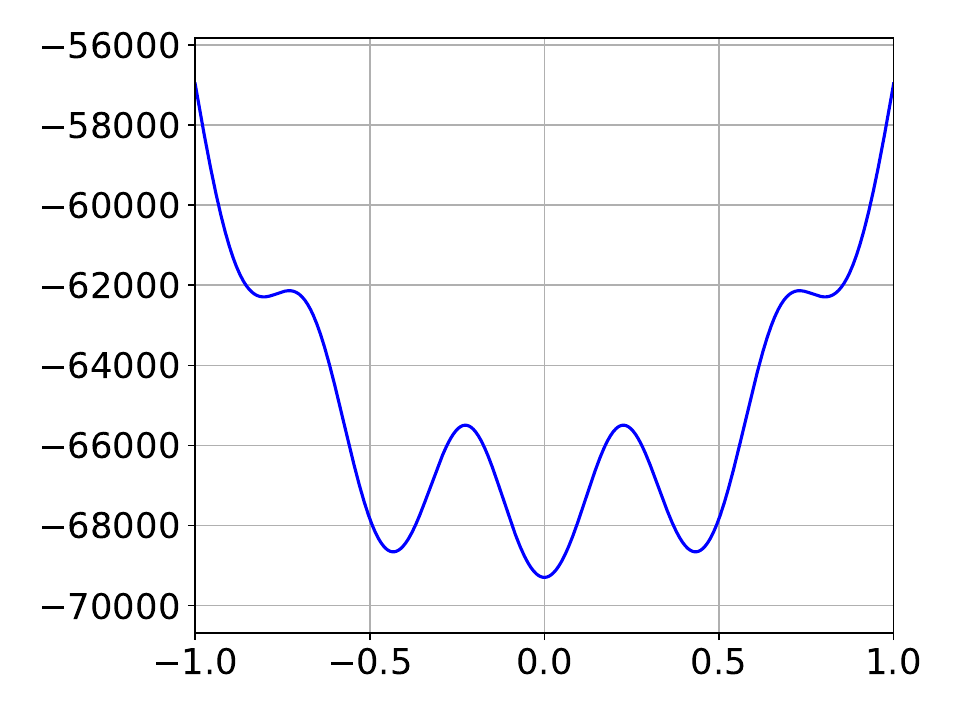}
\hfil
\includegraphics[width=.24\textwidth]{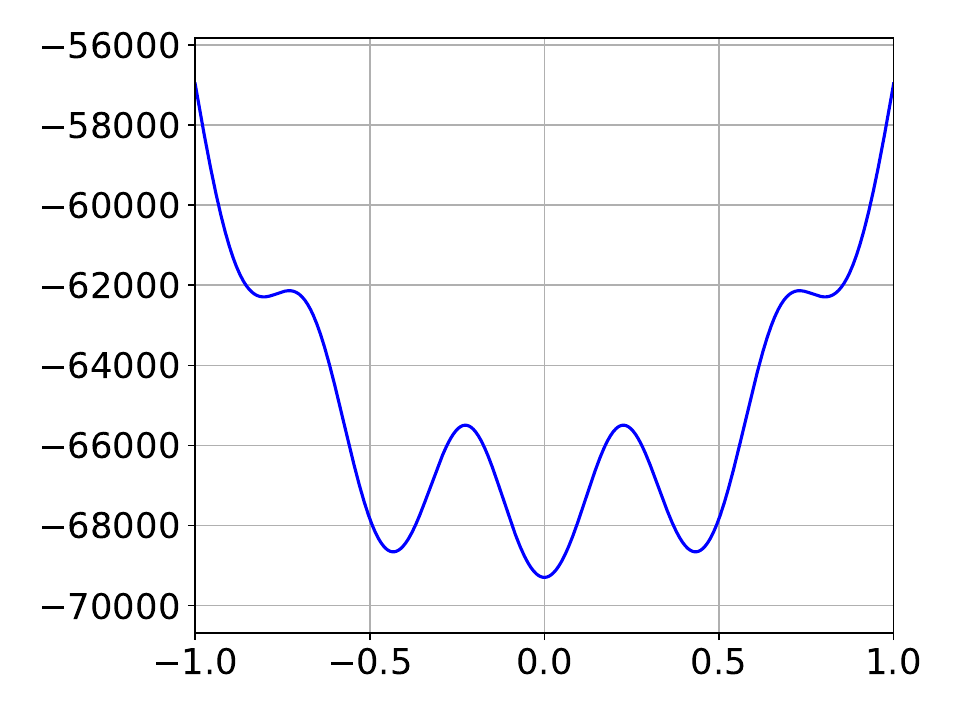}
\hfil
\caption{
The eigenvalues $\lambda_1(t), \lambda_2(t), \lambda_3(t), \lambda_4(t)$ of the coefficient 
matrix for (\ref{experiments5:1}) 
when the parameter $\omega$ is equal to $2^{16}$.  Each column corresponds to one of the eigenvalues,
with the real part appearing in the first row and the imaginary part in the second.}
\label{experiments5:figure3}
\end{figure}

\end{subsection}

%
%
\begin{subsection}{An equation whose coefficient matrix has eigenvalues of small magnitude}
\label{section:experiments:6}

In the experiment described in this section, we solved the equation
\begin{equation}
y'''(t)-i\omega \left(1+t^2\right) y''(t) + \frac{2+t}{1+t^2}\, y'(t) +i\omega \log\left(\frac{3}{2}+t\right) y(t) = 0
\label{experiments6:1}
\end{equation}
over the interval $[-1,1]$  subject to the conditions
\begin{equation}
\begin{aligned}
y(0)   = 1,  \ \ \ y'(0)  = -i \omega,\ \ \  y''(0) = -\omega^2.
\end{aligned}
\label{experiments6:2}
\end{equation}
The coefficient matrix corresponding to (\ref{experiments6:1}) has eigenvalues of small magnitude,
and, as discussed in the introduction, the global Levin method encounters difficulties in such cases.  
In particular, the obtained functions can be discontinuous across
the boundaries of the subintervals set by the  adaptive discretization scheme.
This is readily apparent in Figure~\ref{experiments6:figure2}, which
contains plots of the functions that result
when the global Levin method is applied to (\ref{experiments6:1}) with $\omega=2^{16}$.
Fortunately, the local Levin method has no difficulties in this case;
Figure~\ref{experiments6:figure3} contains plots of the functions obtained
when it is applied to (\ref{experiments6:1}) with $\omega=2^{16}$.


For each $\omega=2^8, 2^9, \ldots, 2^{20}$, we used the local Levin method to solve
(\ref{experiments6:1}) subject to (\ref{experiments6:2}).  
 The error in each obtained solution was measured
at $10,000$ equispaced points on the interval $[-1,1]$.
Figure~\ref{experiments6:figure1} gives the results.
Plots of the eigenvalues $\lambda_1(t),\lambda_2(t),\lambda_3(t),\lambda_4(t)$ of the coefficient matrix
for (\ref{experiments6:1}) when $\omega=2^{16}$ can be found in  Figure~\ref{experiments6:figure3} .

\begin{figure}[h]
\hfil
\includegraphics[width=.325\textwidth]{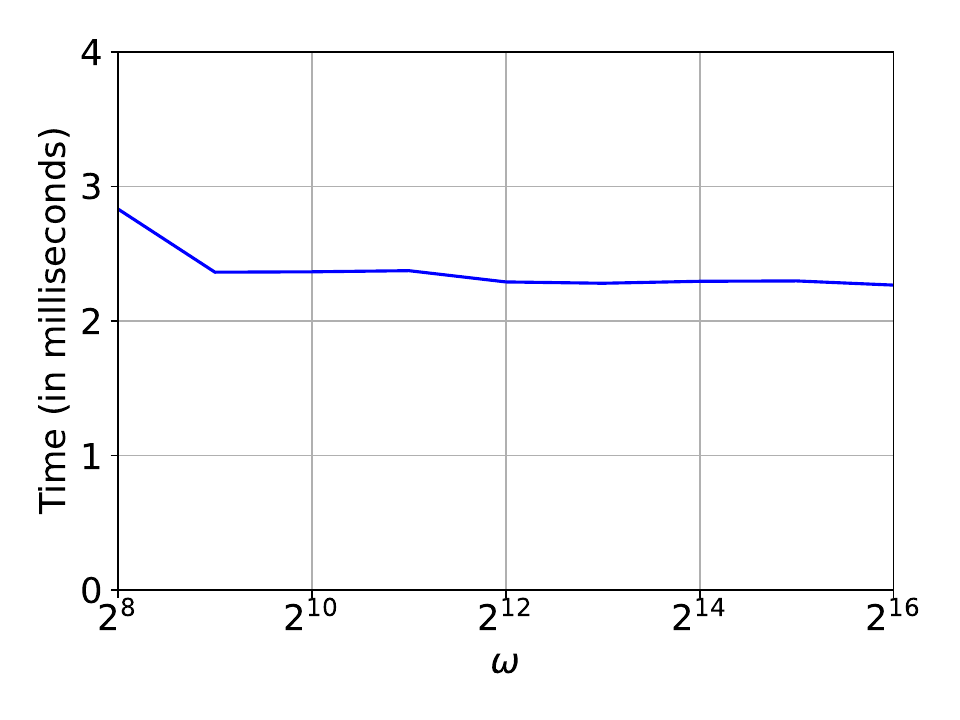}
\hfil
\includegraphics[width=.325\textwidth]{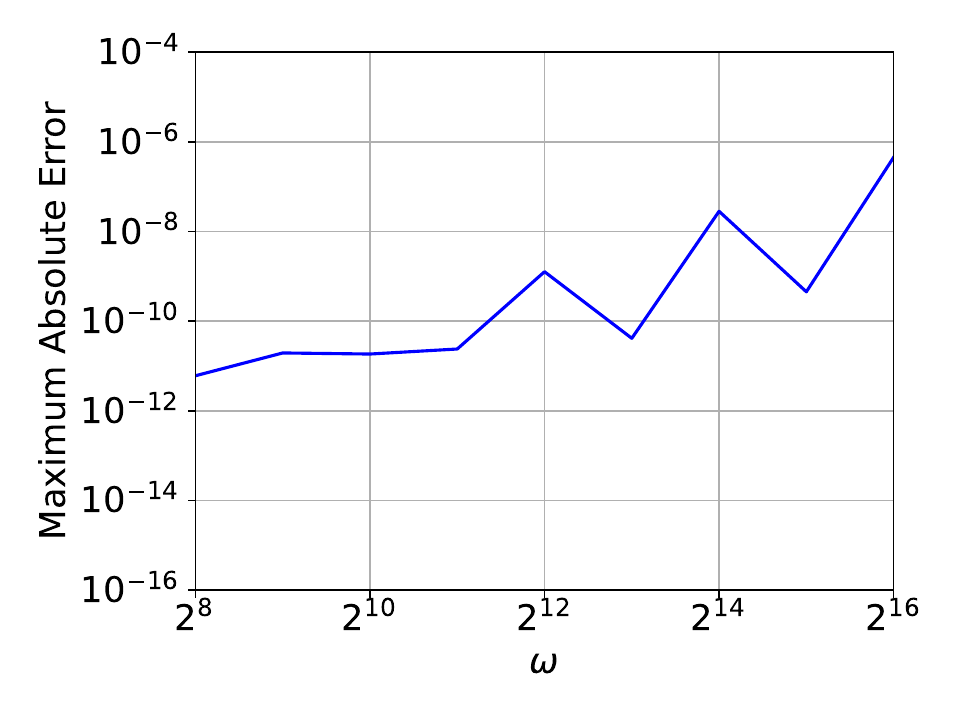}
\hfil
\includegraphics[width=.325\textwidth]{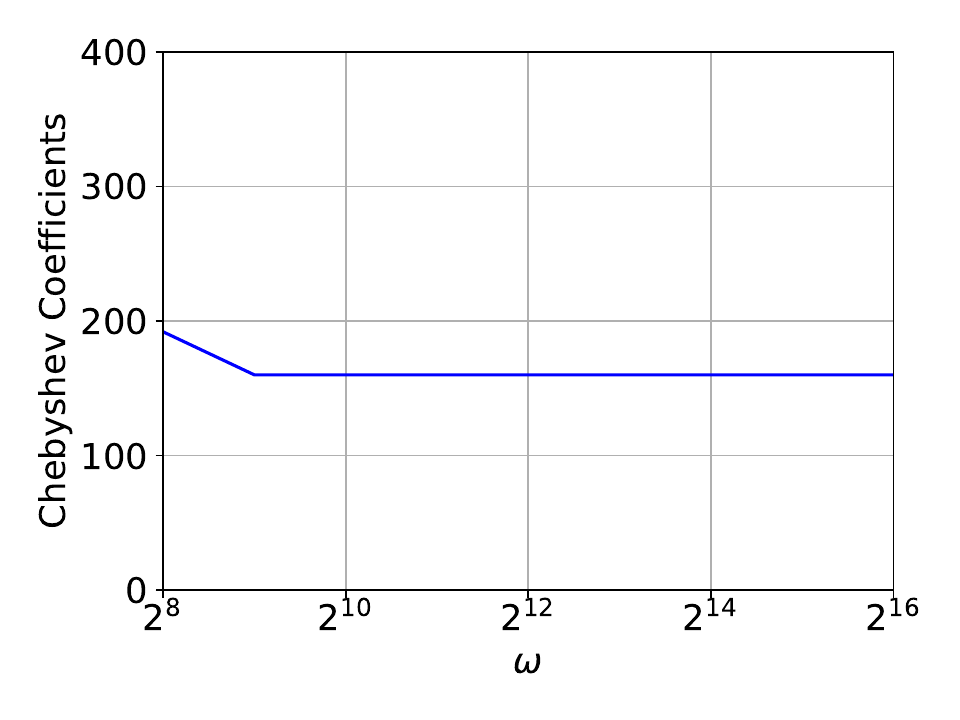}
\hfil
\caption{The results of the experiment of Subsection~\ref{section:experiments:6}.
The leftmost plot gives the time required by the local Levin method
as a function of the parameter $\omega$.  The plot in the middle reports the largest observed
absolute error in the solution of the initial value problem for
(\ref{experiments6:1}) as a function of $\omega$.  The rightmost plot  shows the total number
of Chebyshev coefficients required to represent the slowly-varying phase functions, again
as a function of $\omega$. }
\label{experiments6:figure1}
\end{figure}

\begin{figure}[h]
\hfil
\includegraphics[width=.325\textwidth]{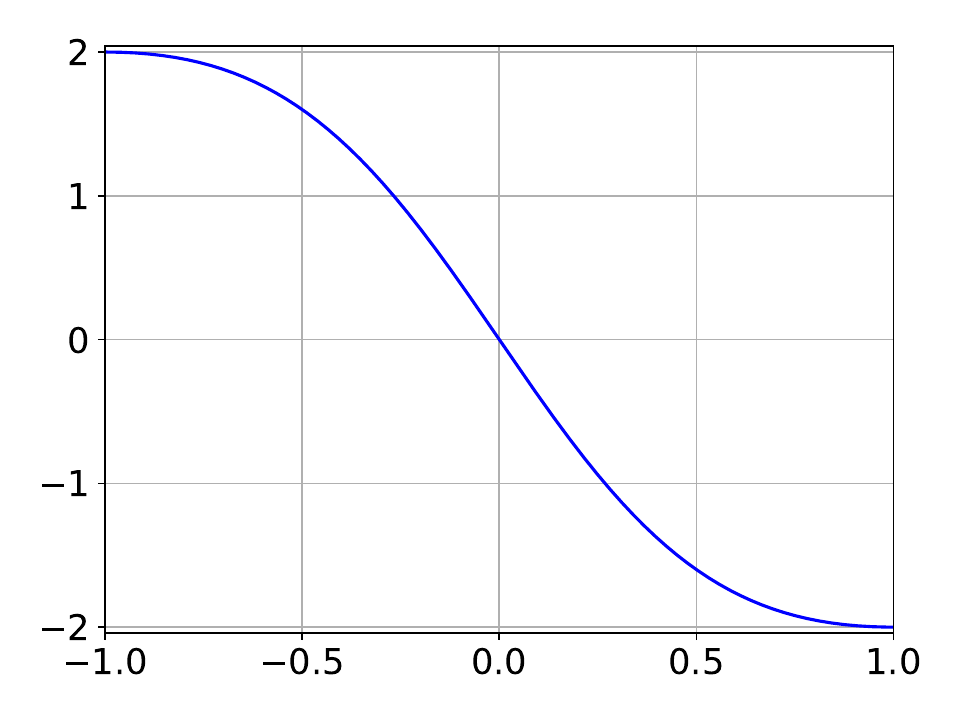}
\hfil
\includegraphics[width=.325\textwidth]{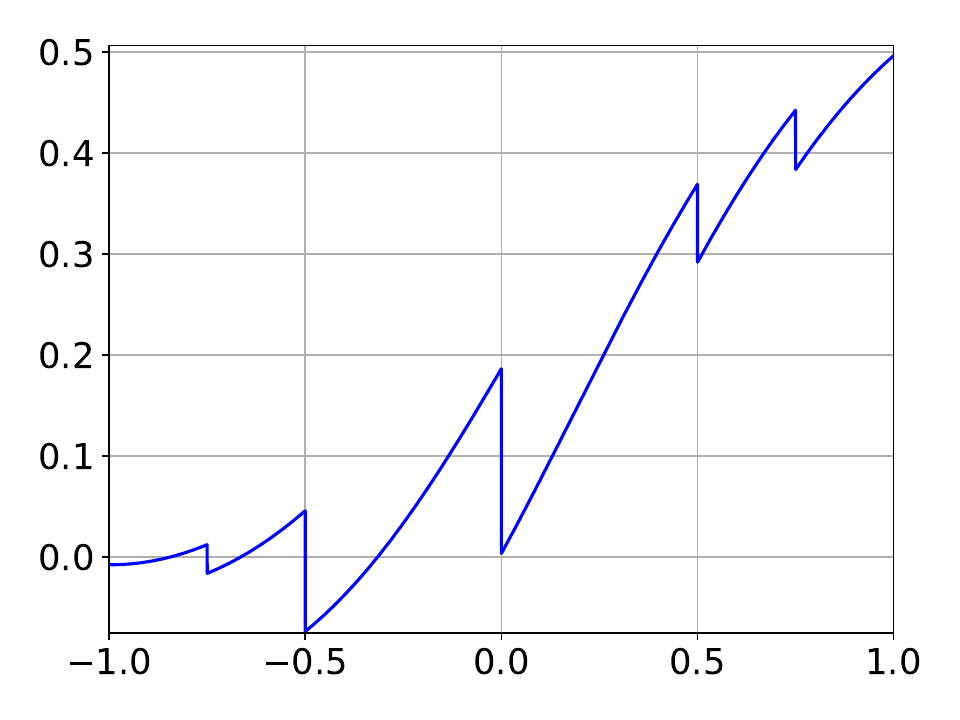}
\hfil
\includegraphics[width=.325\textwidth]{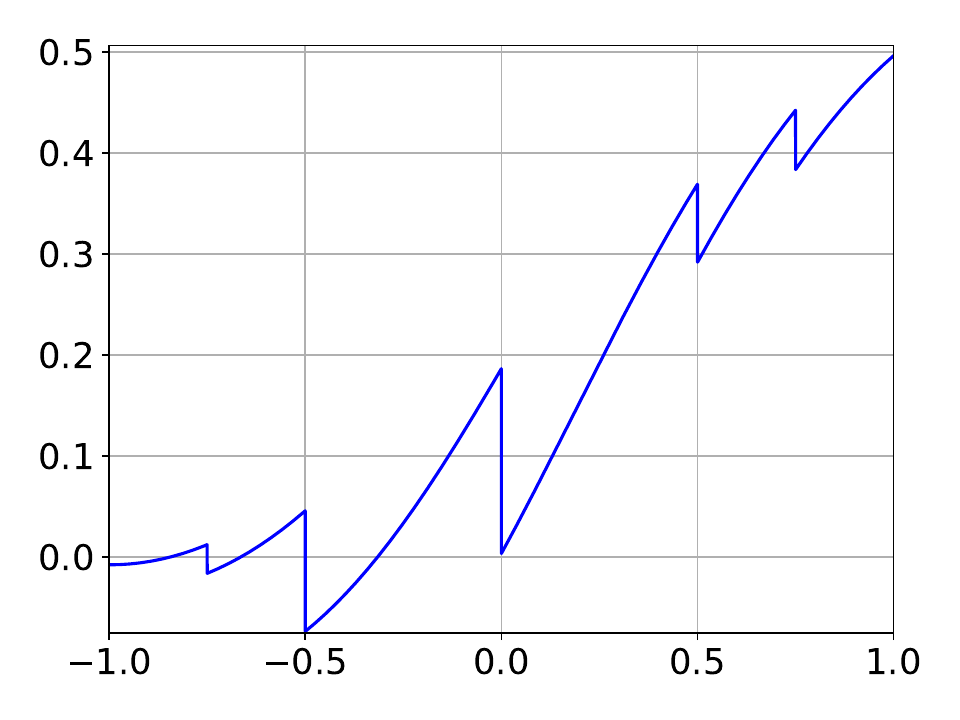}
\hfil

\hfil
\includegraphics[width=.325\textwidth]{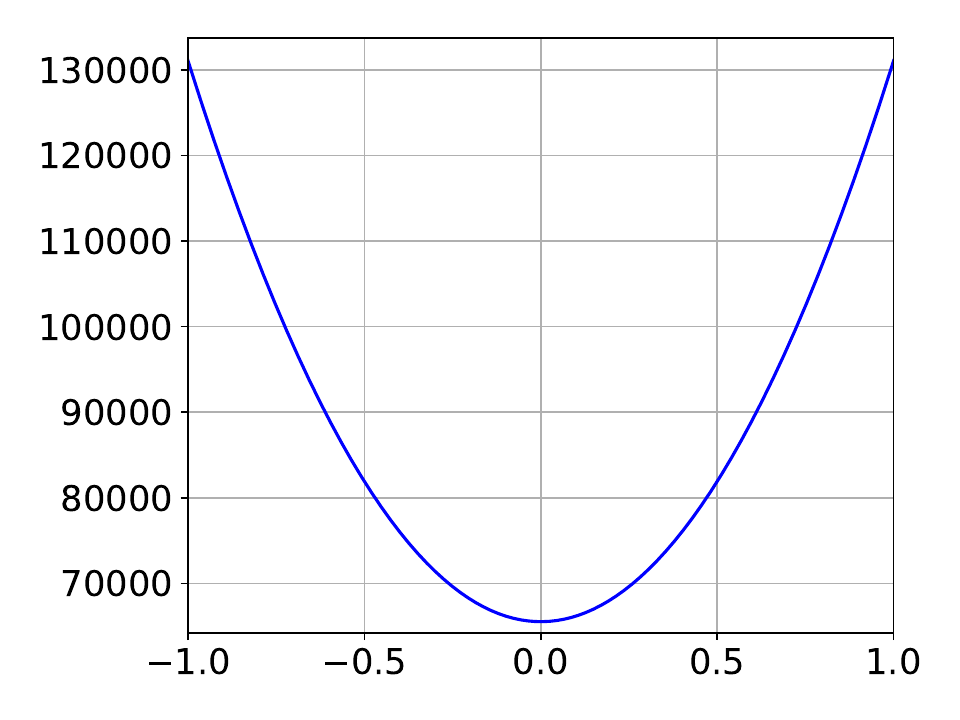}
\hfil
\includegraphics[width=.325\textwidth]{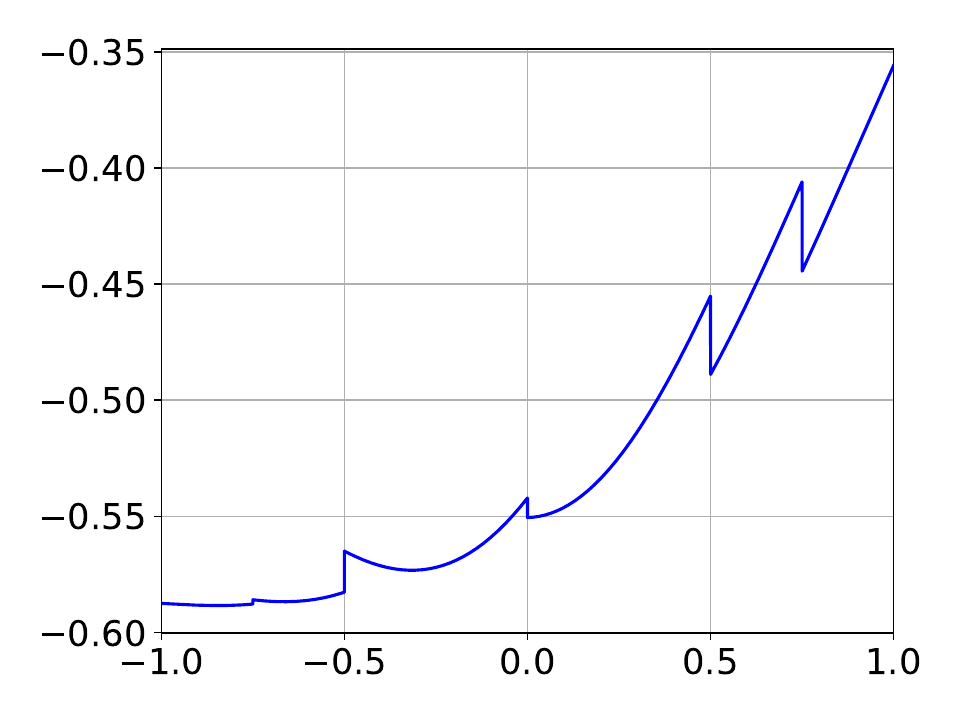}
\hfil
\includegraphics[width=.325\textwidth]{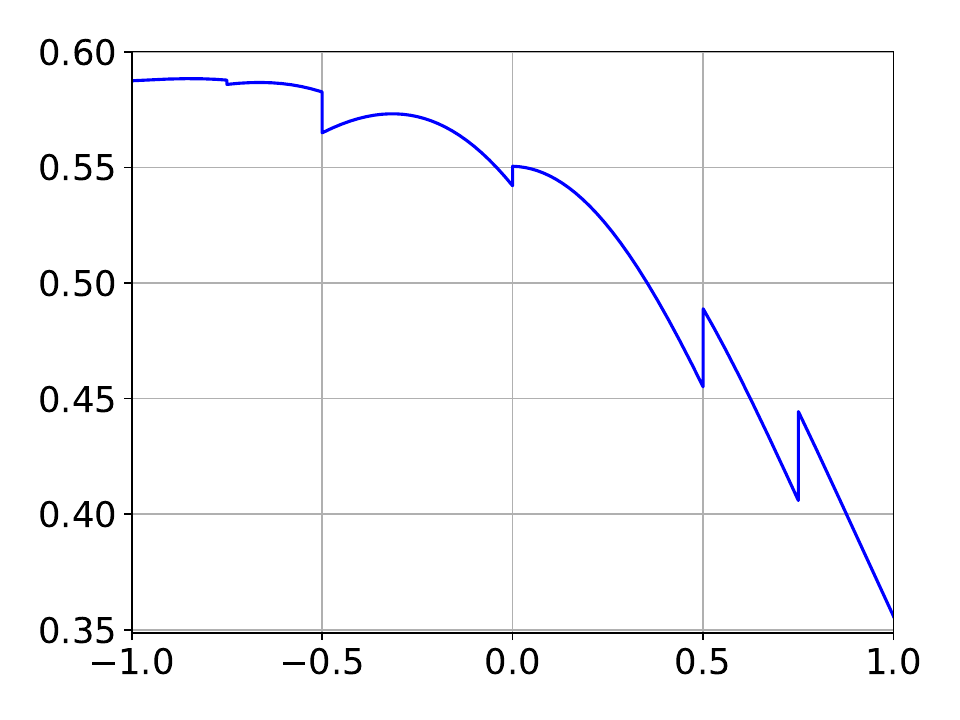}
\hfil

\caption{The degenerate results obtained 
when the  global Levin method is applied to Equation~(\ref{experiments6:1}) of Subsection~\ref{section:experiments:6}
with $\omega=2^{16}$.  In this case, the coefficient matrix admits eigenvalues
of small magnitude and the functions obtained by the global Levin method are contaminated by 
elements of the nullspaces of the linearized operators which arise in the course of the Newton procedure.
Each column corresponds to one of the phase functions,
with the real part appearing in the first row and the imaginary part in the second.}
\label{experiments6:figure2}
\end{figure}

\begin{figure}[h]
\hfil
\includegraphics[width=.325\textwidth]{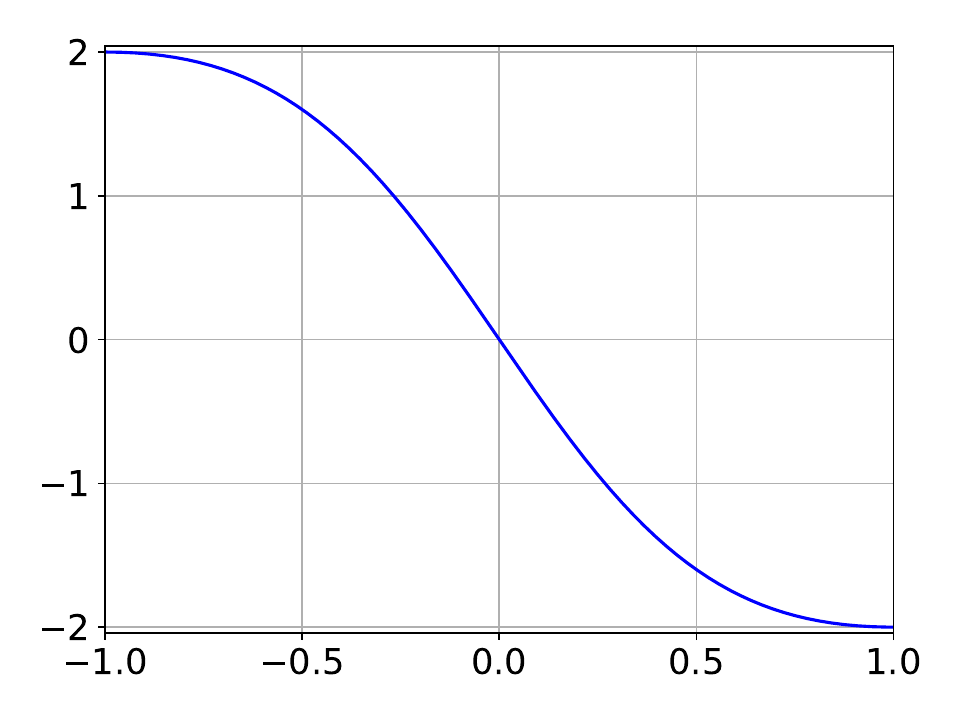}
\hfil
\includegraphics[width=.325\textwidth]{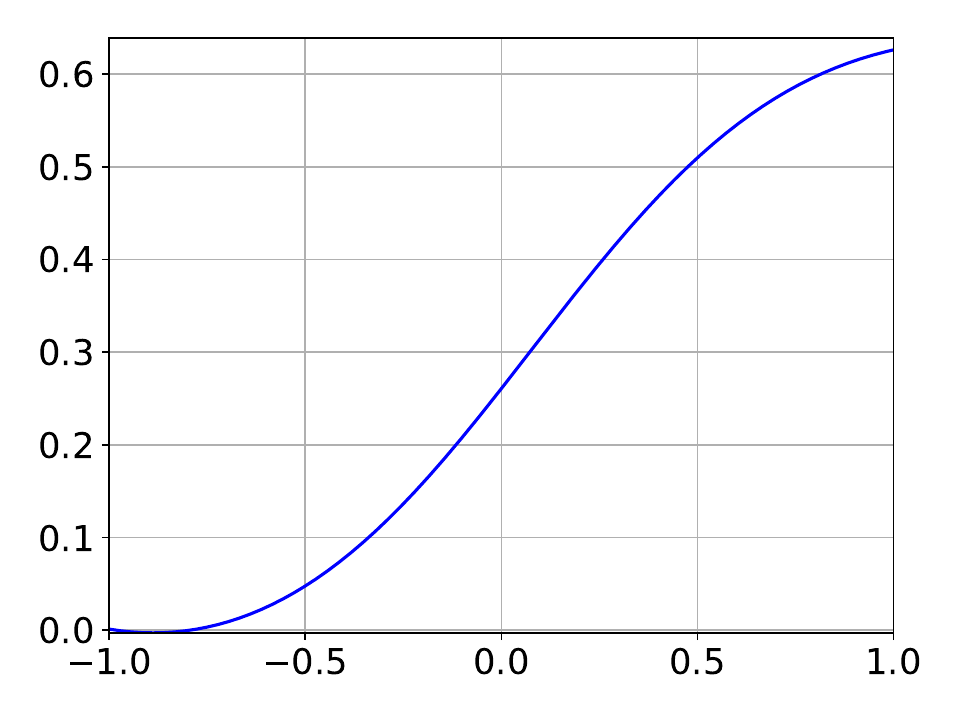}
\hfil
\includegraphics[width=.325\textwidth]{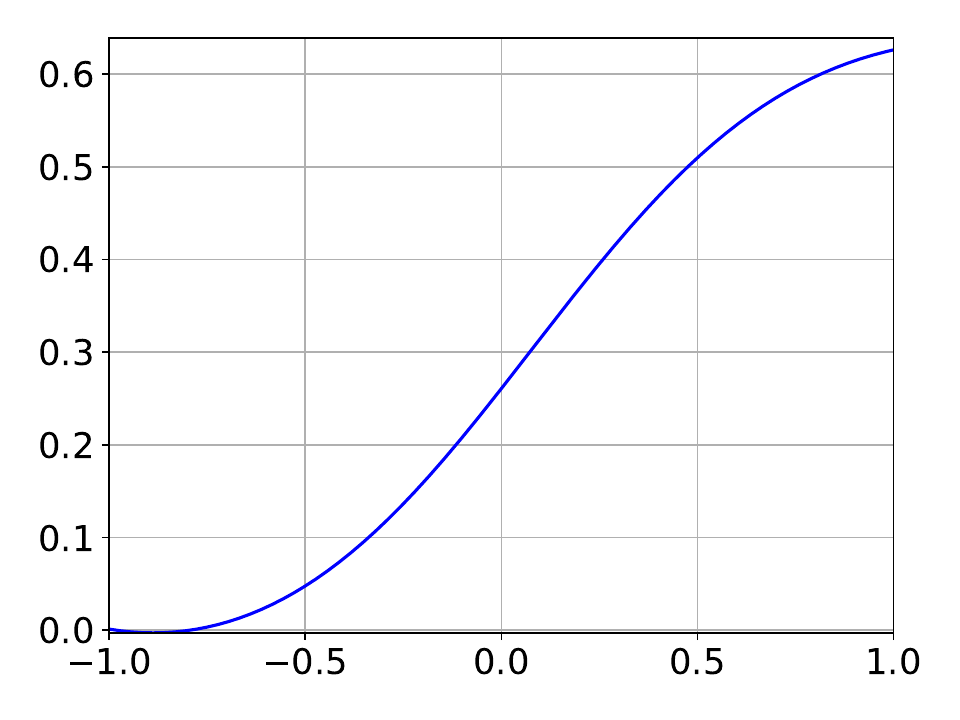}
\hfil

\hfil
\includegraphics[width=.325\textwidth]{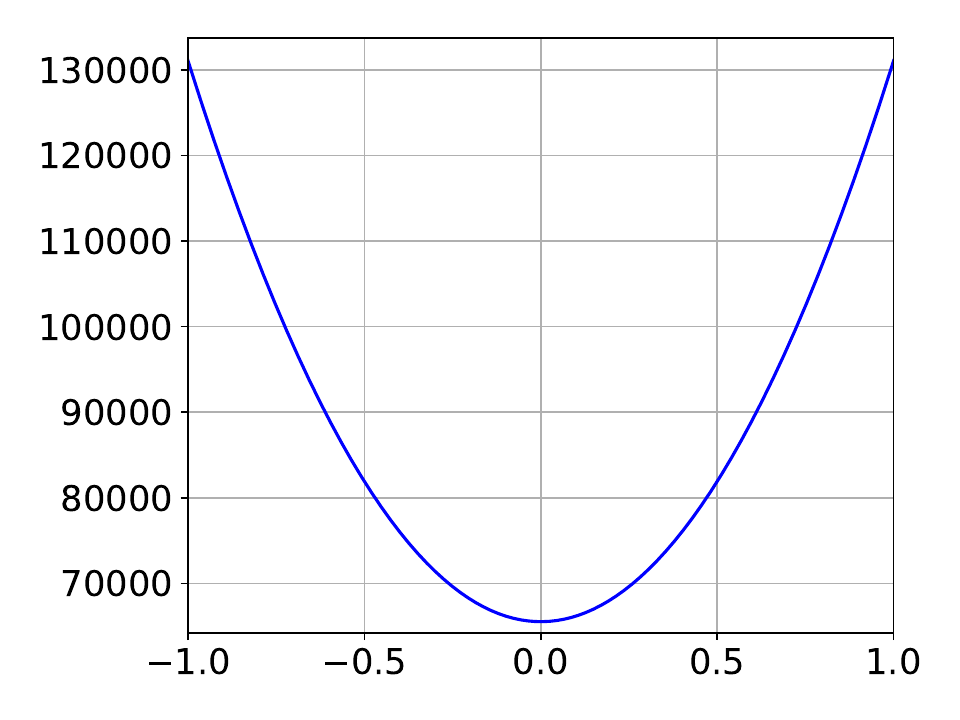}
\hfil
\includegraphics[width=.325\textwidth]{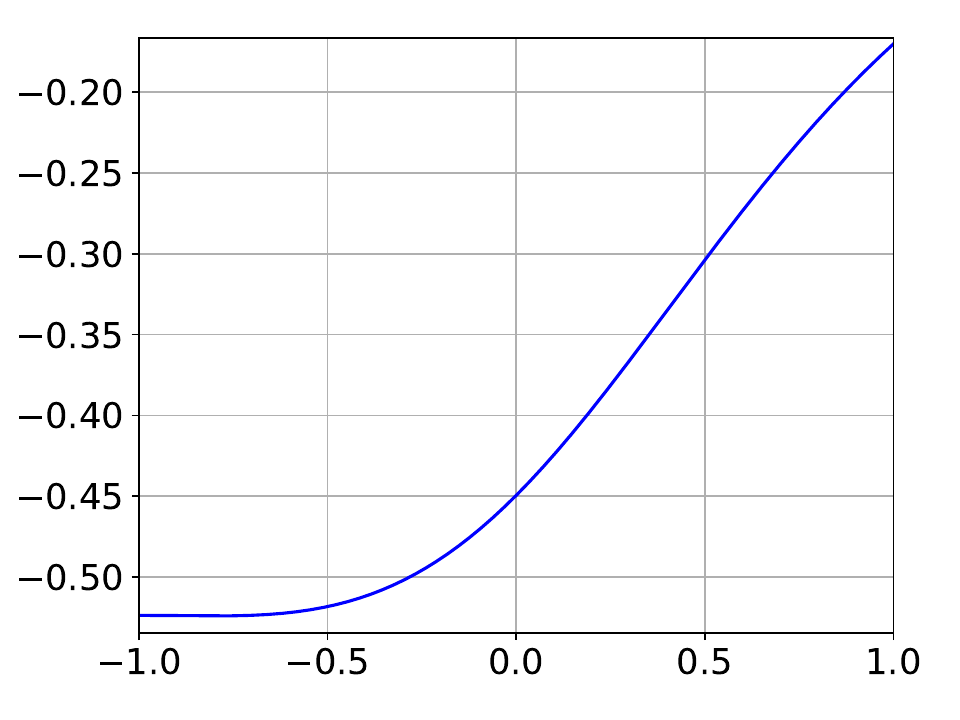}
\hfil
\includegraphics[width=.325\textwidth]{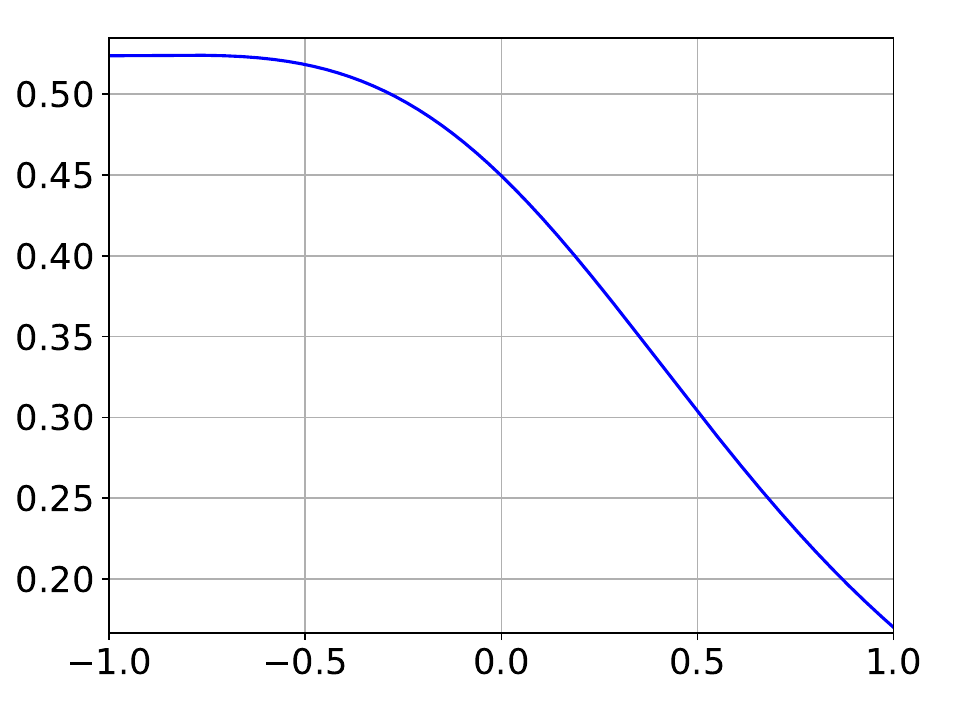}
\hfil

\caption{The derivatives of the slowly-varying phase functions obtained 
when the local  Levin method is applied to Equation~(\ref{experiments6:1}) of Subsection~\ref{section:experiments:6}
with $\omega=2^{16}$.  Even though the coefficient matrix admits eigenvalues
of small magnitude, the local Levin method is still able to construct the desired slowly-varying phase functions.  
Each column corresponds to one of the phase functions,
with the real part appearing in the first row and the imaginary part in the second.}
\label{experiments6:figure3}
\end{figure}

\end{subsection}

%
%

\begin{subsection}{A system of two differential equations in two unknowns}
\label{section:experiments:7}

It is well known that essentially any system of $n$ linear ordinary differential equations
in $n$ unknowns can be transformed into an $n^{th}$ order linear scalar differential equation
(see, for instance, \cite{Kolchin}).
As a consequence, the algorithms of this paper can be used to solve many such  systems
in time independent of the magnitudes of the eigenvalues of their coefficient matrices.    
In the experiment of this section, we solved the system of differential equations
\begin{equation}
\mathbf{y}'(t) = 
\left(
\begin{array}{cc}
1+t^2            & \frac{1}{1+t^4} \\
\frac{-\omega}{1+t^2} & \frac{-i\omega (2+t)}{5+t}
\end{array}
\right)
\mathbf{y}(t)
\label{experiments7:system}
\end{equation}
on the interval $[-1,1]$ subject to the condition
\begin{equation}
\mathbf{y}(0) = 
\left(
\begin{array}{c}
1\\
1
\end{array}
\right).
\end{equation}
We found that letting
\begin{equation}
\mathbf{w}(t) = \Phi(t) \mathbf{y}(t)
\end{equation}
with
\begin{equation}
\Phi(t) = \left(\begin{array}{cc}
0 & 1 \\
\frac{-\omega}{1+t^2} & \frac{-i\omega(2+t)}{5+t}
\end{array}\right)
\label{experiments7:phi}
\end{equation}
transforms (\ref{experiments7:system}) into the system
\begin{equation}
\mathbf{w}'(t) = 
\left(
\begin{array}{cc}
0            & 1 \\
-q_0(t)      & -q_1(t)
\end{array}
\right)
\mathbf{w}(t),
\label{experiments7:companion}
\end{equation}
where
\begin{equation*}
\begin{aligned}
q_0(t) =  
\frac{1}{(t+5)^2\left(t^2+1\right) \left(t^4+1\right)}
\left(
-i \omega t^{10}-7 i \omega t^9-12 i \omega t^8-12 i \omega t^7-5 i \omega t^6+  \right. \\  \left.
6 i \omega t^5-19
   i \omega t^4-12 i \omega t^3+(1-4 i) \omega t^2+(10+13 i) \omega t+(25-7 i) \omega \right)
\end{aligned}
\end{equation*}
and
\begin{equation*}
q_1(t) = \frac{i \omega (t+2)}{t+5}-t^2+\frac{2 t}{t^2+1}-1.
\end{equation*}
Of course, (\ref{experiments7:companion}) is equivalent to the scalar equation
\begin{equation}
z''(t) + q_1(t) z'(t) + q_0(t) z(t) = 0
\label{experiments7:scalar}
\end{equation}
with $z$ related to $\mathbf{w}=\left(\begin{array}{cc} w_1(t) & w_2(t) \end{array}\right)$ via the formulas
\begin{equation}
z(t) = w_1(t) \ \ \ \mbox{and} \ \ \ z'(t) = w_2(t).
\end{equation}

\begin{figure}[h]
\hfil
\includegraphics[width=.325\textwidth]{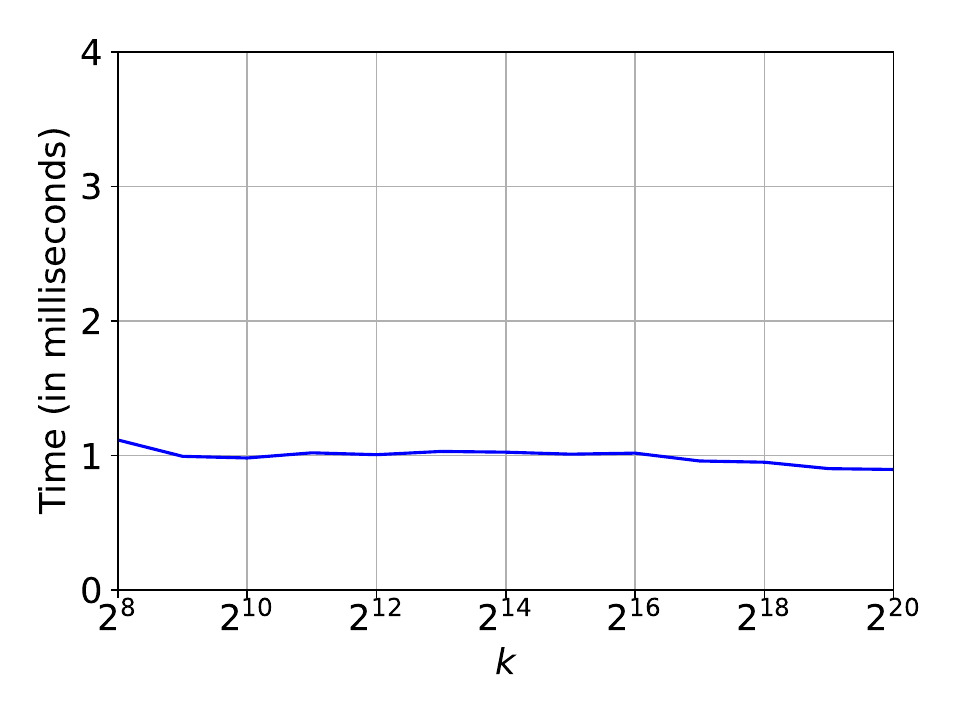}
\hfil
\includegraphics[width=.325\textwidth]{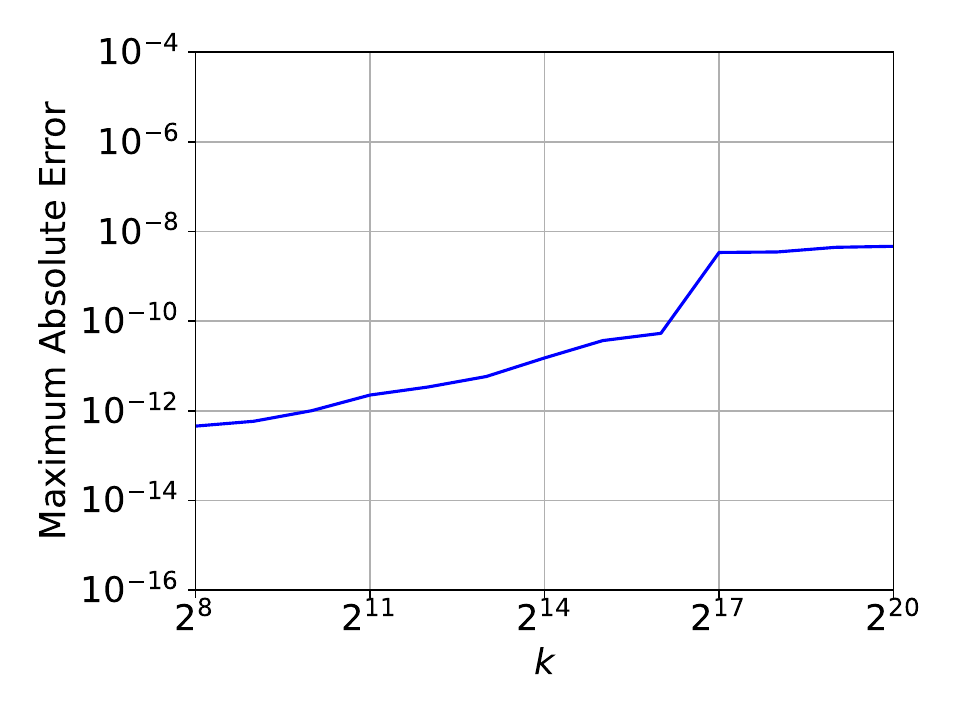}
\hfil
\includegraphics[width=.325\textwidth]{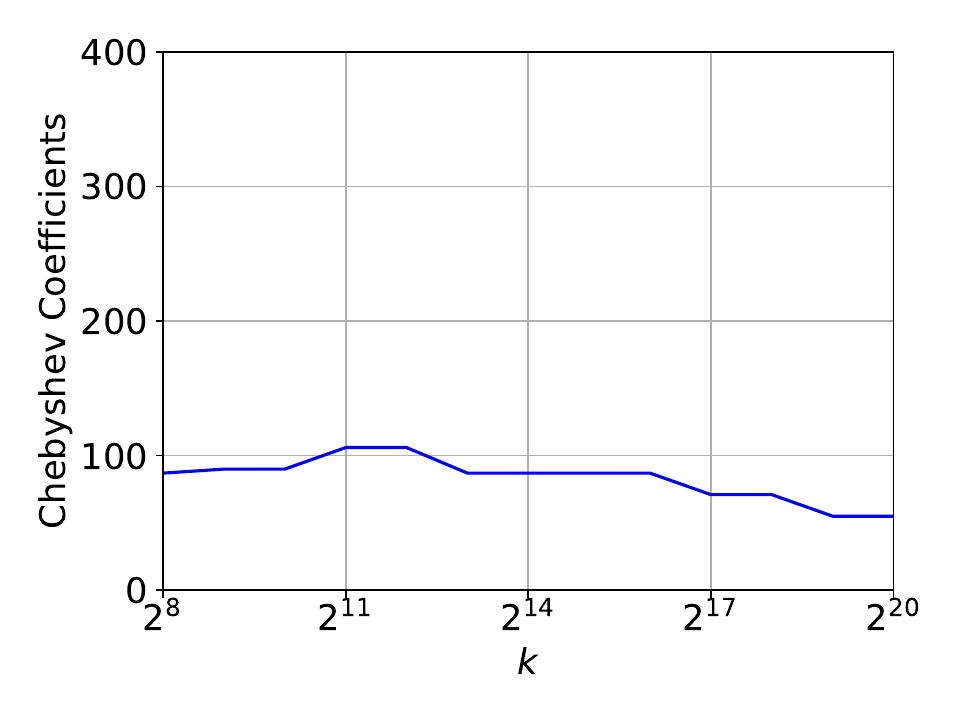}
\hfil
\caption{The results of the experiment of Subsection~\ref{section:experiments:7}.
The leftmost plot gives the time required by each of our methods
as a function of the parameter $\omega$.  The plot in the middle reports the largest observed
error in the solution as a function of $\omega$.  The plot on the right shows the total number
of Chebyshev coefficients required to represent the slowly-varying phase functions, again
as a function of $\omega$. }
\label{experiments7:figure1}
\end{figure}

\begin{figure}[h]
\centering
\hfil
\begin{subfigure}{.45\textwidth}
\hfil
\includegraphics[width=.45\textwidth]{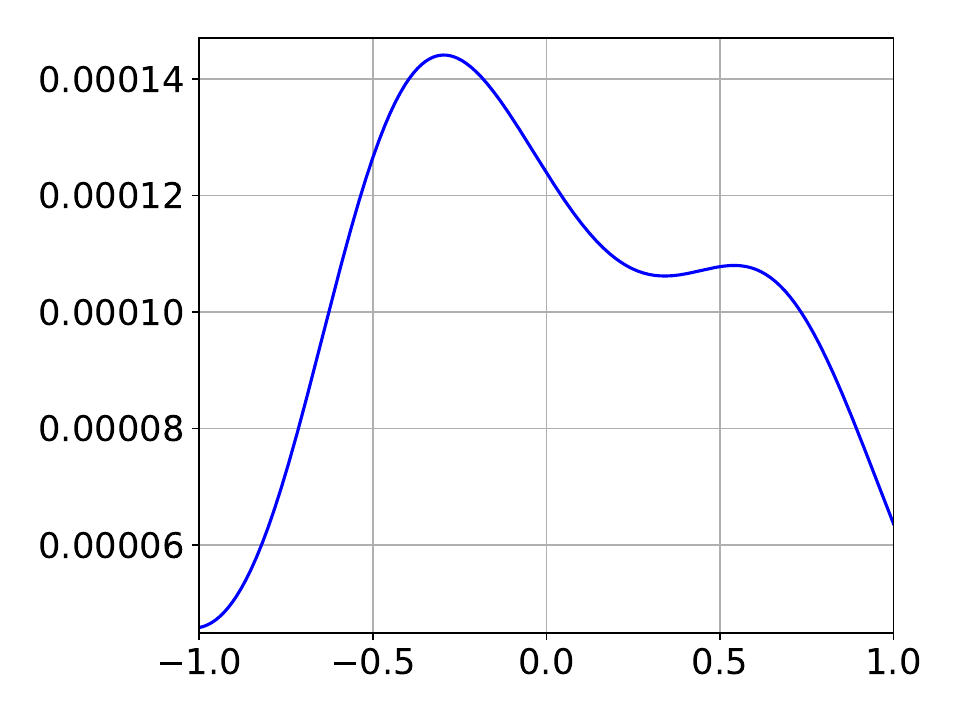}
\hfil
\includegraphics[width=.45\textwidth]{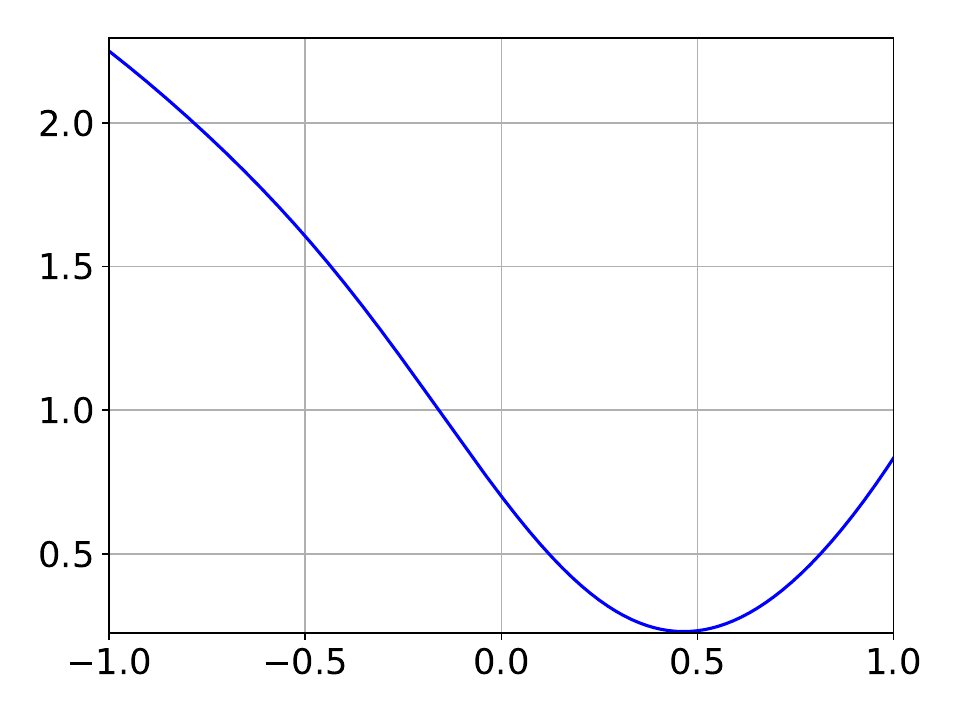}
\hfil

\hfil
\includegraphics[width=.45\textwidth]{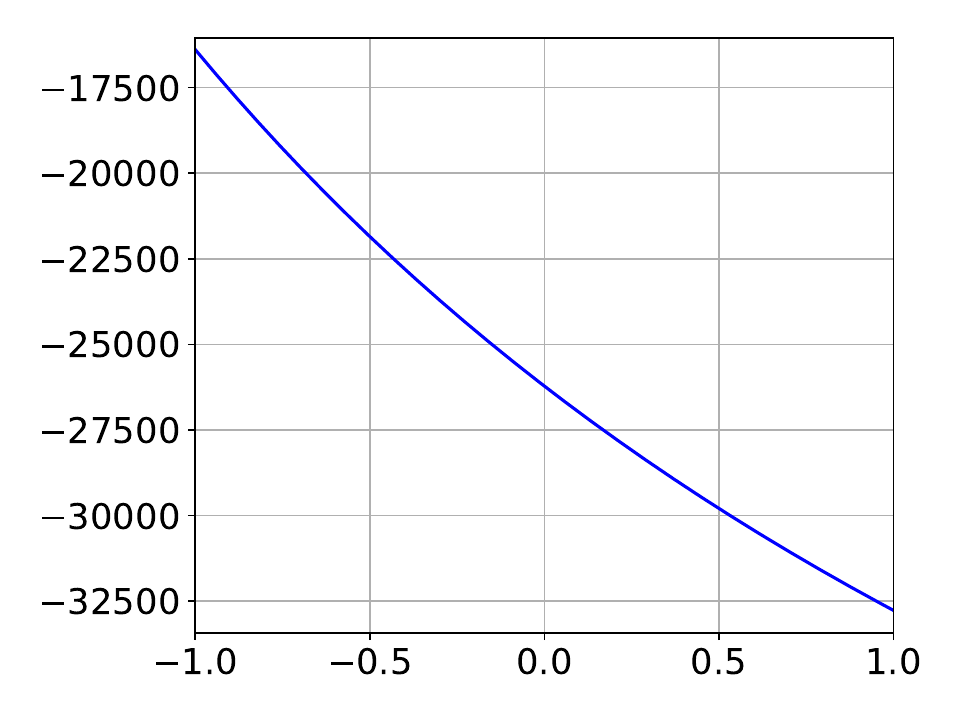}
\hfil
\includegraphics[width=.45\textwidth]{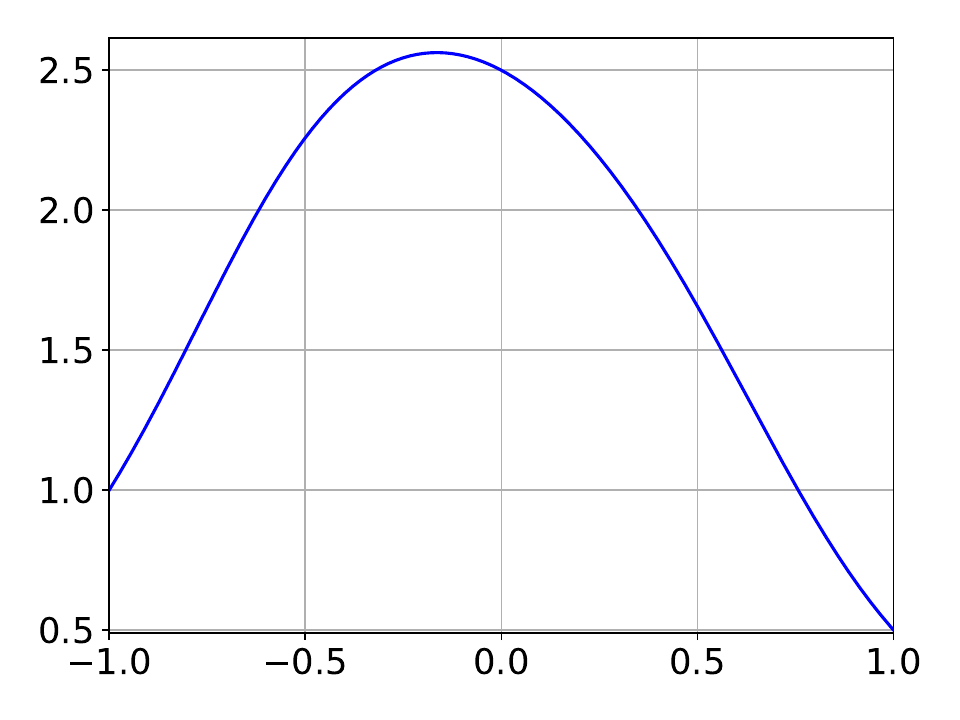}
\hfil
\caption{The derivatives of the two slowly-varying phase functions produced
by applying the global Levin method to the scalar equation (\ref{experiments7:scalar}) obtained
from the system (\ref{experiments7:system}) when $\omega=2^{16}$.  
Each column corresponds to one of the phase functions,
with the real part appearing in the first row and the imaginary part in the second.}
\end{subfigure}
\hfill
\begin{subfigure}{.45\textwidth}
\hfil
\includegraphics[width=.45\textwidth]{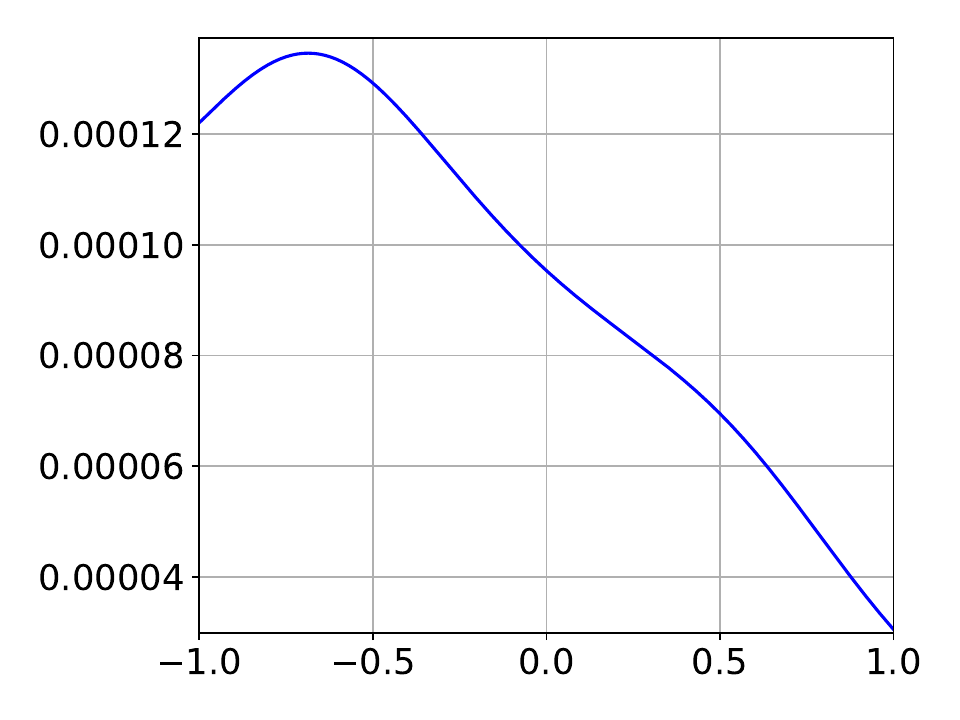}
\hfil
\includegraphics[width=.45\textwidth]{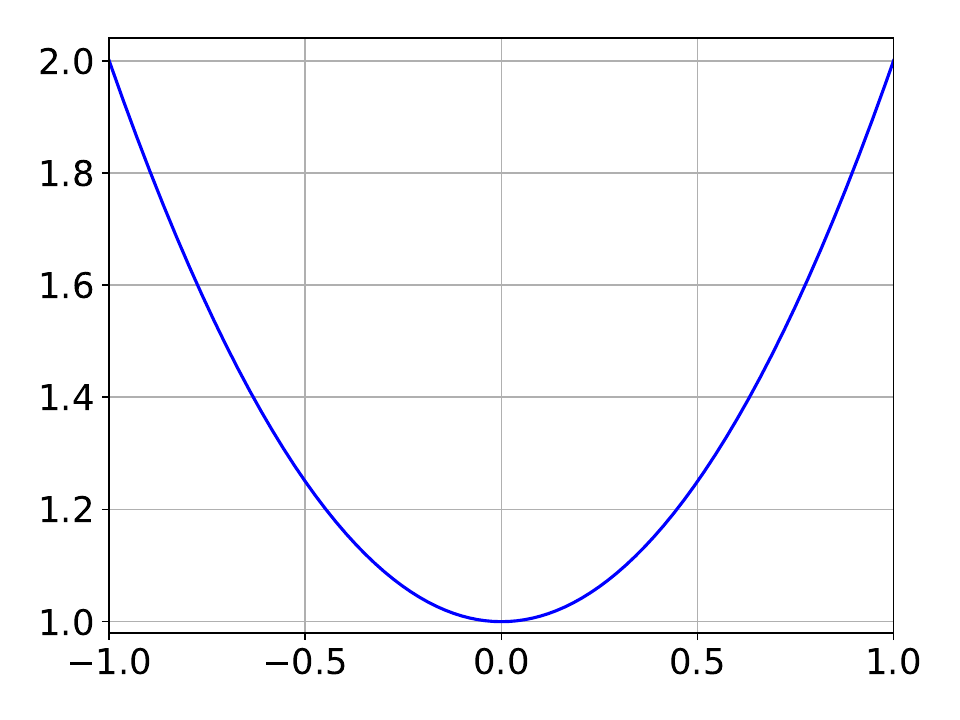}
\hfil

\hfil
\includegraphics[width=.45\textwidth]{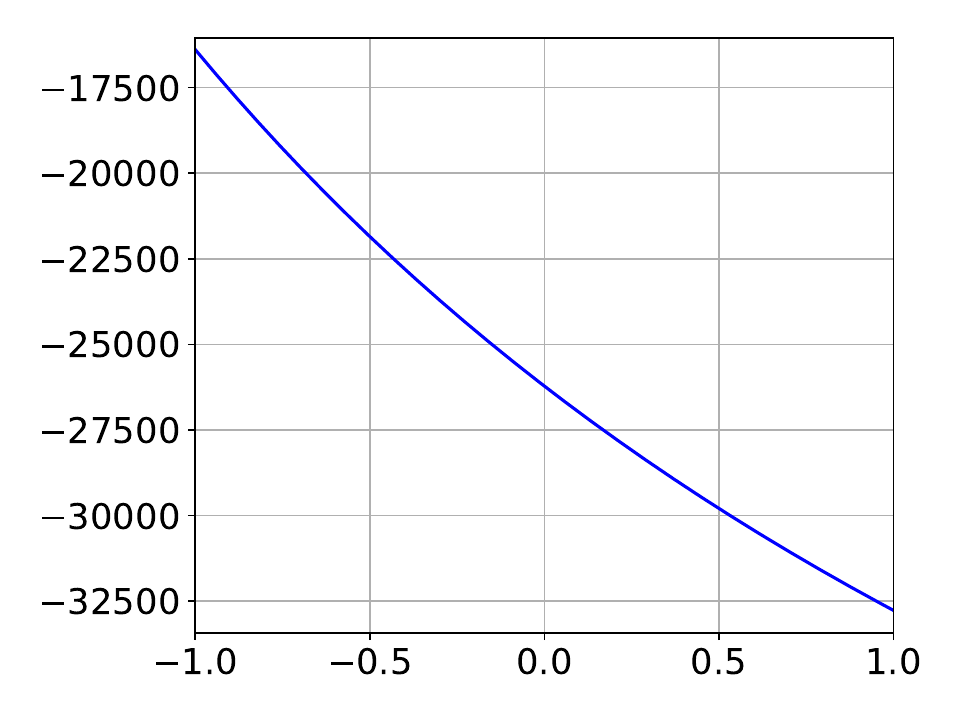}
\hfil
\includegraphics[width=.45\textwidth]{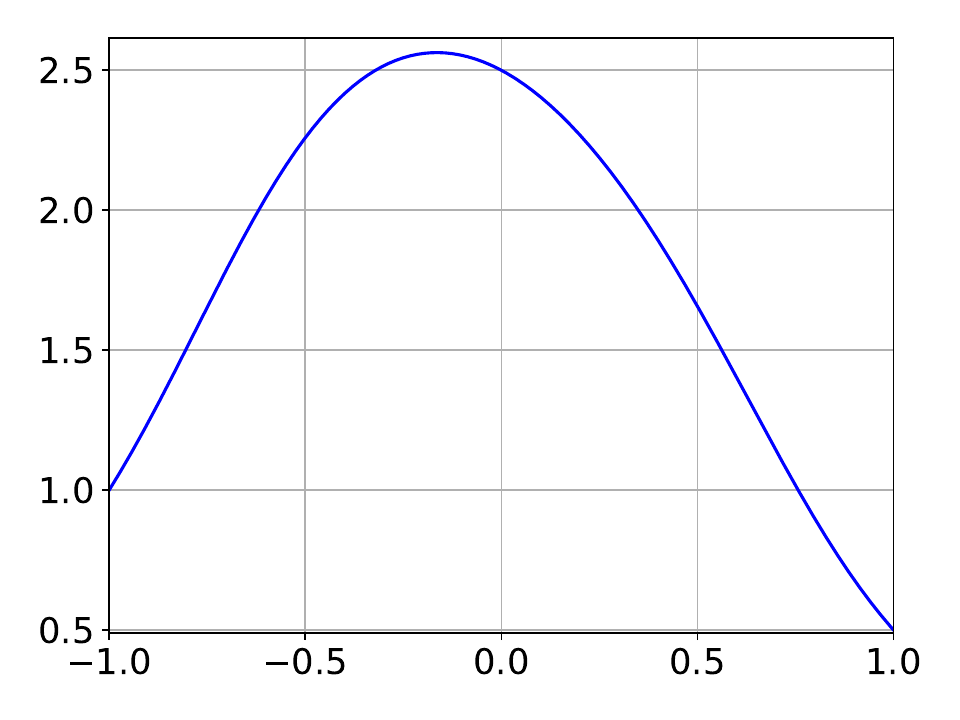}
\hfil
\caption{
The eigenvalues $\lambda_1(t), \lambda_2(t)$ of the 
coefficient  matrix of the system (\ref{experiments7:system})
when $\omega=2^{16}$.  
Each column corresponds to one of the eigenvalues,
with the real part appearing in the first row and the imaginary part in the second.\\
}
\end{subfigure}
\hfil
\caption{Plots of some of the phase functions and eigenvalues
which arose in the course of the  experiments of Subsection~\ref{section:experiments:7}.}
\label{experiments7:figure2}
\end{figure}

For each $\omega=2^8,2^9,\ldots,2^{20}$, we followed the following procedure.
First, we used  the local Levin method
to construct slowly-varying phase functions $\psi_1,\psi_2$ such that
\begin{equation}
z_1(t) = \exp(\psi_1(t)) \ \ \ \mbox{and} \ \ z_2(t) = \exp(\psi_2(t))
\end{equation}
form a basis in the space of solutions of (\ref{experiments7:scalar}) over the interval $[-1,1]$.  
Next, we let formed the two corresponding solutions 
\begin{equation}
\mathbf{y}_1(t) = 
\Phi^{-1}(t)
\left(\begin{array}{c}
z_1(t) \\
z_1'(t) 
\end{array}
\right)
\ \ \ \mbox{and} \ \ \ 
\mathbf{y}_2(t) = 
\Phi^{-1}(t)
\left(\begin{array}{c}
z_2(t) \\
z_2'(t) 
\end{array}
\right)
\end{equation}
of (\ref{experiments7:system}) and found constants $c_1$ and $c_2$ such that
\begin{equation}
c_1 \mathbf{y}_1(0)   + c_2 \mathbf{y}_2(0)   =
\left(
\begin{array}{c}
1\\
1
\end{array}
\right)
\label{experiments7:solution}
\end{equation}
by solving the obvious  system of linear algebraic equations.
Finally, we constructed  a  reference solution for the problem using the  standard
solver of Appendix~\ref{section:appendix}
and compared its value with that of $c_1\mathbf{y}_1(t) + c_2\mathbf{y}_2(t)$
at $10,000$ equispaced points on the interval $[-1,1]$.
Figure~\ref{experiments7:figure1} gives the results.
Figure~\ref{experiments7:figure2}(a)  contains plots of the slowly-varying phase
functions for (\ref{experiments7:scalar}) which were constructed when $\omega=2^{16}$.
Figure~\ref{experiments7:figure2}(b) contains plots of the eigenvalues
of the matrix (\ref{experiments7:system}) when $\omega=2^{16}$.




\end{subsection}

\end{section}

\begin{section}{Conclusions}
\label{section:conclusion}

We have introduced two related approaches  for solving initial and boundary
value problems for a large class of scalar ordinary differential equations. 
They are both based on the principle which underpins the classical Levin method for
calculating oscillatory integrals, namely, that inhomogeneous linear ordinary 
differential equations with slowly-varying coefficients and a slowly-varying
right-hand side admit slowly-varying solutions regardless of the magnitudes
of the coefficients.
 Using Newton's method, we were able to apply
this principle to the nonlinear Riccati equation in order to rapidly
compute slowly-varying phase functions for scalar ordinary differential
equations.

Both of our methods for computing phase functions achieve high-accuracy
and high order convergence, and run in time independent of the magnitude of the equation's coefficients.  
One of the approaches, the global Levin method, fails when the coefficient
matrix has eigenvalues of small magnitude.  The other approach, the local Levin method, is slightly less
accurate and slightly slower than the global method in most cases in which both apply, but it does not suffer from this defect.  

As is clear from our plots, the solutions of the Riccati equation possess various symmetries.
We have made no attempt to exploit them to accelerate our algorithm, or to make a number
of other obvious speed improvements.  Indeed, our implementations are reference codes are meant to 
demonstrate the  basic properties of our algorithms and are not designed to achieve the greatest possible speed.


We encountered difficulties in testing the accuracy of our algorithms.  
For the most part, we did so by using the obtained phase functions to solve an initial or boundary
value problem for a differential equation  and comparing the result  to that 
obtained by a standard solver.  This approach has serious limitations, however, due to the poor 
performance exhibited by standard methods  when applied to most of the problems discussed here.
The problems were so severe, in fact, that we resorted to also measuring the accuracy of our phase 
functions by comparison with results obtained by
running our algorithms using extended precision arithmetic.

There are many obvious applications of this work to evaluation of special functions, calculation
of special function transforms and the solution of physical problems which should be explored.
Moreover,  because essentially any system of $n$ ordinary differential equations in $n$ unknowns
can be transformed into an $n^{th}$ order scalar equation, it should be possible
to use the algorithm of this paper to  solve a large class of systems of differential equations in time
time independent of the magnitudes of the eigenvalues of their coefficient matrices.
We have given one such example in this paper, but further development of this approach is required
in order to obtain a solver which can be applied to a large class of equations.

\end{section}

\begin{section}{Acknowledgments}
JB  was supported in part by NSERC Discovery grant  RGPIN-2021-02613.
\end{section}

\begin{section}{Data availability statement}
The datasets generated during and/or analysed during the current study are available from the corresponding author on reasonable request.
\end{section}

\bibliographystyle{acm}
\bibliography{scalar.bib}

\appendix
\begin{section}{An adaptive spectral solver for ordinary differential equations}
\label{section:appendix}

In this appendix, we detail a standard  adaptive  spectral method for solving ordinary
differential equations.
It is used by the local Levin method, and to calculate reference
solutions in our numerical experiments. We describe its operation in the 
case of the initial value problem
\begin{equation}
\left\{
\begin{aligned}
\bm{y}'(t) &= F(t,\bm{y}(t)), \ \ \ a < t < b,\\
\bm{y}(a) &= \bm{v}
\end{aligned}
\right.
\label{algorithm:system}
\end{equation}
where $F:\mathbb{R}^{n+1} \to \mathbb{C}^n$ is smooth and $\bm{v} \in \mathbb{C}^n$.
However, the solver can be easily modified to produce a solution with a specified value
at any point $\eta$ in $[a,b]$.

The solver takes as input a positive integer $k$, a tolerance parameter $\epsilon$, an interval $(a,b)$, 
the vector $\bm{v}$ and a 
 subroutine for evaluating the function $F$.  It outputs $n$ piecewise $(k-1)^{st}$ order Chebyshev expansions,
one for each of the components $y_i(t)$ of the solution $\bm{y}$ of (\ref{algorithm:system}).

The solver maintains two lists of subintervals of $(a,b)$: one consisting of what we term ``accepted subintervals''
and the other of subintervals which have yet to be processed.  A subinterval is accepted if the solution
is deemed to be adequately represented by a $(k-1)^{st}$ order Chebyhev expansion on that subinterval.
Initially, the list of accepted subintervals is empty and the list of 
subintervals to process contains the single interval $(a,b)$.
It then operates as follows until the list of subintervals to process is empty:
\begin{enumerate}

\item
Find, in the list of subinterval to process, the interval $(c,d)$ such that
$c$ is as small as possible and remove this subinterval from the list.

\item
Solve the initial value problem
\begin{equation}
\left\{
\begin{aligned}
\bm{u}'(t) &= F(t,\bm{u}(t)), \ \ \ c< t < d,\\
\bm{u}(c) &= \bm{w}
\end{aligned}
\right.
\label{algorithm:ivp2}
\end{equation}
If $(c,d) = (a,b)$, then we take $\bm{w}=\bm{v}$.  Otherwise,
the value of the solution at the point $c$ has already been approximated, and we use that estimate
for $\bm{w}$ in (\ref{algorithm:ivp2}).

If the problem is linear, a straightforward Chebyshev integral equation method (see, for instance, \cite{GreengardSolver})
is used to solve (\ref{algorithm:ivp2}).  Otherwise, 
the trapezoidal method is first used to produce an initial
approximation $\bm{y_0}$ of the solution and then Newton's method is applied to refine it.
The linearized problems are solved using a Chebyshev integral equation method.

In any event, the result is a set of $(k-1)^{st}$ order Chebyshev expansions 
\begin{equation}
u_i(t)  \approx \sum_{j=0}^{k-1} \lambda_{ij}\ T_j\left(\frac{2}{d-c} t + \frac{c+d}{c-d}\right),\ \ \ i=1,\ldots,n,
\label{algorithm:exps}
\end{equation}
approximating  the components $u_1,\ldots,u_n$ of the solution of (\ref{algorithm:ivp2}).

\item
Compute the quantities
\begin{equation}
\frac{\sqrt{\sum_{j=\lfloor k/2 \rfloor+1}^{k-1} \left|\lambda_{ij}\right|^2}}{\sqrt{\sum_{j=0}^{k-1} \left|\lambda_{ij}\right|^2}}, \ \ \ i=1,\ldots,n,
\end{equation}
where the $\lambda_{ij}$ are the coefficients in the expansions (\ref{algorithm:exps}).
If any of the resulting values is  larger than $\epsilon$,
then we split the subinterval into two halves $\left(c,\frac{c+d}{2}\right)$ and 
$\left(\frac{c+d}{2},d\right)$ and place them on the list of subintervals to process.  Otherwise, we place the subinterval
$(c,d)$ on the list of accepted subintervals.

\end{enumerate} 

At the conclusion of this procedure,  we have $(k-1)^{st}$ order piecewise Chebyshev expansions
for each component of the solution, with the list of accepted subintervals determining the
partition for each expansion.

\end{section}

\end{document}